\documentclass[10pt,reqno,oneside]{amsart}
 \usepackage{verbatim,amsmath,amssymb,cite,xspace}
\usepackage{color,cite,graphicx}

\theoremstyle{plain}
\newtheorem{theorem}{Theorem}[section]
\newtheorem{claim}[theorem]{Claim}
\newtheorem{lemma}[theorem]{Lemma}
\newtheorem{prop}[theorem]{Proposition} 
\newtheorem{corollary}[theorem]{Corollary}
\theoremstyle{definition}

\theoremstyle{remark}
\newtheorem{remark}[theorem]{Remark}

\usepackage{color}

\usepackage{amssymb}
\usepackage{amssymb, latexsym}
\usepackage{amsmath}
\usepackage{amscd}
\usepackage{enumerate}
\usepackage{amsfonts}
\usepackage{graphicx}
\usepackage{epsfig}

\newcommand{\HL}{\dot{H}^1\times L^2}
\newcommand{\OR}{\overrightarrow}
\newcommand{\E}{\mathcal{E}}

\newcommand{\dual}{2^{\ast}}
\newcommand{\Snorm}{L_t^{\frac{d+2}{d-2}}L_x^{2\frac{d+2}{d-2}}}

\newcommand{\spartial}{ {\slash\!\!\! \partial} }
\newcommand{\eps}{\varepsilon}

\numberwithin{equation}{section} 

\DeclareMathOperator{\supp}{supp}
\DeclareMathOperator{\opdiv}{div}

\title[]{Soliton resolution along a sequence of times for the focusing energy critical wave equation}
\author{Thomas Duyckaerts$^1$}
\thanks{$^1$LAGA, Universit\'{e} Paris 13, Sorbonne Paris Cit\'{e}, UMR 7539. Partially supported by ERC Grant Blowdisol 291214, French ANR Grant SchEq ANR-12-JS01-0005-01, and NSF Grant DMS-1440140 while in residence at the MSRI in Berkeley, during the Fall 2015.}

\author{Hao Jia$^2$}
\thanks{$^2$University of Chicago.}

\author{Carlos Kenig$^3$}
\thanks{$^3$University of Chicago. Partially supported by NSF Grants DMS-1265429 and DMS-1463746.}

\author{Frank Merle$^4$}
\thanks{$^4$Cergy-Pontoise (UMR 8088), IHES. Supported by ERC Grant Blowdisol 291214.}%

\begin{document}

\begin{abstract}In this paper, we prove that any solution of the energy-critical wave equation in space dimensions $3$, $4$ or $5$, which is bounded in the energy space decouples asymptotically, for a sequence of times going to its maximal time of existence, as a finite sum of modulated solitons and a dispersive term. This is an important step towards the full soliton resolution in the nonradial case and without any size restrictions. The proof uses a Morawetz estimate very similar to the one known for energy-critical wave maps, a virial type identity and a new channels of energy argument based on a lower bound of the exterior energy for well-prepared initial data.
\end{abstract}

\maketitle


\section{Introduction}
Consider the energy critical wave equation:
\begin{equation}\label{eq:main}
\partial_{tt}u-\Delta u=|u|^{\dual-2}u,\,\,{\rm in}\,\, R^d\times [0,T_+)
\end{equation}
with initial data $(u_0,u_1)\in \HL(R^d)$, $3\leq d\leq 5$. We use the notation $\dual=\frac{2d}{d-2}$.\\

\hspace{.6cm}The equation is well-posed on $\HL(R^d)$: for each $(u_0,u_1)\in\HL(R^d)$, there exists a unique maximal solution $\OR{u}(t)=(u,\partial_tu)$ satisfying $\OR{u}\in C([0,T_+),\HL)$ and $u\in \Snorm(R^d\times [0,T))$ for any $T<T_+$. Moreover, the energy 
\begin{equation}
\E(\OR{u}(t)):=\int_{R^d}\,\left[\frac{|\partial_tu|^2}{2}+\frac{|\nabla u|^2}{2}-\frac{|u|^{\dual}}{\dual}\right](x,t)\,dx
\end{equation}
is independent of time.\\ 

\hspace{.6cm}Equation (\ref{eq:main}) is invariant under space-time translations and the scaling
\begin{equation*}
u(x,t)\to u_{\lambda}(x,t):=\lambda^{\frac{d}{2}-1}u(\lambda x,\lambda t),
\end{equation*}
 for each $\lambda>0$. 
 The energy is preserved by this scaling.\\

\hspace{.6cm}To put our result in perspective, we will start by a discussion of the extensive background for our work.
As many other nonlinear dispersive equations, equation (\ref{eq:main}) admits solitary waves (or solitons) that are by definition well-localized solutions that travel at a fixed speed. These equations are believed to satisfy the soliton resolution conjecture, namely that any global solution decomposes asymptotically as a sum of decoupled solitons, a dispersive part and a term going to zero. For finite time blow-up solutions, a result in the same spirit is expected, depending on the nature of the blow-up. The soliton resolution conjecture arose in the 70s and 80s from the integrability theory for Korteweg de Vries equations (see \cite{Miura76}, \cite{Segur73}, \cite{Schuur86BO}) and various numerical experiments \cite{IvancevicIvancevic10BO}.\\

\hspace{.6cm}This conjecture is considered quite challenging, even in the completely integrable case. Until recently, to our knowledge, it was only proved for the KdV (see \cite{EcSc83}, \cite{Eckhaus86}) and mKdV (see \cite{Schuur86BO}) equations, which are both completely integrable, by the method of inverse scattering  (see also \cite{ZaSh71}, \cite{SegurAblowitz76}, \cite{Segur76}, \cite{Novoksenov80} for cubic NLS in one space dimension). More precisely, for KdV, a solution with smooth initial data decaying sufficiently at infinity decomposes, for positive $x$, as a finite sum of solitary waves and a term going to zero at infinity. The dispersive component (localized in $\{x<0\}$) is not completely described, and the conjecture is not solved yet for solutions that do not decay fast at infinity. For these solutions, the dynamics is more complicated, as infinite sums of solitons may appear, with time-dependent space shifts that do not converge as $t$ goes to $+\infty$, see \cite[Theorem 2]{MaMe05}.\\

\hspace{.6cm}Very few complete results exist for equations that are not completely integrable. 
However, weaker theorems are available for many dispersive equations, including wave maps and dispersive equations with a power-like linearity with a sufficiently large exponent. For these equations, when solitons exist, the ground state plays an important role as a threshold for the dynamics, and solutions smaller than the ground states, in a suitable sense, are global and scatter to a linear solution (see \cite{KMacta} for equation \eqref{eq:main}, and also \cite{KeMe06} for the energy critical NLS, \cite{Dodson12,Dodson1D,Dodson2D} for the mass-critical NLS, \cite{HoRo08,DuHoRo08,IbMaNa11} for mass-supercritical, energy-subcritical NLS and Klein-Gordon equations etc...). For many of these equations, the dynamics below the ground state can be completed by a complete description of the dynamics in a neighborhood of the ground-state (see \cite{KNS} for equation \eqref{eq:main} and also \cite{NaSc11bo} and references therein for Klein-Gordon and NLS and \cite{MeRa05a} for blow-up solutions of $L^2$-
critical NLS). 
All these works imply local versions of the soliton resolution, where at most one soliton (the ground state) appears.\\

\hspace{.6cm}In \cite{Tao07DPDE,Tao08DPDE}, Tao has proved the existence of an attractor which is compact (modulo the space translations) for $L^2$-supercritical, energy-subcritical NLS in high dimensions, up to a dispersive term. This reduces the proof of a weak variant of the soliton resolution to the proof of a rigidity theorem, namely that any solution with the compactness property (i.e. that has a compact trajectory up to the symmetries of the equation), is a solitary wave. This is a very difficult problem which has been solved so far only in the two regimes cited above (below and close to the ground state), except for KdV, using the complete integrability \cite{LaurentMartel04}, and for equation \eqref{eq:main} where the rigidity theorem was proved by the first, third and fourth authors without any smallness assumption but under an additional nondegeneracy assumption \cite{DKM3}. We also refer to \cite{DKMprofile} where the importance of this type of solution for general dispersive equations is 
highlighted.\\

\hspace{.6cm}The soliton resolution is believed to hold unconditionally for energy-critical wave maps. A first result in this direction is that any solution that blows up in finite time converges locally in space, up to the transformations of the equation, for a sequence of times, to a solitary wave (see \cite{Chri,Struwe2} for the equivariant case, \cite{Tataru4} for the general case). In the  equivariant case, using techniques developed by the first, third and fourth authors for the energy-critical wave equation (\ref{eq:main}), one can prove stronger statements, namely that the soliton resolution holds with a condition on the energy ruling out a multi-soliton configuration \cite{Cotemap1, Cotemap2}, and that it holds for a sequence of times without this condition (see \cite{Cotesoliton} and \cite{JiaKenig}). Recently the method developed in the present article was used by the authors to prove 
the resolution into solitons without any symmetry assumption for blow-up solutions close to the ground-state soliton when the target is $S^2$: see \cite{DuJiKeMe16}, and also \cite{Grinis16} for partial results in the large data case. This showcases the wide scope of applicability of the techniques developed in this paper. For the parabolic analog of wave maps, these questions have been studied earlier (see e.g \cite{Struwe85,Qing95,QingTian97,Topping97}). In full generality, in the parabolic case, the soliton resolution is only known for a sequence of times. An example of Topping \cite{Topping97} shows that for a general target manifold, two sequences of times might lead to different decompositions. In the case of $S^2$ target, it is conjectured that the analog of the soliton resolution conjecture holds, but this has not been completely proved yet.


\medskip

\hspace{.6cm}Going back to equation (\ref{eq:main}), it is known that if $\|(u_0,u_1)\|_{\HL}$ is sufficiently small, then $T_+=\infty$ and the solution scatters, i.e., $u\in \Snorm(R^d\times [0,\infty))$. It is also well known that in general finite energy solutions to equation (\ref{eq:main}) may blow up in finite time, i.e., $T_+<\infty$. Indeed, using the finite speed of propagation for equation (\ref{eq:main}) to localize ODE type blow up solutions, one can easily construct finite time blow up solutions $\OR{u}$ with $\|\OR{u}(t)\|_{\HL}\to\infty$, as $t\to T_+$. These solutions are called \emph{type I} blow up solutions. It is expected that any type I solution decomposes  
into a finite sum of explicit nonlinear objects related to this ODE blow-up, as in the 1-d case (see \cite{MeZa12} and references therein),  however very little is known in the energy-critical case (see \cite{Donninger15P} for a local study).\\

\hspace{.6cm}To rule out the ODE type behavior, one can focus on solutions that are bounded in the energy space, i.e. such that
\begin{equation}
\label{bounded}
\sup_{t\in [0,T_+)}\,\|\OR{u}(t)\|_{\HL}<\infty.
\end{equation}
The dynamics of these solutions is very rich. As we mentioned earlier, small solutions are global and scatter. 
Equation (\ref{eq:main}) admits also various types of finite energy steady states $Q\in\dot{H}^1$, i.e.
\begin{equation}
\label{Ell}
-\Delta Q=|Q|^{\dual-2}Q\,,\quad {\rm in}\,\,R^d,
\end{equation}
(see \cite{Ding86}, \cite{dPMPP11}, \cite{dPMPP13}).
Among them, a distinguished role is played by the {\it ground state}
$$W(x)=\frac{1}{\left( 1+\frac{|x|^2}{N(N-2)} \right)^{\frac{N-2}{2}}}, $$
which is the unique (up to scaling symmetry and sign change) radial $\dot{H}^1$ steady state. $W$ can be characterized as the minimizer of the Sobolev embedding $\dot{H}^1\subset L^{\dual}(R^d)$, and as the non-zero solution of \eqref{Ell} with least energy see \cite{Talenti}.\\

\hspace{.6cm}Recall that equation (\ref{eq:main}) is invariant by the Lorentz transformation:
\begin{equation}
u(x,t)\,\to u_{\ell}(x,t):=u\left(x-\frac{x\cdot \ell}{|\ell |^2}\ell +\frac{\frac{x\cdot\ell}{|\ell |}\frac{\ell}{|\ell |}-\ell t}{\sqrt{1-|\ell |^2}},\,\frac{t-x\cdot\ell}{\sqrt{1-|\ell|^2}}\right),
\end{equation}
for each $\ell\in R^d$, with $|\ell |<1$: if $u$ is a classical solution to equation (\ref{eq:main}), then $u_{\ell}$ is also a solution where it is defined. 
Taking the Lorentz transform of a steady state, we obtain traveling wave solutions for $\ell\in R^d$ with $|\ell |<1$:
\begin{equation}
Q_{\ell}(x,t):=Q\left(x-\frac{x\cdot \ell}{|\ell |^2}\ell +\frac{\frac{x\cdot\ell}{|\ell |}\frac{\ell}{|\ell |}-\ell t}{\sqrt{1-|\ell |^2}}\right).
\end{equation}
Note that $Q_{\ell}(x,t)=Q_{\ell}(x-t\ell,0)$, so that $Q_{\ell}$ travels in the direction of $\ell$ with speed $|\ell |$. 
There are more complicated global bounded solutions, that are asymptotically the sum of a radiation term and one or several decoupled traveling waves (see \cite{KrSc07}, \cite{DoKr13}, \cite{MaMe16}, and \cite{JendrejTwo}). \\

\hspace{.6cm}There also exist solutions blowing up in finite time that are bounded in the energy space. These solutions are called \emph{type II blow up solutions}.
In \cite{KSTwave}, \cite{HillRap}, \cite{KrSc14} and \cite{Jendrej}, type II blow up solutions of the form of a rescaled ground state plus a small dispersive term  were constructed. More precisely the solution is given by 
\begin{equation*}
u(x,t)=\lambda(t)^{-\frac{d}{2}+1}W\left(\frac{x}{\lambda(t)}\right)+\epsilon(x,t),
\end{equation*}
where $\frac{\lambda(t)}{T_+-t}\to 0+$ as $t\to T_+$, and $\OR{\epsilon}(t)=(\epsilon,\,\partial_t\epsilon)$ is small in the energy space. It is expected that multi-soliton concentration is also possible for type II blow up solutions, and it is an open problem to construct such a solution. Similar blow up solutions have been constructed for the energy critical equivariant wave maps and the radial energy critical Yang Mills equation, see \cite{KSTwavemap},  \cite{KSTmills}, \cite{RodSter} and \cite{RapRod}.\\

\hspace{.6cm} The soliton resolution conjecture for 
equation (\ref{eq:main}) predicts that any bounded solution should asymptotically decouple into a finite sum of modulated solitons, a regular part in the finite time blow up case or a free radiation in the global case, plus a residue term that vanishes asymptotically in the energy space as time approaches the maximal time of existence. We note that all the solutions mentioned above are in accordance with this conjecture.\\ 

\hspace{.6cm}In the radial setting,
$W$ is the unique steady state, and thus the only soliton, since the Lorentz transformation destroys radiality.  The soliton resolution conjecture was proved in \cite{DKM} by the first, third and fourth authors, for $d=3$. For other dimensions (still in the radial case), soliton resolution is only known along a sequence of times, see \cite{Kenig4dWave}, \cite{Casey} and \cite{JiaKenig}.\\

\hspace{.6cm}The proof of \cite{DKM} uses techniques that are very specific to radial solutions in $3$ space dimensions. The core of the proof is a characterization of $W$ as the only solution that does not satisfy a particular bound from below of the outer energy (\emph{channel of energy inequality}). This property is proved using an analogous bound for the linear wave equation, valid on a subspace of $\dot{H}^1\times L^2$ which is of codimension $1$. In the nonradial setting and in even dimension \cite{CoKeSc14,KeLaLiSc15}, there is an infinite dimensional subspace of $\dot{H}^1\times L^2$ for which this linear property does not hold, and we do not expect that this method of proof can be adapted.\\

\hspace{.6cm}In the nonradial case, the soliton resolution was proved for type II blow-up solutions in \cite{DKMsmallnonradial} for $d=3,\,5$, under an extra smallness condition. For general large data, it was proved in \cite{DKMprofile} that along a sequence of times, the solution converges locally, after an appropriate rescaling, to a modulated soliton. \\

\hspace{.6cm} In this paper, we prove the soliton resolution along a sequence of times for general bounded solutions $\OR{u}(t)$ of (\ref{eq:main}).
\begin{theorem}\label{th:main}
Let $\OR{u}\in C([0,T_+),\HL(R^d))$, with $u\in \Snorm(R^d\times [0,T))$ for any $T<T_+$, be a solution to equation (\ref{eq:main}) that satisfies (\ref{bounded}). Here $T_+$ denotes the maximal existence time of $u$.  \\

{\bf Case I: $T_+<\infty$}.  Define the singular set 
\begin{equation}
\mathcal{S}:=\left\{x_{\ast}\in R^d:\,\|u\|_{\Snorm\left(B_{\epsilon}(x_{\ast})\times [T_+-\epsilon,\,T_+)\right)}=\infty,\quad{\rm for\,\,any\,\,}\epsilon>0\right\}.
\end{equation}
Then $\mathcal{S}$ is a finite set. Let $x_{\ast}\in\mathcal{S}$ be a singular point. Then there exist an integer $J_{\ast}\ge 1$, $r_{\ast}>0$, $\OR{v}_0\in \HL$, a time sequence $t_n\uparrow T_+$,  scales $\lambda_n^j$ with $0<\lambda_n^j\ll T_+-t_n$, positions $c_n^j\in R^d$ satisfying $c_n^j\in B_{\beta (T_+-t_n)}(x_{\ast})$ for some $\beta\in(0,1)$ with $\ell_j=\lim\limits_{n\to\infty}\frac{c_n^j-x_{\ast}}{T_+-t_n}$ well defined, and traveling waves $Q_{\ell_j}^j$, for $1\leq j\leq J_{\ast}$, such that inside the ball $B_{r_{\ast}}(x_{\ast})$ we have
\begin{eqnarray}
\OR{u}(t_n)&=&\OR{v}_0+\sum_{j=1}^{J_{\ast}}\,\left((\lambda_n^j)^{-\frac{d}{2}+1}\, Q_{\ell_j}^j\left(\frac{x-c_n^j}{\lambda_n^j},\,0\right),\,(\lambda_n^j)^{-\frac{d}{2}}\, \partial_tQ_{\ell_j}^j\left(\frac{x-c_n^j}{\lambda_n^j},\,0\right)\right)\nonumber\\
\nonumber\\
 &&\quad\quad \quad\quad\vspace{1.5cm}+\,\,o_{\HL}(1), \quad\,\,{\rm as}\,\,n\to\infty.\label{eq:maindecompositionhaha}
\end{eqnarray}
In addition, the parameters $\lambda_n^j,\,c_n^j$ satisfy the pseudo-orthogonality condition
\begin{equation}
\label{ortho}
\frac{\lambda_n^j}{\lambda_n^{j'}}+\frac{\lambda_n^{j'}}{\lambda_n^j}+\frac{\left|c^j_n-c^{j'}_n\right|}{\lambda^j_n}\to\infty,
\end{equation}
as $n\to\infty$, for each $1\leq j\neq j'\leq J_{\ast}$.\\

{\bf Case II: $T_+=\infty$.}
There exist a finite energy solution $u^L$ to the linear wave equation
$$ \partial_{tt}u^L-\Delta u^L=0\text{ in }R^d\times R,$$
an integer $J_{\ast}\ge 0$ \footnote{If $J_{\ast}=0$, then there is no soliton in the decomposition below and the solution scatters.}, a time sequence $t_n\uparrow \infty$,  scales $\lambda_n^j$ with $\lambda_n^j>0$ and $\lim\limits_{n\to\infty}\frac{\lambda^j_n}{t_n}=0$, positions $c_n^j\in R^d$ satisfying $c_n^j\in B_{\beta t_n}(0)$ for some $\beta\in(0,1)$ with $\ell_j=\lim\limits_{n\to\infty}\frac{c_n^j}{t_n}$ well defined, and traveling waves $Q_{\ell_j}^j$, for $1\leq j\leq J_{\ast}$, such that 
\begin{eqnarray}
\OR{u}(t_n)&=&\OR{u}^L(t_n)+\sum_{j=1}^{J_{\ast}}\,\left((\lambda_n^j)^{-\frac{d}{2}+1}\, Q_{\ell_j}^j\left(\frac{x-c_n^j}{\lambda_n^j},\,0\right),\,(\lambda_n^j)^{-\frac{d}{2}}\, \partial_tQ_{\ell_j}^j\left(\frac{x-c_n^j}{\lambda_n^j},\,0\right)\right)\nonumber\\
\nonumber\\
 &&\quad\quad \quad\quad\vspace{1.5cm}+\,\,o_{\HL}(1), \quad\,\,{\rm as}\,\,n\to\infty.\label{eq:maindecompositionhahaGlobal}
\end{eqnarray}
In addition, the parameters $\lambda_n^j,\,c_n^j$ satisfy the pseudo-orthogonality condition (\ref{ortho}).
\end{theorem}

\smallskip

\begin{remark}
The radiation term $u^L$ appearing in Case II was constructed in \cite{DKMsp} and does not depend on the choice of the sequence $t_n$ (see Lemma \ref{lm:radiation} below).
\end{remark}
\begin{remark}
 \label{R:time}
 Following carefully the choice of the sequences of times in the proof of Theorem \ref{th:main}, we can prove that if $t_n^1$ and $t_n^2$ are any two sequences such that
 $$\lim_{n\to+\infty} t_n^1=\lim_{n\to+\infty} \frac{t_n^2}{t_n^1}=+\infty,$$
 then, extracting subsequences, we can find a sequence $t_n$ satisfying the conclusion of Theorem \ref{th:main} and such that 
 $$\forall n, \quad t_n^1\leq t_n\leq t_n^2.$$
\end{remark}

\begin{remark}
 We note that the asymptotic decomposition (\ref{eq:maindecompositionhaha}) in the finite time blow-up case may hold locally only, as there may be more than one blow up points. Thus the residue term vanishes asymptotically only near the blow up point in the  energy space. In this case, the local norm $\dot{H}^1(B_r)$ for $r>0$ is slightly ambiguous. To be precise, we adopt the convention that 
$$\|f\|_{\dot{H}^1(B_r)}^2:=\|\nabla f\|_{L^2(B_r)}^2+\left\|\frac{f}{|x|}\right\|_{L^2(B_r)}^2.$$ 
Let us also mention that in the finite time blow-up case, it follows easily from the proof that one may choose the same sequence $t_n$ for all singular points.
\end{remark}

\medskip

\hspace{.6cm}Below we outline the main ideas of the proof.\\

\hspace{.6cm}Let us firstly consider the finite time blow-up case, $T_+<\infty$. The first new input in the proof of Theorem \ref{th:main} is the observation that the quantity 
\begin{equation}\label{eq:mainpartofenergyflux}
\int_t^{T_+}\int_{|x|=T_+-s}\,\left[\frac{|\nabla u|^2}{2}+\frac{|\partial_tu|^2}{2}-\frac{x}{|x|}\cdot\nabla u\,\partial_tu\right](x,s)\,d\sigma ds
\end{equation}
can be controlled using the energy flux
\begin{equation*}
\int_t^{T_+}\int_{|x|=T_+-s}\,\left[\frac{|\nabla u|^2}{2}+\frac{|\partial_tu|^2}{2}-\frac{x}{|x|}\cdot\nabla u\,\partial_tu-\frac{|u|^{\dual}}{\dual}\right](x,s)\,d\sigma ds
\end{equation*}
and a simple trace theorem,
despite the fact that the energy flux may not be positive in general. The control of (\ref{eq:mainpartofenergyflux}) allows us to use a Morawetz estimate, similar to the one used for the energy critical wave maps, see \cite{Chri,GrillakisEnergy,taoIII,Tataru4}. We prove this estimate by a self-similar change of variables. For conformally-subcritical nonlinearities, the same change of variables yields a Lyapunov functional, the self-similar energy, which is the main tool to study type I blow-up solutions (see \cite{MeZa12} and re\-fe\-ren\-ces therein). In the setting of equation \eqref{eq:main}, the self-similar energy is not well-defined, but one can prove a related approximate monotonicity formula (close to the one used to rule out compact self-similar blow-up in \cite{KMacta}) which implies the desired Morawetz estimate. The estimate on the energy flux is crucial to control the boundary terms in self-similar variables.\\ 

\hspace{.6cm}We next consider a profile decomposition of $\OR{u}(t_n)$ (see \cite{BaGe}) along a time sequence $\{t_n\}_n$ going to $T_+(u)$ as $n$ goes to infinity. If the sequence $(t_n)_n$ is well-chosen, the Morawetz estimate implies that all the profiles that are not translated in time (i.e. that converge up to space translation and scaling) are solitary waves. These profiles are exactly the ones for which the $L^{\dual}$ norm does not goes to $0$ as $n$ goes to infinity. As a result, we obtain a preliminary decomposition similar to the decomposition (\ref{eq:maindecompositionhaha}) but with the residue term vanishing in $L^{\dual}$ only, instead of vanishing in the energy space.\\

\hspace{.6cm}We remark that the Morawetz estimate does not seem sufficient to rule out the profiles (specific to the wave equation in comparison with the elliptic situation) with a nontrivial time-translation, i.e. the profiles that are asymptotically close to an outgoing or an incoming wave.  
To prove that these profiles do not exist, we use a virial identity. In addition to the characterization of all profiles, the virial identity allows us to describe more precisely the residue term, proving that it has certain concentration and vanishing properties. This will play a fundamental role in the last step of the proof.\\

\hspace{.6cm}The idea of using an additional virial identities to eliminate dispersive energy was first introduced in \cite{JiaKenig}. In our case, an interesting new feature is that the virial identity is useful only after we have obtained the preliminary decomposition (integrating between two sequences of times for which this decomposition holds), unlike in the radial case in \cite{JiaKenig} where it is possible to use the virial identity directly. The difference is that in the radial case, we know that there is asymptotically no energy in the self similar region, see \cite{Kenig4dWave, JiaKenig}, while this is not true in general in the nonradial case.\\


\hspace{.6cm}It remains to prove that the residue term goes to $0$ in the energy space. For this we use a new channel of energy inequality for certain well-prepared initial data. These initial data are in some sense outgoing, asymptotically radial, and concentrating energy in a thin annulus.   
The residue term satisfies these properties thanks to the vanishing conditions established in the previous step. Arguing by contradiction, we obtain that any energy coming from the residue term would travel back in an arbitrarily small neighborhood of the wave cone, yielding a nontrivial concentration of energy at a fixed time in the domain of existence of $u$, which is absurd.\\

\hspace{.6cm}This particular channel of energy argument does not depend on dimension, and appears to be among the first successful applications of channel of energy inequalities for even dimensions in the nonradial case.\footnote{A different channel of energy argument in the nonradial case, which is also independent of dimension, was used in the characterization of solutions with the compactness property. See \cite{DKM3}. } We note that this new type of channel of energy argument seems much more robust than the one 
used in \cite{DKM}, which is specific to the radial setting. It was recently adapted to the critical wave maps equation (see \cite{DuJiKeMe16}).\\

\hspace{.6cm}In the global case, besides the steps outlined above in the finite time blow-up case, an additional major step is to extract the radiation term from the solution. This has been done in \cite{DKMsp}. After this step, we can perform the same Morawetz estimates and follow similar arguments as in the finite time blow-up case to prove Theorem \ref{th:main}.\\

\hspace{.6cm}With this method of proof, relying on monotonicity laws giving convergence only after an averaging in time, we cannot hope for more than a decomposition for particular sequences of times (see Remark \ref{R:time} above). The difficulty of obtaining the resolution for all times is illustrated by the heat flow equation, for which the analog of Theorem \ref{th:main} is known, but the soliton resolution for all times does not hold in full generality because of the example of Topping \cite{Topping97} mentioned above.
However we see our result as an important step toward the proof of the resolution  for all times. Indeed, we are now reduced to study the dynamics close to a sum of solitons plus a dispersive term, rather than the general large data dynamics. The main challenge in this study is to prove that the collision of two or more solitons produces dispersion (see \cite{MartelMerle11} for gKdV). Note that in the radial case if $d=3$, this property is a consequence of the rigidity result \cite[Section 2]{DKM} that states that any nonstationary solution produces dispersion.

\begin{remark}
The results in this paper should remain valid for $d\geq 6$. The mo\-difications to the argument that are needed in order to establish this are technical in nature, having to do with the modifications that are needed to establish the local theory of the Cauchy problem, the profile decomposition in Section 2, and the ``approximation theorem" Lemma \ref{lm:nonlinearprofiledecomposition}, due to the limited smoothness of the nonlinearity $F(u)=|u|^{\frac{4}{d-2}}u$.

In space dimension $d=6$, the second derivative of the nonlinearity $|u|u$ is bounded and the Cauchy theory that we use here, based on the space $\Snorm(R^d\times R)=L^2_tL^4_x(R^6\times R)$, works the same way. However this space is an ``endpoint Strichartz space" for which the profile decomposition is not valid anymore. More precisely, the $\Snorm$ norm of the remainder of the profile decomposition (the first term in the left hand side of (\ref{righthandside5})), need not tend to $0$ in general, unlike in the case $d=3,4,5$, as we prove in Appendix \ref{A:residue}. This is a small technical problem, and can be dealt with using other Strichartz norms (such as the $S$ norm defined below) in the profile decomposition.

In dimension $d\geq 7$, the technicalities introduced by the limited smoothness of the nonlinearity become more serious. Furthermore, since $\frac{d+2}{d-2}<2$, the space $\Snorm$ is not a Strichartz space anymore, and the local Cauchy theory must be modified. These problems are dealt with in \cite{BCLPZ}  (see also \cite{Casey}), where an appropriate local Cauchy theory, based on the following norms (introduced in \cite{KMacta}):
$$ \|u\|_{S}=\|u\|_{L^{\frac{2(d+1)}{d-2}}(R^d\times R)},\quad \|u\|_{W}=\||D|^{1/2} u\|_{L^{\frac{2(d+1)}{d-1}}(R^d\times R)}$$
 is developed.
The introduction of these spaces 
allows one to get away from the ``endpoint Strichartz estimate" and thus restores the validity of (\ref{righthandside5}), using the $S$ norm instead of the $\Snorm$ norm. The price to pay for this is that one also needs the $W$ norm, which introduces fractional derivatives. This makes it necessary to use the Leibniz rule for fractional derivatives in many places, and in particular, in ``truncation" arguments outside light-cones, as are used, for instance, in Lemma 2.9. Corresponding difficulties were met successfully in \cite{DKM}, Appendix A, and the same techniques can be used here. Similar comments apply to the results in \cite{DKMsp}. 

We have refrained from carrying out this in detail here, in order not to overburden the exposition with technicalities.
\end{remark}

\hspace{.6cm}This paper is an extension of the arXiv preprint 1510:00075 (see \cite{JiaDispersiveError}) by the second author (which will not be submitted for publication), where Theorem \ref{th:decompbetter} is proved, that is, the case $T_+<\infty$ of Theorem \ref{th:main}, with the error converging to $0$ in the weaker sense of (\ref{eq:refinedintro}). Using the main result in \cite{DKMsp}, this is extended to the case $T_+=\infty$. Moreover, using the new channel of energy arguments in Section 8 (which were inspired by the proof of Claim A.5 in \cite{DKMsp}) the full convergence in energy norm is established here in both the cases $T_+<\infty$ and $T_+=\infty$, in Sections 8 and 9.\\

Our paper is organized as follows.\\

\noindent
\begin{itemize}
 \item In Section 2 we recall basic properties of profile decompositions for wave equations and properties of solutions to the linear wave equation;\\
\item
In Section 3 we derive the Morawetz estimate;\\
\item
In Section 4 we obtain a vanishing condition from the Morawetz estimate;\\
\item
 In Section 5 using the vanishing condition we obtain a preliminary decomposition;\\
\item
 In Section 6 we use the virial identity to derive further vanishing conditions;\\
\item
 In Section 7 we rule out profiles coming from time infinity;\\
\item
 In Section 8 we show that the dispersive term can not contain any non\-tri\-vial amount of energy using the channel of energy argument and finish the proof in the finite blow-up case;\\
\item
 In Section 9 we briefly explain the steps for proving Theorem \ref{th:main} in the global case.
\end{itemize}

\medskip

\noindent

\section{Preliminaries on Profile decompositions}
We briefly recall the profile decompositions, which was firstly introduced to wave equations by Bahouri and G\'erard in $R^3$ (see \cite[Section III]{BaGe}), and was then extended to general dimensions by Bulut (see \cite[Theorem 1.1]{Bulut}).\\

\hspace{.6cm}Let $(u_{0n},u_{1n})\in \HL$ satisfy $\sup\limits_{n}\|(u_{0n},u_{1n})\|_{\HL}\leq M<\infty$. Then passing to a subsequence there exist scales $\lambda^j_n>0$, positions $c^j_n\in R^d$, time translations $t^j_n\in R$, and finite energy solutions $U_j^L$ to the linear wave equation for each $j$, such that we have the following decomposition
\begin{eqnarray}
&&\hspace{2in}(u_{0n},\,u_{1n})=\nonumber\\
&&\quad\quad\quad\sum_{j=1}^J\,\left(\left(\lambda^j_n\right)^{-\frac{d}{2}+1}U_j^L\left(\frac{x-c^j_n}{\lambda^j_n}, \,-\frac{t^j_n}{\lambda^j_n}\right),\,\left(\lambda^j_n\right)^{-\frac{d}{2}}\partial_tU_j^L\left(\frac{x-c^j_n}{\lambda^j_n}, \,-\frac{t^j_n}{\lambda^j_n}\right)\right)\nonumber\\
&&\nonumber\\
&&\hspace{1in}\quad\quad\quad\quad\quad+\left(w^J_{0n},\,w^J_{1n}\right),\label{eq:profiledecomposition}
\end{eqnarray}
where the parameters satisfy
\begin{equation}\label{eq:profileorthogonality1}
t^j_n\equiv 0,\,\,{\rm for\,\,all\,\,}n,\,\,{\rm or}\,\,\lim_{n\to\infty}\frac{t^j_n}{\lambda^j_n}\in\{\pm \infty\},
\end{equation}
and for $j\neq j'$
\begin{equation}\label{eq:profileorthogonality2}
\lim_{n\to\infty}\,\left(\frac{\lambda^j_n}{\lambda_n^{j'}}+\frac{\lambda^{j'}_n}{\lambda^j_n}+\frac{\left|c^j_n-c^{j'}_n\right|}{\lambda^j_n}+\frac{\left|t^j_n-t^{j'}_n\right|}{\lambda^j_n}\right)=\infty.
\end{equation}
In addition, let $\OR{w}^J_n$ be the solution to the linear wave equation with $\OR{w}^J_n(0)=\left(w^J_{0n},\,w^J_{1n}\right)$, then $w^J_n$ vanishes in the sense that
\begin{equation}\label{righthandside5}
\lim_{J\to\infty}\,\limsup_{n\to\infty}\left(\|w^J_n\|_{\Snorm(R^d\times R)}+\|w^J_n\|_{L^{\infty}_tL^{\dual}_x(R^d\times R)}\right)=0.
\end{equation}
 (See (3.48) in \cite{BaGe} and Appendix \ref{A:residue}  for the convergence to $0$ of the $L^{\infty}_tL^{\dual}_x(R^d\times R)$ norm.)

Moreover, if we write for $j\leq J$
\begin{equation}\label{eq:tildewj}
w^J_n(x,t)=(\lambda^j_n)^{-\frac{d}{2}+1}\widetilde{w}^J_{jn}\left(\frac{x-c^j_n}{\lambda^j_n},\,\frac{t-t^j_n}{\lambda^j_n}\right),
\end{equation}
then $\OR{\widetilde{w}}^J_{jn}(\cdot,0)=\left(\widetilde{w}^J_{jn},\,\partial_t\widetilde{w}^J_{jn}\right)(\cdot,0)\rightharpoonup 0$, as $n\to\infty$ for each $j\leq J$. \\

\hspace{.6cm}For later applications, we also need to recall some properties of free radiations (i.e, finite energy solutions to the linear wave equation) and of the profile decomposition. We begin with the following lemma which describes the concentration property of free radiations.
\begin{lemma}\label{lm:concentrationoffreeradiation}
Let $U$ be the solution to the $d+1$ dimensional linear wave equation in $R^d\times R$ with initial data $(U_0,\,U_1)\in\HL(R^d)$. Then
\begin{equation}\label{eq:concentrationforfreeradiation}
\lim_{\lambda\to\infty}\,\sup_{t\in R}\int_{\left||x|-|t|\right|>\lambda}\left[|\partial_tU|^2+|\nabla U|^2\right](x,t)\,dx=0.
\end{equation}
\end{lemma}

\smallskip
\noindent
{\it Proof.}  This is proved in \cite{DKMsp} (see Claim A.5) using a virial-type identity. We give here another proof for the sake of completeness. 

For any $\epsilon>0$, there exists Schwartz class initial data $(u_0,\,u_1)$ such that
\begin{equation}
\|(u_0,\,u_1)-(U_0,\,U_1)\|_{\HL}<\epsilon.
\end{equation}
Let $\OR{u}$ be the solution to the $d+1$ dimensional linear wave equation with initial data $(u_0,\,u_1)$. Then by Theorem 1.1 in \cite{weightedestimate},
\begin{equation}\label{eq:boundforucon}
\sup_{(x,t)\in R^d\times R}(1+|x|+|t|)^{\frac{d-1}{2}}(1+||x|-|t||)^{\frac{d-1}{2}}|\nabla_{x,t}u(x,t)|<\infty.
\end{equation}
Hence for $u$, we can verify directly from (\ref{eq:boundforucon}) that 
\begin{equation}\label{eq:concentrationforu12}
\lim_{\lambda\to\infty}\,\sup_{t\in R}\int_{\left||x|-|t|\right|>\lambda}\left[|\partial_tu|^2+|\nabla u|^2\right](x,t)\,dx=0.
\end{equation}
Note that by energy conservation for linear wave equation, we have for all $t$
\begin{equation*}
\|\OR{u}(t)-\OR{U}(t)\|_{\HL}<\epsilon.
\end{equation*}
Thus,
\begin{equation}
\limsup_{\lambda\to\infty}\,\sup_{t\in R}\int_{\left||x|-|t|\right|>\lambda}\left[|\partial_tU|^2+|\nabla U|^2\right](x,t)\,dx\lesssim \epsilon^2.
\end{equation}
Since $\epsilon>0$ is arbitrary, (\ref{eq:concentrationforfreeradiation}) follows.\\

\hspace{.6cm}We shall also need the following lemma on the ``absolute continuity" of energy distribution for the free radiation on its concentration set $||x|-|t||=O(1)$, adapted to profiles. Denote
\begin{equation}
\OR{U}^L_{jn}(x,t)=\left(\left(\lambda^j_n\right)^{-\frac{d}{2}+1}U_j^L\left(\frac{x-c^j_n}{\lambda^j_n}, \,\frac{t-t^j_n}{\lambda^j_n}\right),\,\left(\lambda^j_n\right)^{-\frac{d}{2}}\partial_tU_j^L\left(\frac{x-c^j_n}{\lambda^j_n}, \,\frac{t-t^j_n}{\lambda^j_n}\right)\right).
\end{equation}
We have
\begin{lemma}\label{lm:absolutecontinuity}
Let $\OR{U}^L_{jn}(0)$ be a linear profile in $\HL$ with $\lim\limits_{n\to\infty}\frac{t^j_n}{\lambda^j_n}\in\{\pm\infty\}$. Suppose that a sequence of sets $E_n\subseteq R^d$ satisfies for all $M>0$ that
\begin{equation}
\label{bound_E}
\Big|\left\{x\in R^d:\,\left||x-c^j_n|-|t^j_n|\right|<M\lambda^j_n\right\}\cap\, E_n\Big|=o\left(\left|t^j_n\right|^{d-1}\lambda^j_n\right),
\end{equation}
as $n\to\infty$. Then we have
\begin{equation}
\lim_{n\to\infty}\,\|\OR{U}^L_{jn}(0)\|_{\HL(E_n)}=0.
\end{equation}
\end{lemma}

\smallskip
\noindent
{\it Proof.} If $\OR{U}^L_{jn}$ has Schwartz class initial data, the lemma follows from the bound (\ref{eq:boundforucon}). The general case follows from approximation.\\

\hspace{.6cm}We shall also need the following orthogonality property for the profiles:
\begin{lemma}\label{eq:profileorthopre}
Suppose that $(u_{0n},\,u_{1n})$ is a bounded sequence in $\HL$ and has the profile decomposition (\ref{eq:profiledecomposition}). Then for each $J$
\begin{equation}\label{eq:almostequalno1}
\|(u_{0n},\,u_{1n})\|^2_{\HL}=\sum_{j=1}^J\|\OR{U}^L_j(0)\|^2_{\HL}+\|(w^J_{0n},w^J_{1n})\|_{\HL}^2+o_n(1),
\end{equation}
as $n\to\infty$, and for all $t$
\begin{eqnarray}
&&\int_{R^d}\,\left[(u_{1n}+\partial_1u_{0n})^2+|\nabla_{x'}u_{0n}|^2\right](x)\,dx\nonumber\\
&&\quad=\sum_{j=1}^J\int_{R^d}\,\left[\left(\partial_tU^L_j+\partial_1U^L_j\right)^2+|\nabla_{x'}U^L_j|^2\right](x,t)\,dx\nonumber\\
&&\quad\quad+\int_{R^d}\left[\left(w^J_n+\partial_1w^J_n\right)^2+|\nabla_{x'} w^J_n|^2\right](x,t)\,dx+o_n(1),\label{eq:almostequalno2}
\end{eqnarray}
as $n\to\infty$. In the above, $x=(x_1,x')$.
\end{lemma}

\smallskip
\noindent
{\it Proof.} The proof follows easily from the fact that for the linear wave $\OR{U}$, the quantities $\|\OR{U}\|^2_{\HL}$ and $\int\,(\partial_tU+\partial_1U)^2+|\nabla_{x'}U|^2\,dx$ are conserved, and the orthogonality condition for the profiles. Indeed, the first quantity is the energy, and the second quantity is a sum of the energy and the momentum in the $x_1$ direction. For the sake of completeness, we briefly outline the proof of (\ref{eq:almostequalno2}). It suffices to prove that for each $1\leq j\neq j'\leq J$\\
\begin{eqnarray}
&&\lim_{n\to\infty}\int_{R^d}\left[\left(\partial_t\,U^L_{jn}+\partial_1\,U^L_{jn}\right)\left(\partial_t\,U^L_{j'n}+\partial_1\,U^L_{j'n}\right)
+\nabla_{x'}\,U^L_{jn}\cdot\nabla_{x'}\,U^L_{j'n}\right](x,t)dx\nonumber\\
&&\hspace{1in}\quad\quad=0;\label{eq:orthogo11}\\
\nonumber\\
&&\lim_{n\to\infty}\int_{R^d}\left[\left(\partial_t\,U^L_{jn}+\partial_1\,U^L_{jn}\right)\left(\partial_t\,w^J_n+\partial_1\,w^J_n\right)
+\nabla_{x'}\,U^L_{jn}\cdot\nabla_{x'}\,w^J_n\right](x,t)dx\nonumber\\
&&\hspace{1in}\quad\quad=0.\label{eq:orthogo12}
\end{eqnarray}
It is easy to check that the left hand sides without taking limit are constant in time. Recalling the definition (\ref{eq:tildewj}) of $\OR{\widetilde{w}}^J_{jn}$ from the profile decompositions, we see by taking $t=t^j_n$ and rescaling that the left hand side of (\ref{eq:orthogo12}) is
\begin{equation*}
\int_{R^d}\,\left[\left(\partial_t\,U^L_j+\partial_1\,U^L_j\right)\left(\partial_t\,\widetilde{w}^J_{jn}+\partial_1\,\widetilde{w}^J_{jn}\right)
+\nabla_{x'}U^L_j\cdot\nabla_{x'}\widetilde{w}^J_{jn}\right](x,0)\,dx\to 0,
\end{equation*}
as $n\to\infty$, since $\OR{\widetilde{w}}^J_{jn}(\cdot,0)\rightharpoonup 0$. The proof of (\ref{eq:orthogo11}) is similar. \\

\hspace{.6cm}The orthogonality conditions can be localized. More precisely we have
\begin{lemma}\label{lm:localizedorthogonality}
Suppose that $(u_{0n},\,u_{1n})$ is a bounded sequence in $\HL$ and has the profile decomposition (\ref{eq:profiledecomposition}). Fix $j\ge 1$ and assume that
\begin{equation*}
\lim_{n\to\infty}\frac{t^j_n}{\lambda^j_n}\in\{\pm\infty\}.
\end{equation*}
Then for any $\beta>1$, 
\begin{equation}\label{eq:localizedenergyexpansion}
\liminf_{n\to\infty}\|(u_{0n},\,u_{1n})\|^2_{\HL\left(\left\{x:\,\frac{|t^j_n|}{\beta}<|x-c^j_n|<\beta |t^j_n|\right\}\right)}\ge \|\OR{U}^L_j(0)\|_{\HL},
\end{equation}
and
\begin{eqnarray}
&&\liminf_{n\to\infty}\int_{\frac{|t^j_n|}{\beta}<|x-c^j_n|<\beta |t^j_n|}\,(u_{1n}+\partial_1u_{0n})^2+|\nabla_{x'}u_{0n}|^2\,dx\nonumber\\
&&\ge\int_{R^d}\,\left[\left(\partial_tU^L_j+\partial_1U^L_j\right)^2+|\nabla_{x'}U^L_j|^2\right](x,0)\,dx.\label{eq:localizedexpansion}
\end{eqnarray}
\end{lemma}
The property (\ref{eq:localizedenergyexpansion}) follows from Lemma \ref{lm:removedcone} below with $E_n=\emptyset$. The proof of (\ref{eq:localizedexpansion}) is very similar, in view of the conservation of the quantity
$\int \left(\partial_tU+\partial_1U\right)^2+|\nabla_{x'}U|^2\,dx$ 
for a solution $U$ of the wave equation, and we omit it. 

\begin{lemma}\label{lm:removedcone}
Suppose that $(u_{0n},\,u_{1n})$ is a bounded sequence in $\HL$ and has the profile decomposition (\ref{eq:profiledecomposition}). Fix $j\ge 1$ and assume that
\begin{equation*}
\lim_{n\to\infty}\frac{t^j_n}{\lambda^j_n}\in\{\pm\infty\}.
\end{equation*}
Assume that a sequence of sets $E_n\subseteq R^d$ satisfies (\ref{bound_E}) for all $M>0$.
Then for any $\beta>1$, 
\begin{equation}\label{eq:coneenergyexpansion}
\liminf_{n\to\infty}\,\|(u_{0n},\,u_{1n})\|^2_{\HL\left(\left\{x:\,\frac{|t^j_n|}{\beta}<|x-c^j_n|<\beta |t^j_n|\right\}{\displaystyle \backslash}\,E_n\right)}\ge \|\OR{U}^L_j(0)\|_{\HL}.
\end{equation}
\end{lemma}

\smallskip
\noindent
{\it Proof.} Assume that $\lim\limits_{n\to\infty}\frac{t^j_n}{\lambda^j_n}=\infty$, the other case is identical. Fix $\beta>1$, denote
\begin{equation}
F^j_n:=\left\{x:\,|x-c^j_n|>\beta\,t^j_n\,\,{\rm or}\,\,|x-c^j_n|<\frac{t^j_n}{\beta}\right\}\bigcup\,E_n.
\end{equation}
By Lemma \ref{lm:concentrationoffreeradiation} and Lemma \ref{lm:absolutecontinuity}, we see that
\begin{equation}\label{eq:VanishingOnFj}
\left\|\OR{U}^L_{jn}(0)\right\|_{\HL\left(F^j_n\right)}\to 0,\,\,{\rm as}\,\,n\to\infty.
\end{equation}
Write 
\begin{equation}
(\widetilde{u}_{0n},\,\widetilde{u}_{1n})=(u_{0n},\,u_{1n})-\OR{U}^L_{jn}(\cdot,0).
\end{equation}
Then by Lemma \ref{lm:absolutecontinuity} and (\ref{eq:VanishingOnFj}), we get that
\begin{eqnarray*}
&&\int_{\{x:\,\frac{t^j_n}{\beta}<|x-c^j_n|<\beta t^j_n\} \backslash\,E_n}\,|u_{1n}|^2+|\nabla u_{0n}|^2\,dx\\
&&=\int_{\{x:\,\frac{t^j_n}{\beta}<|x-c^j_n|<\beta t^j_n\}\backslash\,E_n}\,\left[\left(\partial_tU^L_{jn}\right)^2+|\nabla U^L_{jn}|^2\right](x,0)\,dx+\\
&&\quad\quad+\int_{\{x:\,\frac{t^j_n}{\beta}<|x-c^j_n|<\beta t^j_n\} \backslash\,E_n}|\widetilde{u}_{1n}|^2+|\nabla \widetilde{u}_{0n}|^2\,dx+\\
&&\quad\quad\quad+2\int_{R^d}\,\nabla U^L_{jn}(0)\,\nabla\widetilde{u}_{0n}+\partial_tU^L_{jn}(0)\,\partial_t\widetilde{u}_{1n}\,dx\\
&&\quad\quad\quad\quad\quad+\,\,O\left(\|\OR{U}^L_{jn}(0)\|_{\HL(F^j_n)}\right)\\
&&\ge \int_{R^d}\left[\left(\partial_tU^L_{jn}\right)^2+|\nabla_{x'}U^L_{jn}|^2\right](x,0)\,dx+o_n(1).
\end{eqnarray*}
In the above, we used the pseudo-orthogonality of the profile $\OR{U}^L_{jn}$ and $(\widetilde{u}_{0n},\,\widetilde{u}_{1n})$. The lemma is proved.\\

To apply the linear profile decompositions to the nonlinear equation (\ref{eq:main}), we  need the following perturbation lemma.
\begin{lemma}\label{lm:longtimeperturbation}
Let $I$ be a time interval with $0\in I$. Let $\OR{U}\in C(I,\,\HL)$, with 
\begin{equation*}
\|U\|_{\Snorm(R^d\times I)}\leq M<\infty,
\end{equation*}
 be a solution to 
\begin{equation}\label{eq:equationwitherror}
\partial_{tt}U-\Delta U=|U|^{\dual-2}U+f,\,\,{\rm in}\,\,R^d\times I.
\end{equation}
Then we can find $\epsilon_{\ast}=\epsilon_{\ast}(M,d)>0$ sufficiently small, such that if
\begin{equation}\label{eq:errorissmall}
\|f\|_{L^1_tL^2_x(R^d\times I)}+\|(u_0,u_1)-\OR{U}(0)\|_{\HL}\leq \epsilon\leq \epsilon_{\ast},
\end{equation}
then the unique solution $\OR{u}$ to equation (\ref{eq:main}) with initial data $(u_0,u_1)\in\HL$ exists in $R^d\times I$. Moreover, $\OR{u}$ verifies the estimate
\begin{equation}
\sup_{t\in I}\,\|\OR{u}(t)-\OR{U}(t)\|_{\HL}+\|u-U\|_{\Snorm(R^d\times I)}\leq C(M,d)\epsilon.
\end{equation} 
\end{lemma}

\smallskip
\noindent
{\it Remark.} This perturbation lemma is often termed ``long time perturbation" lemma, due to the fact there is no restriction on the size of the time interval $I$. Such results are well known, see for example Theorem 2.14 in \cite{KMacta} and Lemma 2.1 in \cite{JiaLiuXu}. We omit the standard proof, and refer readers to the presentations in \cite{KMacta,JiaLiuXu}. \\

The perturbation lemma \ref{lm:longtimeperturbation} has important applications to the nonlinear profile decompositions. Assume that the sequence of initial data $(u_{0n},\,u_{1n})$ is uniformly bounded in the energy space with respect to $n$ and that $(u_{0n},\,u_{1n})$ has the profile decomposition (\ref{eq:profiledecomposition}). For each $j$, introduce the nonlinear profile $U_j$ as follows. 
\begin{itemize}
\item if $\lim\limits_{n\to\infty}\frac{t^j_n}{\lambda^j_n}=L\in\{\pm\infty\}$, then define $\OR{U}_j$ as the unique solution to equation (\ref{eq:main}) in a neighborhood of $-L$, with 
\begin{equation*}
\lim_{t\to -L}\|\OR{U}_j(t)-\OR{U}_j^L(t)\|_{\HL}=0.
\end{equation*}
The existence of $\OR{U}_j$ follows from standard perturbative arguments.
\item if $t^j_n\equiv 0$ for all $n$, then define $\OR{U}_j$ as the solution to equation (\ref{eq:main}) with initial data $\OR{U}_j(0)=\OR{U}_j^L(0)$.\\
\end{itemize}

We have the following nonlinear approximation lemma for the linear profile decomposition, see Bahouri and G\'erard \cite{BaGe} for details, see also \cite{DKMsmall}.
\begin{lemma}\label{lm:nonlinearprofiledecomposition}
Let $(u_{0n},\,u_{1n})$ be a sequence of initial data that are uniformly bounded in the energy space. Assume that $(u_{0n},\,u_{1n})$ has the profile decomposition (\ref{eq:profiledecomposition}). Let $\OR{U}_j$ be the nonlinear profile associated with $\OR{U}_j^L$, $\lambda^j_n$, $t^j_n$. Denote $T_+(U_j)$ as the maximal time of existence for $U_j$. 
Let $\theta_n\in (0,\infty)$. Assume that for all $j\ge 1$,  $\frac{\theta_n-t^j_n}{\lambda^j_n}<T_+(U_j)$ for large $n$ and
\begin{equation}
\limsup_{n\to\infty}\left\|U^j\right\|_{\Snorm\left(R^d\times \left(-\frac{t^j_n}{\lambda^j_n},\,\frac{\theta_n-t^j_n}{\lambda^j_n}\right)\right)}<\infty.
\end{equation}
Let $\OR{u}_n$ be the solution to equation (\ref{eq:main}) with initial data $(u_{0n},\,u_{1n})$.
Then for sufficiently large $n$, $\OR{u}_n$ has the following decomposition
\begin{eqnarray*}
\OR{u}_n&=&\sum_{j=1}^J\,\left(\left(\lambda^j_n\right)^{-\frac{d}{2}+1}U_j\left(\frac{x-c^j_n}{\lambda^j_n}, \,\frac{t-t^j_n}{\lambda^j_n}\right),\,\left(\lambda^j_n\right)^{-\frac{d}{2}}\partial_tU_j\left(\frac{x-c^j_n}{\lambda^j_n}, \,\frac{t-t^j_n}{\lambda^j_n}\right)\right)+\\
\\
&&\quad\quad\quad+\,\OR{w}^J_n+\OR{r}^J_n,
\end{eqnarray*}
in $R^d\times [0,\,\theta_n)$, where $\OR{r}_n$ vanishes in the sense that
\begin{equation*}
\lim_{J\to\infty}\,\limsup_{n\to\infty}\left(\sup_{t\in[0,\theta_n)}\|\OR{r}^J_n(t)\|_{\HL}+\|r^J_n\|_{\Snorm(R^d\times [0,\,\theta_n))}\right)=0.
\end{equation*}
\end{lemma}

\medskip

The principle of finite speed of propagation plays an essential role in the study of wave equations. Below we shall use the following version of this principle. Let us set $a\wedge b:=\min\{a,\,b\}$ for any $a,\,b\in R$.
\begin{lemma}
Let $\OR{u},\,\OR{v}\in C([-T,T],\,\HL)$, with $u,\,v\in \Snorm(R^d\times[-T,T])$, be two solutions to equation (\ref{eq:main}), with initial data $(u_0,\,u_1)\in\HL$ and $(v_0,\,v_1)\in\HL$ respectively. Assume that for some $R>0$,
\begin{equation*}
(u_0(x),\,u_1(x))=(v_0(x),\,v_1(x))\,\,\,{\rm in}\,\,B_R.
\end{equation*}
Then 
\begin{equation}
u(x,t)=v(x,t)\,\,\,{\rm for}\,\,(x,t)\in\{(x,t):\,|x|<R-|t|,\,\,|t|<R\wedge T\}.
\end{equation}

\end{lemma}

\medskip

The perturbation lemma \ref{lm:longtimeperturbation} combined with finite speed of propagation implies the following local-in-space approximation result.\\

\begin{lemma}\label{lm:localinspaceapproximation}
Let $\OR{u}_n$ and $\OR{v}_n$ be two sequences of solutions defined in $R^d\times [-T,T]$ with initial data $(u_{0n},\,u_{1n})\in\HL$ and $(v_{0n},\,v_{1n})\in\HL$ respectively. Assume that 
\begin{equation*}
\|u_n\|_{\Snorm(R^d\times [-T,T])}\leq M<\infty, 
\end{equation*}
and that 
\begin{equation*}
\|(u_{0n},\,u_{1n})-(v_{0n},\,v_{1n})\|_{\HL(B_{R+\epsilon})}+\|u_{0n}-v_{0n}\|_{L^{\dual}(B_{R+\epsilon})}\to 0,
\end{equation*}
as $n\to\infty$, for some $R>0$ and $\epsilon>0$ . Then we have for $\tau=R\wedge T$ that

\begin{eqnarray*}
&&\sup_{t\in(-\tau,\,\tau)}\,\|\OR{u}_n(t)-\OR{v}_n(t)\|_{\HL(|x|<R-|t|)}+\\
\\
\quad\quad&&\,\,\,\,\,\,\,\,\,\quad\quad\quad\vspace{2cm}+\,\|u_n-v_n\|_{\Snorm(\{(x,t):\,|x|<R-|t|,\,|t|<\tau\})}\to 0,
\end{eqnarray*}
as $n\to\infty$.
\end{lemma}  

\smallskip
\noindent
{\it Proof.} We define the following modified initial data $(\tilde{v}_{0n},\,\tilde{v}_{1n})$, with
\begin{equation}
(\tilde{v}_{0n},\,\tilde{v}_{1n}):=\left\{\begin{array}{ll}
                              (u_{0n},\,u_{1n})&\,\,{\rm if}\,\,|x|>R+\epsilon,\\
                   &\\
                               \frac{|x|-R}{\epsilon}(u_{0n},\,u_{1n})+\frac{R+\epsilon-|x|}{\epsilon}(v_{0n},\,v_{1n}),&\,\,{\rm if}\,\,R<|x|<R+\epsilon,\\
&\\
                      (v_{0n},\,v_{1n}) &\,\,{\rm if}\,\,|x|<R.
                              \end{array}\right.
\end{equation}
Then it is straightforward to verify that
\begin{equation*}
\|(\tilde{v}_{0n},\,\tilde{v}_{1n})-(u_{0n},\,u_{1n})\|_{\HL(R^d)}\to 0,\,\,{\rm as}\,\,n\to\infty,
\end{equation*}
and that
\begin{equation*}
(v_{0n},\,v_{1n})\equiv (\tilde{v}_{0n},\,\tilde{v}_{1n}),\,\,\,{\rm in}\,\,B_R.
\end{equation*}
Let $\OR{\tilde{v}}_n$ be the solution to equation (\ref{eq:main}) with initial data $\OR{\tilde{v}}_n(0)=(\tilde{v}_{0n},\,\tilde{v}_{1n})$. Then by the principle of finite speed of propagation, we see that 
\begin{equation}\label{eq:localequality}
v_n\equiv \tilde{v}_n, \,\,\,\,{\rm in}\,\, \{(x,t):\,|x|<R-|t|,\,|t|<\tau\}.
\end{equation}
By the perturbation Lemma \ref{lm:longtimeperturbation}, we get that
\begin{equation}\label{eq:differencegoestozero}
\lim_{n\to\infty}\,\left(\sup_{t\in[-T,\,T]}\|\OR{u}_n(t)-\OR{\tilde{v}}_n(t)\|_{\HL(R^d)}+\|u_n-\tilde{v}_n\|_{\Snorm(R^d\times [-T,T])}\right)=0.
\end{equation}
Combining (\ref{eq:differencegoestozero}) with (\ref{eq:localequality}) finishes the proof of the lemma.\\

\section{Control of energy flux and Morawetz identity}
\subsection{Generalities on type II blow-up solutions}
\label{SS:typeII}
In this section we begin the study of the type II singular solutions to equation (\ref{eq:main}). Assume without loss of generality (by a rescaling argument) that $\OR{u}\in C((0,1],\HL)$, with $u\in \Snorm(R^d\times (\epsilon,1])$ for all $\epsilon>0$, is a type II solution to equation (\ref{eq:main}) in $R^d\times (0,1]$ with Cauchy data prescribed at time $t=1$, and that $t=0$ is the blow-up time for $u$, i.e., $u\notin \Snorm(R^d\times (0,1])$. Following \cite[Section 3]{DKMsmall}, we define the set of singular points
\begin{equation}
\mathcal{S}:=\left\{x_{\ast}\in R^d:\,\|u\|_{\Snorm(B_{\epsilon}(x_{\ast})\times (0,\epsilon])}=\infty,\,\,{\rm for\,\,all\,\,}\epsilon>0\right\}.
\end{equation}
By small data global existence theory for equation (\ref{eq:main}) and finite speed of propagation, we see that if $x_{\ast}\in\mathcal{S}$, then there exists a $\epsilon_0(d)>0$, with
\begin{equation}\label{eq:nontrivialcon}
\forall \,\epsilon>0,\,\,\liminf_{t\to 0+}\,\int_{|x-x_{\ast}|<\epsilon}\,\left[\frac{|\nabla u|^2}{2}+\frac{|\partial_tu|^2}{2}+\frac{|u|^{\dual}}{\dual}\right](x,t)\,dx>\epsilon_0.
\end{equation}
We deduce from the assumption that $\OR{u}$ is type II 
that $\mathcal{S}$ is a non-empty finite set $\{x_1,\dots,x_{K}\}$. Let $A_{\eps}=R^d\backslash \bigcup_{j=1}^{K}B_{\epsilon}(x_j)$.
The finite speed of propagation, together with local Cauchy theory for equation (\ref{eq:main}), implies that for any $\epsilon>0$,
 $$u\in  \Snorm\big(A_{\eps}\times (0,1]\big),\quad \OR{u}\in C\left([0,1],\,\HL(A_{\eps})\right).$$
 Thus 
\begin{equation}
(f,g):=w-\lim_{t\to 0+}\,\OR{u}(t)
\end{equation}
is well defined as a function in $\HL(R^d)$. See again \cite{DKMsmall} for the details. \\

\hspace{.6cm}Let $\OR{v}\in C([0,\delta],\,\HL)$ with $v\in \Snorm(R^d\times [0,\delta])$ be the local in time solution to equation (\ref{eq:main}) with initial data $\OR{v}(0)=(f,g)$. Using finite speed of propagation again,  we have that
\begin{equation}\label{eq:regularsingular}
u(x,t)=v(x,t)\,\,{\rm for}\,\,(x,t)\in \Big(R^d\times (0,\delta]\Big)\backslash \bigcup_{j=1}^{K}\{(x,t):\,|x-x_j|\leq t\}.
\end{equation}
$v$ is called the regular part of $u$. (\ref{eq:regularsingular}) shows that there is a sharp division of the singular and regular regions for type II solutions. We again refer the reader to the Section 3 of \cite{DKMsmall} for the details.\\

Below we shall assume that $0\in \mathcal{S}$ and the main focus is to study the asymptotics of $u$ near this singular point. More precisely, we will prove: 
\begin{prop}\label{pr:morawetz}
Let $u$ be as above. Then for any $0<10\,t_1<t_2<\delta$, we have
\begin{equation}
\int_{t_1}^{t_2}\int_{|x|<t}\left(\partial_tu+\frac{x}{t}\cdot\nabla u+\left(\frac{d}{2}-1\right)\frac{u}{t}\right)^2\,dx\,\frac{dt}{t}\leq C(u)\, \left(\log{\frac{t_2}{t_1}}\right)^{\frac{1}{2}},
\end{equation}
for some $C(u)$ depending on $u$ but independent of $t_1,\,t_2$.
\end{prop}
The proof of Proposition \ref{pr:morawetz} uses the energy flux identity:
\begin{multline}
\frac{1}{\sqrt{2}}\int_{t_1}^{t_2}\int_{|x|=t}\,\frac{|\nabla_{x,t}u|^2}{2}+\partial_tu\,\frac{x}{|x|}\cdot\nabla u-\frac{|u|^{\dual}}{\dual}\,d\sigma dt\\
=\int_{|x|\leq t_2}\left[\frac{|\nabla_{x,t}u|^2}{2}-\frac{|u|^{\dual}}{\dual}\right](x,t_2)\,dx
-\int_{|x|\leq t_1}\left[\frac{|\nabla_{x,t}u|^2}{2}-\frac{|u|^{\dual}}{\dual}\right](x,t_1)\,dx.
\label{eq:energyfluxidentity}
\end{multline}
One obstacle in the study of equation (\ref{eq:main}) is that, unlike for linear or defocusing wave equations and for wave maps,  the energy flux (the left-hand side of \eqref{eq:energyfluxidentity}) is in general not positive, because the term $|u|^{\dual}$ comes with a minus sign. The main new input in the proof of Proposition \ref{pr:morawetz} is the following simple observation. Since $v\in \Snorm(R^d\times [0,\delta])$ and $\OR{v}\in L^{\infty}([0,\delta],\,\HL)$, we see that 
\begin{equation}
\left|\nabla_{x,t}\,|v|^{\dual}\right|\lesssim_d |\nabla_{x,t} v||v|^{\dual-1}\in L^1(R^d\times[0,\delta]).
\end{equation}
This regularity property implies that $v$ is well defined on any sufficiently regular hypersurface of $R^d\times R$ as an $L^{\dual}$ function. As a consequence $|u|^{\dual}=|v|^{\dual}$ is integrable on the boundary of the lightcone, the energy-flux term can be controlled, and we can use a variant of the Morawetz identity as in the study of energy critical wave maps \cite{GrillakisEnergy,taoIII,Tataru4}. 

We now turn to the details of the proof of Proposition \ref{pr:morawetz}. In \S \ref{SS:regular}, we estimate boundary terms using the idea described in the preceding paragraph. We introduce self-similar variables and the approximate self-similar energy in \S \ref{SS:SS}. The proof is completed in \S \ref{SS:Morawetz}.

\subsection{Estimates on regular solutions}
\label{SS:regular}
We prove here:
\begin{lemma}
\label{L:trace}
Let $v$ be a solution of \eqref{eq:main} defined in a neighborhood of $R^d\times [0,1]$. Let 
$$\Gamma=\left\{(x,t):\; t=|x|,\; |x|\leq 1\right\}.$$
Let 
$$ M_1=\sup_{0\leq t\leq 1} \|\vec{v}(t)\|_{\HL},\quad M_2=\|v\|_{L^{\frac{d+2}{d-2}}\left([0,1],L^{\frac{2(d+2)}{d-2}}\right)}.$$
Then the restrictions of $|v|^{\frac{2d}{d-2}}$, $\frac 12|\nabla_{x,t}v|^2+\partial_tv\frac{x}{|x|}\nabla v$ and $\frac{1}{|x|^2}v$ to $\Gamma$ are well-defined and the following estimates hold:
\begin{align}
 \label{trace_i} \int_{\Gamma} |v|^{\dual}d\sigma &\leq C(M_1,M_2)\\
\label{trace_ii} \int_{\Gamma}\left[\frac{1}{2}|\nabla_{x,t} v|^2+\partial_tv\frac{x}{|x|}\nabla v\right]d\sigma&\leq C(M_1,M_2)\\
\label{trace_iii} \int_{\Gamma}\frac{1}{|x|^2}|v|^2\,d\sigma&\leq C(M_1,M_2).
\end{align}
\end{lemma}
\begin{remark}
Note that the integrand in \eqref{trace_ii} can be rewritten as:
$$ \frac{1}{2}\left|\spartial v\right|^2+\frac 12 \left| \partial_tv+\frac{x}{|x|}\cdot\nabla v\right|^2,$$
where $|\spartial u|^2$ denotes the tangential derivative $|\nabla u|^2-\left|\frac{x}{|x|}\cdot\nabla u\right|^2$ and is thus nonnegative. 
\end{remark}
\begin{remark}
It suffices to prove Lemma \ref{L:trace} for $v$ with smooth, compactly supported initial data, and then pass to the limit.
\end{remark}
\begin{proof}[Proof of Lemma \ref{L:trace}]
 We first prove \eqref{trace_i}. Note that $\partial_t |v|^{\dual}=\dual|v|^{\frac{4}{d-2}}v\partial_t v$, so that
\begin{multline*}
\int_{0}^{1} \int_{R^d} \left|\partial_t |v|^{\dual}\right|\,dx\,dt\leq \dual \int_{0}^{1} \left(\int_{R^d} |v|^{\frac{2(d+2)}{d-2}}\,dx\right)^{\frac{1}{2}}\left( \int_{R^d}|\partial_t v|^2\,dx \right)^{1/2}\,dt\\
\leq \dual M_1M_2^{\frac{d+2}{d-2}}. 
\end{multline*}
Thus
\begin{multline*}
 \int_{\Gamma}|v|^{\dual}d\sigma=\sqrt{2}\int_{|x|<1} |v(x,|x|)|^{\dual}dx\\
\leq \sqrt{2}\int_{|x|<1}\int_0^{|x|}\left| \partial_t |v(x,t)|^{\dual}  \right|\,dt\,dx+\sqrt{2} \int_{|x|<1}|v(x,0)|^{\dual}\,dx\lesssim M_1M_2^{\frac{d+2}{d-2}}+M_1^{\dual}.
\end{multline*}
To prove \eqref{trace_ii}, we use the flux identity \eqref{eq:energyfluxidentity} on the solution $v$, for $0\leq t_1\leq t_2\leq 1$.
Since $\int_{\Gamma} |v|^{\dual}d\sigma$ is controlled by \eqref{trace_i}, \eqref{trace_ii} follows immediately.

To estimate \eqref{trace_iii}, let $\theta$ be a smooth function, $\theta(x)=1$ for $|x|<1$, $\supp \theta \subset\{|x|<1+\eps\}$ (where $\eps>0$ is small, so that $\{0\leq t\leq 1+\eps\}$ is in the domain of existence of $v$). Consider $f(x)=\theta(x) v(x,|x|)$. We apply the Hardy estimate
$$\int_{R^d} \frac{1}{|x|^2}|f(x)|^2\,dx\lesssim \int_{R^d} |\nabla_x f(x)|^2\,dx.$$
We note that 
$$ \nabla f(x)=\nabla \theta(x)v(x,|x|)+\theta(x)\nabla v(x,|x|)+\theta(x)\frac{x}{|x|} \partial_tv(x,|x|),$$
and thus
$$|\nabla_xf(x)|^2\leq 2|\nabla \theta(x)|^2|v(x,|x|)|^2+2\theta^2(x)\left( |\nabla_{x,t}v(x,|x|)^2+2\partial_t v\frac{x}{|x|}\cdot\partial_x v(x,|x|)|\right).$$
The first term is controlled because on the support, $|x|\leq 1+\eps$, by \eqref{trace_i} and H\"older's inequality. The control of the second term follows from \eqref{trace_ii}. More accurately, we have used a variant of \eqref{trace_i} and \eqref{trace_ii} where $\Gamma$ is replaced by $\{(x,t):\; t=|x|\leq 1+\eps\}$. As a consequence, we obtain \eqref{trace_iii} since $f(x)=v(x,|x|)$ for $|x|\leq 1$. 
\end{proof} 
\subsection{Self-similar variables and energy}
\label{SS:SS}
We consider, as in Subsection \ref{SS:typeII}, a type II blow-up solution $u$, defined for $t\in (0,1]$ and blowing up at $t=0$, such that the origin of $R^d$ is in the singular set of $u$. We denote by $v$ the regular part of $u$, defined in Subsection \ref{SS:typeII}. Rescaling if necessary, we can assume without loss of generality that $[0,1]$ is in the domain of existence of $v$, and that 
$$ |t|<|x|\leq 1\text{ and }0<t\leq 1\Longrightarrow \vec{u}(x,t)=\vec{v}(x,t).$$
Notice that, in light of the fact that $v=u$ on $\Gamma$, Lemma \ref{L:trace} holds for $u$. We introduce, for $|x|<t$, $0<t<1$,
$$ y=\frac{x}{t},\quad s=-\log t, \quad w(y,s)=t^{\frac{d}{2}-1}u(x,t).$$
We let $\rho(y)=(1-|y|^2)^{-1/2}$. As in \cite{KMacta} we obtain the following equation for $s>0$, $|y|<1$ 
\begin{equation}
\label{eq:w}
\partial_s^2w=-\frac{1}{\rho}\opdiv\left( \rho\nabla w-\rho(y\cdot \nabla w) \right)-\frac{d(d-2)}{4}w+|w|^{\dual-2}w-2y\nabla \partial_s w-(d-1)\partial_s w.
\end{equation} 
We note that $\nabla_yw(y,s)=t^{\frac{d}{2}} \nabla_xu(x,t)$, and hence, for $t\in ]0,1]$
$$ \int_{|y|<1}\left|\nabla_yw(y,s)\right|^2\,dy=\int_{|x|<t}|\nabla_xu(x,t)|^2\,dx\leq C(u).$$
Next,
\begin{equation}
\label{dsw}
-\partial_sw(y,s)=\left( \frac{d}{2}-1 \right)t^{\frac{d}{2}-1}u(x,t)+t^{\frac{d}{2}}\partial_t u(x,t)+t^{\frac d2} \frac{x}{t}\cdot \nabla u(x,t). 
\end{equation} 
Observe that 
$$ \int_{|y|<1}t^{d-2}|u(x,t)|^2\,dy=\int_{|x|<t} |u(x,t)|^2\,\frac{dx}{t^2}\lesssim \frac{1}{t^2}\left(\int_{|x|<t}|u(x,t)|^{\dual}\right)^{\frac{d-2}{d}} \left(t^{d}\right)^{\frac{2}{d}}\leq C(u). $$
Thus 
\begin{equation}
 \label{bnd_nabla_w}
\sup_{s>0} \|\nabla_{y,s}w(s)\|_{L^2(\{|y|<1\})}<\infty.
\end{equation} 
We next claim the following trace estimate:
\begin{claim}
 \label{C:trace}
$$\int_0^{\infty}\int_{|y|=1} \left(\partial_sw(y,s)\right)^2d\sigma(y)\,ds<\infty.$$
\end{claim}
\begin{proof}
We first consider the term $t^{\frac{d}{2}-1}u(x,t)$ in \eqref{dsw}. By a change of variable, we have:
$$ \int_0^{\infty}\int_{|y|=1}t^{\frac{d}{2}-1}u^2(x,t)d\sigma(y)\,ds=\int_{\Gamma}\frac{1}{|x|^2} u^2\,d\sigma,$$
so this term is controlled by \eqref{trace_iii} in Lemma \ref{L:trace}. 

Next consider the term $t^{\frac{d}{2}}\partial_t u(x,t)+t^{\frac d2} \frac{x}{t}\cdot \nabla u(x,t)$ in \eqref{dsw}. By the same change of variable, we are reduced to bound
$$ \int_{\Gamma} \left(\partial_tu +\frac{x}{|x|}\cdot\nabla_xu\right)^2\,d\sigma$$
which is finite by \eqref{trace_ii} in Lemma \ref{L:trace}. This concludes the proof of the claim. 
\end{proof}
Define $\rho_{\eps}(y)=(1+\eps^2-|y|^2)^{-1/2}$ and, for $s>0$, the \emph{regularized self-similar energy}:
$$ I_{\eps}(s)=\int_{|y|<1} \left(\frac{1}{2}(\partial_sw)^2+\frac 12\left(|\nabla w|^2-(y\cdot \nabla w)^2\right)+\frac{d(d-2)}{8}w^2-\frac{1}{\dual}|w|^{\dual}\right)\rho_{\eps}\,dy.$$
By direct computation, multiplying equation \eqref{eq:w} by $\partial_s w\,\rho_{\eps}$ and integrating by parts, we obtain:
\begin{claim}
 \label{C:selfsim}
Let $0<s_2<s_1$. Then
\begin{multline}
 \label{se_si_energy}
I_{\eps}(s_1)-I_{\eps}(s_2)
=(1+\eps^2)\int_{s_2}^{s_1} (\partial_sw)^2\rho_{\eps}^3\,ds\\-\frac{1}{\eps}\int_{s_2}^{s_1}\int_{|y|=1}(\partial_sw)^2d\sigma(y)\,ds +\eps^2\int_{s_2}^{s_1}\int_{|y|<1} y\cdot \nabla w\partial_sw\rho_\eps^3\,dy\,ds.
\end{multline} 
\end{claim}

\subsection{Proof of the Morawetz estimate}
\label{SS:Morawetz}
We are now ready to prove Proposition \ref{pr:morawetz}. We use the self-similar change of variables and the notations of the preceding subsection. Let $s_1=-\log t_1$ and $s_2=-\log t_2$. In view of formula \eqref{dsw}, we have
\begin{multline*}
\int_{t_1}^{t_2} \int_{|x|<t} \left( \partial_tu+\frac{x}{t}\nabla u(x,t)+\frac{1}{2 t}u(x,t)\right)^2\,dx\,\frac{dt}{t}
\\=\int_{s_2}^{s_1} \int_{|y|<1} \left( \partial_sw(y,s) \right)^2\frac{1}{t^N}\,dx\,ds=\int_{s_2}^{s_1} \int_{|y|<1} \left( \partial_sw(y,s) \right)^2dy\,ds. 
\end{multline*}
We now let $\eps>0$, to be chosen later. We use Claim \ref{C:selfsim}.

By Claim \ref{C:trace}, 
$$\int_{s_2}^{s_1}\int_{|y|=1}(\partial_sw)^2d\sigma(y)\,ds\leq C(u).$$
Furthermore, if $0<\eps<1$,
\begin{equation}
\label{rho_eps}
\frac{1}{\sqrt{2}}\leq \rho_{\eps}\leq \frac{1}{\eps}. 
\end{equation} 
By the bound \eqref{bnd_nabla_w}, we have $|I_{\eps}(s)|\leq \frac{C(u)}{\eps}$. In view of Claim \ref{C:selfsim}, we obtain
\begin{multline*}
\int_{s_2}^{s_1} \int_{|y|<1} (\partial_sw)^2\rho_{\eps}^3\leq \frac{C(u)}{\eps}+\eps^2\int_{s_2}^{s_1}\int_{|y|<1} |\nabla_y w||\partial_sw|\rho_{\eps}^3\\
\leq \frac{C(u)}{\eps}+\frac{\eps^4}{2}\int_{s_2}^{s_1}|\nabla_y w|^2\rho_{\eps}^3+\frac 12 \int_{s_2}^{s_1} \int_{|y|<1} (\partial_sw)^2\rho_{\eps}^3 ,
\end{multline*}
so that, using \eqref{bnd_nabla_w} and \eqref{rho_eps} again,
$$ \int_{s_2}^{s_1} \int_{|y|<1} (\partial_sw)^2\rho_{\eps}^3\leq C(u)\left( \frac{1}{\eps}+\eps(s_1-s_2) \right).$$
We now choose $\eps=\frac{1}{(s_1-s_2)^{1/2}}$ and obtain
$$\int_{s_2}^{s_1} (\partial_sw)^2\,dy\,ds\leq C(u)(s_1-s_2)^{1/2},$$
which is exactly the conclusion of Proposition \ref{pr:morawetz} since
$$s_1-s_2=\log t_2-\log t_1=\log \frac{t_2}{t_1}.$$
\qed
\section{Applications of the Morawetz inequality}
The Morawetz estimate implies the following vanishing condition.
\begin{lemma}\label{lm:keyvanishingMorawetz}
Let $u$ be as in the last section. Then there exist $\mu_j\downarrow 0$, $t_j\in (\frac{4}{3}\mu_j,\,\frac{13}{9}\mu_j)$ and $t_j'\in (\frac{14}{9}\mu_j,\,\frac{5}{3}\mu_j)$, such that
$$\frac{1}{\mu_j}\int_{\mu_j}^{2\mu_j}\int_{|x|<t}\,\left(\partial_tu+\frac{x}{t}\cdot\nabla u+\left(\frac{d}{2}-1\right)\frac{u}{t}\right)^2\,dx \,dt\to 0,\,\,{\rm as}\,\,j\to \infty,$$
and that, for all $C>0$,
\begin{eqnarray*}
&&\sup_{0<\tau<\frac{t_j}{16}}\frac{1}{\tau}\int_{|t_j-t|<\tau}\,\int_{|x|<Ct}\,\left(\partial_tu+\frac{x}{t}\cdot\nabla u+\left(\frac{d}{2}-1\right)\frac{u}{t}\right)^2\,dx\,dt\\
&&\quad\quad+\sup_{0<\tau<\frac{t_j'}{16}}\frac{1}{\tau}\int_{|t_j'-t|<\tau}\,\int_{|x|<Ct}\,\left(\partial_tu+\frac{x}{t}\cdot\nabla u+\left(\frac{d}{2}-1\right)\frac{u}{t}\right)^2\,dx\,dt\\
\\
&&\quad\quad\quad\,\,\,\,\to 0,\,\,{\rm as}\,\,j\to \infty.
\end{eqnarray*}
\end{lemma}

\smallskip
\noindent
{\it Proof.} Recall the following bound for type II solution $\OR{u}(t)$ in $R^d\times (0,1]$ with $0\in\mathcal{S}$, proved in the last section:
\begin{equation}\label{eq:Morawetz2}
\int_{t_1}^{t_2}\int_{|x|<t}\,\left(\partial_tu+\frac{x}{t}\cdot\nabla u+\left(\frac{d}{2}-1\right)\frac{u}{t}\right)^2\,dx\frac{dt}{t}\leq C \left(\log{\frac{t_2}{t_1}}\right)^{\frac{1}{2}}.
\end{equation}
For each large natural number $J$, applying inequality (\ref{eq:Morawetz2}) to $t_1=4^{-J}$ and $t_2=2^{-J}$,  we get that
\begin{equation*}
\sum_{j=0}^{J-1}\,\int_{2^j\,4^{-J}}^{2^{j+1}\,4^{-J}}\int_{|x|<t}\,\left(\partial_tu+\frac{x}{t}\cdot\nabla u+\left(\frac{d}{2}-1\right)\frac{u}{t}\right)^2\,dx\,\frac{dt}{t}\leq C J^{\frac{1}{2}}.
\end{equation*}
Hence, there exists $0\leq j\leq J-1$, such that 
\begin{equation*}
\int_{2^j\,4^{-J}}^{2^{j+1}\,4^{-J}}\int_{|x|<t}\,\left(\partial_tu+\frac{x}{t}\cdot\nabla u+\left(\frac{d}{2}-1\right)\frac{u}{t}\right)^2\,dx\,\frac{dt}{t}\leq C J^{-\frac{1}{2}}.
\end{equation*}
For a decreasing subsequence $\mu_j\to 0$ of $2^j\,4^{-J}$, we get that
\begin{equation}
\frac{1}{\mu_j}\int_{\mu_j}^{2\mu_j}\int_{|x|<t}\,\left(\partial_tu+\frac{x}{t}\cdot\nabla u+\left(\frac{d}{2}-1\right)\frac{u}{t}\right)^2\,dx \,dt\to 0,\,\,{\rm as}\,\,j\to \infty.
\end{equation}
Let 
\begin{equation}
g(t)=\int_{|x|<t}\,\left(\partial_tu+\frac{x}{t}\cdot\nabla u+\left(\frac{d}{2}-1\right)\frac{u}{t}\right)^2\,dx,
\end{equation}
then $\frac{1}{\mu_j}\int_{\mu_j}^{2\mu_j}\,g(t)\,dt\to 0+$ as $j\to\infty$. Denote $M(g\chi_{(\mu_j,\,2\mu_j)})$ as the Hardy-Littlewood maximal function of $g\chi_{(\mu_j,\,2\mu_j)}$. Passing to a subsequence and renumbering the indices, we can assume that 
\begin{equation*}
\frac{1}{\mu_j}\int_{\mu_j}^{2\mu_j}\,g(t)\,dt\leq 4^{-j}.
\end{equation*}
Then using the fact that the Hardy-Littlewood maximal operator $M$ is bounded from $L^1$ to $L^{1,\infty}$, we get that
\begin{equation*}
\left|\{t\in (\mu_j,\,2\mu_j):\,M(g\chi_{(\mu_j,\,2\mu_j)})(t)>2^{-j}\}\right|\leq C 2^{-j}\mu_j.
\end{equation*}
Thus we can find sequences $t_j$ and $t_j'$ such that $t_j\in (\frac{4}{3}\mu_j,\,\frac{13}{9}\mu_j)$ and $t_j'\in (\frac{14}{9}\mu_j,\,\frac{5}{3}\mu_j)$, such that $M(g\chi_{(\mu_j,\,2\mu_j)})(t_j)\to 0+$ and $M(g\chi_{(\mu_j,\,2\mu_j)})(t_j')\to 0+$ as $j\to\infty$. Consequently along sequences $0<t_j<t_j'$, we have 
\begin{equation}
\sup_{0<\tau<\frac{t_j}{16}}\frac{1}{\tau}\int_{|t_j-t|<\tau}\,g(t)\,dt+\sup_{0<\tau<\frac{t_j'}{16}}\frac{1}{\tau}\int_{|t_j'-t|<\tau}\,g(t)\,dt\to 0,\,\,{\rm as}\,\,j\to \infty.
\end{equation}
That is 
\begin{eqnarray*}
&&\sup_{0<\tau<\frac{t_j}{16}}\frac{1}{\tau}\int_{|t_j-t|<\tau}\,\int_{|x|<t}\,\left(\partial_tu+\frac{x}{t}\cdot\nabla u+\left(\frac{d}{2}-1\right)\frac{u}{t}\right)^2\,dx\,dt\\
&&\quad\quad+\sup_{0<\tau<\frac{t_j'}{16}}\frac{1}{\tau}\int_{|t_j'-t|<\tau}\,\int_{|x|<t}\,\left(\partial_tu+\frac{x}{t}\cdot\nabla u+\left(\frac{d}{2}-1\right)\frac{u}{t}\right)^2\,dx\,dt\\
\\
&&\quad\quad\quad\,\,\to 0,\,\,{\rm as}\,\,j\to \infty.
\end{eqnarray*}
Fix $C>0$. Since $u(x,t)=v(x,t)$ for $(x,t)\in\{(x,t):\,x\in B_{\delta},\,|x|>t,\,t\in(0,\delta)\}$ and $v$ is regular, we can replace $|x|<t$ by $|x|<Ct$ in the domains of integration above. The lemma is proved.\\

\hspace{.6cm}The reason why we need two sequences $t_j$, $t_j'$ approaching zero will become clear below (see \eqref{eq:secondvanishing} and after).\\

\section{Characterization of solutions along a time sequence with the vanishing condition}

In this section, we use the asymptotic vanishing of the quantity 
\begin{equation}
\int_{|x|<Ct}\,\left(\partial_tu+\frac{x}{t}\cdot\nabla u+\left(\frac{d}{2}-1\right)\frac{u}{t}\right)^2\,dx,
\end{equation}
to obtain a preliminary decomposition along a sequence of times in the spirit of decomposition (\ref{eq:maindecompositionhaha}), albeit with a remainder term that vanishes only in $L^{\dual}$.  Suppose that $\OR{u}\in C((0,1],\HL)$, with $u\in \Snorm(R^d\times (\epsilon,1])$ for any $\epsilon>0$, is a type II solution to equation (\ref{eq:main}). We assume that $x_{\ast}=0$ is a singular point. Then by Lemma \ref{lm:keyvanishingMorawetz}, there exists a sequence of times $t_n\downarrow 0$, such that, for all $C>0$,
\begin{equation}\label{eq:vanishingcharacterization}
\lim_{n\to\infty}\,\sup_{0<\tau<\frac{t_n}{16}}\frac{1}{\tau}\int_{|t_n-t|<\tau}\,\int_{|x|<Ct}\,\left(\partial_tu+\frac{x}{t}\cdot\nabla u+\left(\frac{d}{2}-1\right)\frac{u}{t}\right)^2\,dx\,dt = 0.
\end{equation}

\medskip

Our main goal in this section is to prove the following theorem.
\begin{theorem}\label{th:mainpreliminary}
Let $\OR{u}$, $t_n$ be as above. Let $\OR{v}\in C([0,\,\delta],\HL)$, with $v\in \Snorm(R^d\times [0,\delta])$ for some $\delta>0$, be the regular part of $\OR{u}$. Then passing to a subsequence there exist an integer $J_0\ge 0$\,\footnote{If $J_0=0$, then there is no sum of solitons. We shall show in Section 7 that $J_0\ge 1$.}, $r_0>0$, scales $\lambda_n^j$ with $0<\lambda_n^j\ll t_n$, positions $c_n^j\in R^d$ satisfying $c_n^j\in B_{\beta\, t_n}$ for some $\beta\in(0,1)$, with $\ell_j=\lim\limits_{n\to\infty}\frac{c_n^j}{t_n}$ well defined, and traveling waves $Q_{\ell_j}^j$, for $1\leq j\leq J_0$, such that
\begin{eqnarray}
\OR{u}(t_n)&=&\sum_{j=1}^{J_0}\,\left((\lambda_n^j)^{-\frac{d}{2}+1}\, Q_{\ell_j}^j\left(\frac{x-c_n^j}{\lambda_n^j},\,0\right),\,(\lambda_n^j)^{-\frac{d}{2}}\, \partial_tQ_{\ell_j}^j\left(\frac{x-c_n^j}{\lambda_n^j},\,0\right)\right)+\nonumber\\
\nonumber\\
&&\quad\quad\quad\,\,+\,\OR{v}(t_n)+(\epsilon_{0n},\epsilon_{1n}),\label{eq:maindecompositionsec5}
\end{eqnarray}
where $(\epsilon_{0n},\epsilon_{1n})$ vanishes asymptotically in the following sense:
\begin{equation*}
\|\epsilon_{0n}\|_{L^{\dual}(|x|\leq r_0)}\to 0,\quad {\rm as}\,\,n\to\infty;
\end{equation*}
if we write  
\begin{equation*}
(\epsilon_{0n},\,\epsilon_{1n})(x)
=\left((\lambda^j_n)^{-\frac{d}{2}+1}\widetilde{\epsilon}^{\,j}_{0n}\left(\frac{x-c^j_n}{\lambda^j_n}\right),\,(\lambda^j_n)^{-\frac{d}{2}}\widetilde{\epsilon}^{\,j}_{1n}\left(\frac{x-c^j_n}{\lambda^j_n}\right)\right),
\end{equation*}
then $(\widetilde{\epsilon}^{\,j}_{0n},\widetilde{\epsilon}^{\,j}_{1n})\rightharpoonup 0$ as $n\to\infty$, for each $j\leq J_0.$
In addition, the parameters $\lambda_n^j,\,c^j_n$ satisfy the pseudo-orthogonality condition
\begin{equation}\label{eq:pseudo}
\frac{\lambda_n^j}{\lambda_n^{j'}}+\frac{\lambda_n^{j'}}{\lambda_n^j}+\frac{\left|c_n^j-c_n^{j'}\right|}{\lambda_n^j}\to\infty,
\end{equation}
as $n\to\infty$, for $1\leq j\neq j'\leq J_0.$
\end{theorem}

\medskip
\noindent
{\it Proof.} Fix $\phi\in C_c^{\infty}(B_4)$ with $\phi|_{B_3}\equiv 1$. Let 
\begin{equation*}
(u_{0n},\,u_{1n})=\OR{u}(t_n)\phi\left(\frac{x}{t_n}\right).
\end{equation*}
Clearly $(u_{0n},\,u_{1n})$ is a bounded sequence in $\HL$.\\

Passing to a subsequence if necessary, we can assume that $(u_{0n},\,u_{1n})$ has the profile decomposition 
\begin{eqnarray*}
&&\hspace{2in}(u_{0n},\,u_{1n})=\\
&&\quad\,\sum_{j=1}^J\,\left(\left(\lambda^j_n\right)^{-\frac{d}{2}+1}U_j^L\left(\frac{x-c^j_n}{\lambda^j_n}, \,-\frac{t^j_n}{\lambda^j_n}\right),\,\left(\lambda^j_n\right)^{-\frac{d}{2}}\partial_tU_j^L\left(\frac{x-c^j_n}{\lambda^j_n}, \,-\frac{t^j_n}{\lambda^j_n}\right)\right)+\\
&&\hspace{1in}\quad\quad\quad\,\,\,\,\,\,\,\,+\,\left(w^J_{0n},\,w^J_{1n}\right),
\end{eqnarray*}
where the parameters satisfy
\begin{equation}\label{profileorthogonality3}
t^j_n\equiv 0\,\,{\rm for\,\,all\,\,}n,\,\,\,\,{\rm or}\,\,\lim_{n\to\infty}\frac{t^j_n}{\lambda^j_n}\in\{\pm \infty\},
\end{equation}
and for $j\neq j'$
\begin{equation}\label{eq:profileorthogonality4}
\lim_{n\to\infty}\,\left(\frac{\lambda^j_n}{\lambda_n^{j'}}+\frac{\lambda^{j'}_n}{\lambda^j_n}+\frac{\left|c^j_n-c^{j'}_n\right|}{\lambda^j_n}+\frac{\left|t^j_n-t^{j'}_n\right|}{\lambda^j_n}\right)=\infty.
\end{equation}
We note that $(u_{0n},\,u_{1n})$ has the following concentration property
\begin{equation}\label{eq:concentrationproperty}
\lim_{n\to\infty}\,\|(u_{0n},\,u_{1n})\|_{\HL(|x|\ge t_n)}=0.
\end{equation}
By the concentration property (\ref{eq:concentrationproperty}) of $(u_{0n},\,u_{1n})$ and the discussion on pages 144-145 of \cite{BaGe}, we have $\lambda^j_n\lesssim t_n$, and $\lim\limits_{n\to\infty}\frac{\lambda^j_n}{t_n}\neq 0$ for at most one $j$, say $j=j_0$, $\left|c^j_n\right|\lesssim t_n$, $\left|t^j_n\right|\lesssim t_n$. By passing to a subsequence, we can choose $\lambda^{j_0}_n=t_n$, and let $\ell_j=\lim\limits_{n\to\infty}\frac{c^j_n}{t_n}$ be well-defined for each $j$. \\

We divide the profiles into three cases.\\

\begin{itemize}

\item {\it Case I:} \,\,$t^{j_0}_n\equiv 0$, $\lambda^{j_0}_n\equiv t_n$. We shall show that the nonlinear profile $U_{j_0}$ is a compactly supported self similar solution, a case ruled out by \cite{KMacta};\\

\item {\it Case II:} \,\, $t^j_n\equiv 0$, $\lambda^j_n\ll t_n$. In this case, we will show that $|\ell_j|<1$ and $\OR{U}_j^L(\cdot,0)=\OR{Q}^j_{\ell_j}$;\\

\item {\it Case III:}\,\, $\lambda^j_n\ll \left|t^j_n\right|$. In this case, we absorb these profiles into the error term. \\

\end{itemize}

\noindent{\it Case I:} To simplify notation, let $U^L=U_{j_0}^L$. Note that we can assume that $c^{j_0}_n\equiv 0$, and that ${\rm supp}\,\OR{U}^L\subseteq \overline{B_1}$. Indeed, $t^{j_0}_n\equiv 0$, $\lambda^{j_0}_n\equiv t_n$ and $\frac{c^{j_0}_n}{t_n}\to \ell_{j_0}$, so that
\begin{equation*}
\frac{1}{(\lambda^{j_0}_n)^{\frac{d}{2}-1}}U^L_{j_0}\left(\frac{x-c^{j_0}_n}{\lambda^{j_0}_n},\,0\right)=\frac{1}{t_n^{\frac{d}{2}-1}}U^L_{j_0}\left(\frac{x-c^{j_0}_n}{t_n},\,0\right).\\
\end{equation*}

Let
\begin{equation*}
\OR{\widetilde{U}}^L_{j_0}(x,0)=\OR{U}^L_{j_0}\left(x-\ell_{j_0},\,0\right).\\
\end{equation*}

Now, 
\begin{eqnarray*}
&&\frac{1}{t_n^{\frac{d}{2}-1}}\widetilde{U}^L_{j_0}\left(\frac{x}{t_n},\,0\right)=\frac{1}{t_n^{\frac{d}{2}-1}}U^L_{j_0}\left(\frac{x}{t_n}-\ell_{j_0},\,0\right)=\\
&&\quad\quad=\frac{1}{t_n^{\frac{d}{2}-1}}U^L_{j_0}\left(\frac{x}{t_n}-\frac{c^{j_0}_n}{t_n},\,0\right)+o_n(1),\quad {\rm in}\,\,\dot{H}^1,
\end{eqnarray*}
and similarly for the time derivatives, so that we can replace $\left\{U^L_{j_0},\,\left(\lambda^{j_0}_n,\,c^{j_0}_n,\,t^{j_0}_n\right)\right\}$ with $\left\{\widetilde{U}^L_{j_0},\,\left(t_n,\,0,\,0\right)\right\}$. Thus, since $\widetilde{U}^L_{j_0},\,\left(t_n,\,0,\,0\right)$ is now a profile, the properties of the profile decomposition imply that 
\begin{equation*}
\left(t_n^{\frac{d}{2}-1}u_{0n}(t_n\,x),\,t_n^{\frac{d}{2}}u_{1n}(t_n\,x)\right)\rightharpoonup \OR{\widetilde{U}}^L_{j_0}(0)
\end{equation*}
weakly in $\HL$ and so (\ref{eq:concentrationproperty}) shows ${\rm supp}\,\OR{\widetilde{U}}^L_{j_0}\subseteq \overline{B_1}$. For convenience, we now ignore the\, $\widetilde{}$\, in $\widetilde{U}^L_{j_0}$ and simply call it $U^L$, with the associated non-linear profile $U$.
 Hence, there exist $T>0$, $M>0$, such that $U$ exists in $R^d\times[-T,T]$ and $\|U\|_{\Snorm(R^d\times[-T,T])}\leq M<\infty$. By the Pythagorean expansion of the linear energy, for $\epsilon>0$ small, chosen from the small data theory for equation (\ref{eq:main}), there exists $J(\epsilon)$ such that
\begin{equation}\label{eq:allsmallprofiles}
\sum_{j\ge J(\epsilon)}\|\OR{U}^L_j(0)\|_{\HL}^2<\epsilon^2.
\end{equation}
We will see that the profiles with $j\ge J(\epsilon)$ can be controlled by the Approximation Lemma \ref{lm:nonlinearprofiledecomposition}, since they are small. For the remaining profiles, passing to a subsequence if necessary, we can assume that $\lim\limits_{n\to\infty}\frac{t^j_n}{t_n}$ exists and is finite. We can also assume, by reordering the profiles, that $j_0=1$. We thus let $j$ such that $1< j<J(\epsilon)$, and consider two categories of $j$'s. \\

\begin{itemize}

\item Category 1:\,\,  $1< j<J(\epsilon)$, $\lim\limits_{n\to\infty}\frac{t^j_n}{t_n}\neq 0$. We then claim that for some $\gamma>1$ sufficiently large, we have
\begin{equation*}
\left\|\left(\lambda^j_n\right)^{-\frac{d}{2}+1}U_j^L\left(\frac{x-c^j_n}{\lambda^j_n}, \,\frac{t-t^j_n}{\lambda^j_n}\right)\right\|_{\Snorm\left(R^d\times(-\frac{t_n}{\gamma},\,\frac{t_n}{\gamma})\right)}\to 0,
\end{equation*}
as $n\to\infty$. \\

To see this, note that $t^j_n\not\equiv 0$, so that $\lim\limits_{n\to\infty}\left|\frac{t^j_n}{\lambda^j_n}\right|=\infty$. Assume that $\lim\limits_{n\to\infty}\frac{t^j_n}{t_n}=\alpha$, and (say) $\alpha>0$. Then, 
$$t_n=\frac{t^j_n}{\alpha}+o_n(1)t^j_n.$$
Thus, 
$$\frac{1}{\lambda^j_n}\left(\frac{t_n}{\gamma}-t^j_n\right)=\left[\frac{1}{\alpha \gamma}-\left(1-\frac{o_n(1)}{\gamma}\right)\right]\frac{t^j_n}{\lambda^j_n},$$
and
$$\frac{1}{\lambda^j_n}\left(-\frac{t_n}{\gamma}-t^j_n\right)=\left[-\frac{1}{\alpha \gamma}-\left(1+\frac{o_n(1)}{\gamma}\right)\right]\frac{t^j_n}{\lambda^j_n},$$
which gives the claim. 
These profiles have asymptotically vanishing interaction with the profile $U$, and can be controlled perturbatively at least for a short time that is comparable to $t_n$, as we will see.\\

\item Category 2:\,\, $1<j<J(\epsilon)$, $t^j_n=o(t_n)$. This includes the $j$'s for which $t^j_n\equiv 0$ and hence $\lambda^j_n\ll t_n$. (Note that since $j_0=1$ is the unique $j_0$ with $\lim\limits_{n\to\infty}\frac{\lambda^{j_0}_n}{t_n}\neq 0$, if $j>1$, then $\lambda^j_n\ll t_n$). Recall that $\lim\limits_{n\to\infty}\frac{c^j_n}{t_n}=\ell_j$ is well defined. These profiles can have nontrivial interaction with the first profile, and we will remove them. We denote this set of $j$ as $\mathcal{J}_1$.\\

\end{itemize}

Recall that
\begin{equation*}
\OR{U}^L_{jn}(\cdot,0):=\left(\left(\lambda^j_n\right)^{-\frac{d}{2}+1}U_j^L\left(\frac{x-c^j_n}{\lambda^j_n}, \,-\frac{t^j_n}{\lambda^j_n}\right),\,\left(\lambda^j_n\right)^{-\frac{d}{2}}\partial_tU_j^L\left(\frac{x-c^j_n}{\lambda^j_n}, \,-\frac{t^j_n}{\lambda^j_n}\right)\right).\\
\end{equation*}

We now claim that for any $\tau>0$ and $j\in \mathcal{J}_1$, we have
\begin{equation}
\left\|\OR{U}^L_{jn}(0)\right\|_{\HL\left(R^d\backslash\,B_{\tau t_n}(\ell_jt_n)\right)}+\|U^L_{jn}(0)\|_{L^{\dual}\left(R^d\backslash\,B_{\tau t_n}(\ell_jt_n)\right)}\longrightarrow 0,\,\,\,\,{\rm as}\,\,n\to\infty.\label{eq:closenesslocal}
\end{equation}
Note that, because of  Lemma \ref{lm:concentrationoffreeradiation}, if we are in the case when $\frac{\left|t^j_n\right|}{\lambda^j_n}\to+\infty$, then $\OR{U}^L_{jn}(\cdot,0)$ is concentrated where 
$$\left||x-c^j_n|-|t^j_n|\right|\leq \overline{\lambda}\,\lambda^j_n,$$
where $\overline{\lambda}$ is large but fixed. Note that if $|x-\ell_jt_n|>\tau \,t_n$, since $\ell_j=\lim\limits_{n\to\infty}\frac{c^j_n}{t_n}$, for $n$ large $$\left|x-c^j_n\right|\ge \frac{\tau}{2}t_n,$$
and since $\frac{|t^j_n|}{t_n}=o_n(1)$, if $\left||x-c^j_n|-|t^j_n|\right|\leq \overline{\lambda}\,\lambda^j_n,$ we see that 
$$\overline{\lambda}\,\lambda^j_n\ge \frac{\tau}{2}t_n-|t^j_n|=\left(\frac{\tau}{2}-o_n(1)\right)t_n,$$
which contradicts $\lambda^j_n\ll t_n$, and gives (\ref{eq:closenesslocal}) in this case. If $t^j_n\equiv 0$, recalling that $\lambda^j_n\ll t_n$, we easily obtain (\ref{eq:closenesslocal}).\\

Fix $J$, and define 
\begin{eqnarray*}
&&\hspace{1in}(h_{0n},\,h_{1n})=\left(t_n^{-\frac{d}{2}+1}U_0\left(\frac{x}{t_n}\right),\,t_n^{-\frac{d}{2}}U_1\left(\frac{x}{t_n}\right)\right)+\\
\\
&&\quad+\sum_{1<j\leq J,\,j\not\in \mathcal{J}_1}\left(\left(\lambda^j_n\right)^{-\frac{d}{2}+1}U_j^L\left(\frac{x-c^j_n}{\lambda^j_n}, \,-\frac{t^j_n}{\lambda^j_n}\right),\,\left(\lambda^j_n\right)^{-\frac{d}{2}}\partial_tU_j^L\left(\frac{x-c^j_n}{\lambda^j_n}, \,-\frac{t^j_n}{\lambda^j_n}\right)\right)\\
&&\\
&&\hspace{1in}\quad\quad\quad\,\,\,\,+\,(w^J_{0n},\,w^J_{1n}).
\end{eqnarray*}
Note that this gives a profile decomposition for $(h_{0n},\,h_{1n})$.\\

Consider the rescaled sequences
\begin{eqnarray*}
&&(\widetilde{u}_{0n},\,\widetilde{u}_{1n})=\left(t_n^{\frac{d}{2}-1}\,u_{0n}(t_n\,x),\,t_n^{\frac{d}{2}}\,u_{1n}(t_n\,x)\right);\\
&&(\widetilde{h}_{0n},\,\widetilde{h}_{1n})=\left(t_n^{\frac{d}{2}-1}\,h_{0n}(t_n\,x),\,t_n^{\frac{d}{2}}\,h_{1n}(t_n\,x)\right).
\end{eqnarray*}
Hence
\begin{eqnarray*}
\left(\widetilde{h}_{0n},\,\widetilde{h}_{1n}\right)&=&(U_0,\,U_1)+\,\left(t_n^{\frac{d}{2}-1}\,w^J_{0n}(t_nx),\,t_n^{\frac{d}{2}}\,w^J_{1n}(t_nx)\right)\\
&&\quad\quad\quad+\,\sum_{1<j\leq J,\,j\not\in \mathcal{J}_1}\left(\left(t_n^{-1}\lambda^j_n\right)^{-\frac{d}{2}+1}U_j^L\left(\frac{x-t_n^{-1}c^j_n}{t_n^{-1}\lambda^j_n}, \,-\frac{t_n^{-1}t^j_n}{t_n^{-1}\lambda^j_n}\right),\right.\\
&&\quad\quad\quad\quad\quad\quad\quad\quad\quad\quad\left.\left(t_n^{-1}\lambda^j_n\right)^{-\frac{d}{2}}\partial_tU_j^L\left(\frac{x-t_n^{-1}c^j_n}{t_n^{-1}\lambda^j_n}, \,-\frac{t_n^{-1}t^j_n}{t_n^{-1}\lambda^j_n}\right)\right).\\
\end{eqnarray*}

Let $\widetilde{u}_n$ and $\widetilde{h}_n$ be the solutions to equation (\ref{eq:main}) with initial data $(\widetilde{u}_{0n},\,\widetilde{u}_{1n})$ and $(\widetilde{h}_{0n},\,\widetilde{h}_{1n})$ respectively. \\

      By (\ref{eq:closenesslocal}), and a rescaling, we see that
\begin{eqnarray}
&&\hspace{.5in}\left\|\left(\widetilde{u}_{0n},\,\widetilde{u}_{1n}\right)-\left(\widetilde{h}_{0n},\,\widetilde{h}_{1n}\right)\right\|_{\HL\left(R^d\backslash\,\bigcup_{j\in\mathcal{J}_1}B_{\tau }(\ell_j)\right)}\nonumber\\
&&\hspace{.7in}+\,\left\|\widetilde{u}_{0,n}-\widetilde{h}_{0,n}\right\|_{L^{\dual}\left(R^d\backslash\,\bigcup_{j\in\mathcal{J}_1}B_{\tau }(\ell_j)\right)}\to 0,\label{eq:rescaledcloseness}
\end{eqnarray}
as $n\to\infty$, for any $\tau>0$.\\

    We now use Lemma \ref{lm:nonlinearprofiledecomposition}, with $\theta_n=T_1$, $T_1$ small and the estimate for $j\not\in\mathcal{J}_1$, $j\leq \mathcal{J}(\epsilon)$ in Category 1, together with (\ref{eq:allsmallprofiles}) and the small data theory, to conclude that $\OR{\widetilde{h}}_n$ exists in $R^d\times [-T_1,\,T_1]$, with bound
\begin{equation*}
\|\widetilde{h}_n\|_{\Snorm\left(R^d\times[-T_1,T_1]\right)}\leq M_1<\infty.
\end{equation*}
Fix any $R>0$, $\sigma>0$ and consider $y$ such that 
\begin{equation*}
\inf\,\{|y-\ell_j|,\,j\in \mathcal{J}_1\}\ge R+\sigma.
\end{equation*}
(\ref{eq:rescaledcloseness}) implies that
\begin{equation*}
\|(\widetilde{u}_{0n},\,\widetilde{u}_{1n})-(\widetilde{h}_{0n},\,\widetilde{h}_{1n})\|_{\HL(B_{R+\frac{\sigma}{2}}(y))}+\|\widetilde{u}_{0n}-\widetilde{h}_{0n}\|_{L^{\dual}(B_{R+\frac{\sigma}{2}}(y))}\to 0,\,\,{\rm as}\,\,n\to\infty.
\end{equation*}
Hence by Lemma \ref{lm:localinspaceapproximation}, we obtain that
\begin{eqnarray}
&&\sup_{|t|\leq T_1\wedge R}\,\|\OR{\widetilde{h}}_n(t)-\OR{\widetilde{u}}_n(t)\|_{\HL(|x-y|\leq R-|t|)}\nonumber\\
&&\nonumber\\
&&\quad\quad\quad\quad\quad+\,\|\widetilde{h}_n-\widetilde{u}_n\|_{\Snorm\left(\left\{(x,t):\,|x|\leq R-|t|,\,|t|\leq T_1\wedge R\right\}\right)}\to 0.\label{eq:localcloseness}
\end{eqnarray}

On the other hand, as pointed out earlier, by Lemma \ref{lm:nonlinearprofiledecomposition}, the solution $\OR{\widetilde{h}}_n$ admits the following expansion, for $|t|\leq T_1\wedge R$:
\begin{eqnarray*}
\OR{\tilde{h}}_n(x,t)&=&\OR{U}(x,t)+\left(t_n^{\frac{d}{2}-1}\,w^J_{n}(t_nx,\,t_nt),\,t_n^{\frac{d}{2}}\,\partial_tw^J_{n}(t_nx,\,t_nt)\right)+\OR{r}^J_n(x,t)\\
&&\quad\quad\quad\quad+\sum_{j\not\in \mathcal{J}_1,\,1<j\leq J}\left(\left(t_n^{-1}\lambda^j_n\right)^{-\frac{d}{2}+1}U_j\left(\frac{x-t_n^{-1}c^j_n}{t_n^{-1}\lambda^j_n}, \,\frac{t-t_n^{-1}t^j_n}{t_n^{-1}\lambda^j_n}\right),\right.\\
&&\quad\quad\quad\quad\quad\quad\quad\quad\quad\left.\left(t_n^{-1}\lambda^j_n\right)^{-\frac{d}{2}}\partial_tU_j\left(\frac{x-t_n^{-1}c^j_n}{t_n^{-1}\lambda^j_n}, \,\frac{t-t_n^{-1}t^j_n}{t_n^{-1}\lambda^j_n}\right)\right),
\end{eqnarray*}
where $\lim\limits_{J\to\infty}\limsup\limits_{n\to\infty}\,\sup_{|t|\leq T_1\wedge R}\,\|\OR{r}^J_n(t)\|_{\HL}=0$.
By the pseudo-orthogonality of parameters (\ref{eq:profileorthogonality4}) and the fact that $t_n=\lambda^1_n$, the approximation (\ref{eq:localcloseness}) and the weak convergence to $0$ of 
$\left(t_n^{\frac{d}{2}-1}w^J_n(t_n\,x,\,t_n\,t),\,t_n^{\frac{d}{2}}\partial_tw^J_n(t_n\,x,\,t_n\,t)\right)$,
which is a consequence of the profile decomposition and $\lambda^1_n=t_n$,  we obtain that
\begin{equation}\label{11111}
 \OR{\widetilde{u}}_n\rightharpoonup \OR{U}\,\,\,\,{\rm in}\,\, \{(x,t):\,|x-y|< R-|t|,\,|t|< T_1\wedge R\},
\end{equation} 
weakly in the $(x,t)$ variables.\\

By finite speed of propagation and rescaling, we have 
\begin{equation}
\widetilde{u}_n(x,t)=t_n^{\frac{d}{2}-1}u(t_nx,\,t_n+t_nt),
\end{equation}
for $|x|<2$ and $|t|<1$. After rescaling, (\ref{11111}) and (\ref{eq:vanishingcharacterization}) with $\tau=\frac{t_n}{20}$ imply that
\begin{equation*}
\int_{\left\{(x,t):\,|x|<2,\,|x-y|\leq R-|t|,\,|t|\leq T_1\wedge R\right\}}\left(\partial_tU+\frac{x}{t+1}\cdot\nabla U+\left(\frac{d}{2}-1\right) \frac{U}{t+1}\right)^2\,dxdt=0.\\
\end{equation*}
Hence 
\begin{equation*}
\partial_tU+\frac{x}{t+1}\cdot\nabla U+\left(\frac{d}{2}-1\right) \frac{U}{t+1}\equiv 0,
\end{equation*}
in $\{(x,t):\,|x|<2,\,|x-y|\leq R-|t|,\,|t|\leq T_1\wedge R\}.$\\

Moving $y$ around and shrinking $\sigma$ to $0$, we can conclude that in fact
\begin{equation}\label{eq:firstorderselfsimilar}
\partial_tU+\frac{x}{t+1}\cdot\nabla U+\left(\frac{d}{2}-1\right) \frac{U}{t+1}\equiv 0,
\end{equation}
in $\left\{(x,t):\,|t|<{\rm dist}\left(x,\,\{\ell_j,\,j\in \mathcal{J}_1\}\right),\,\,|x|<2,\,\,|t|<T_1\right\}$. It is easy to see from (\ref{eq:firstorderselfsimilar}), by integrating along the characteristics, that for some $\Psi$ we have
\begin{equation}\label{eq:selfsimilarrepresentation}
U(x,t)=(t+1)^{-\frac{d}{2}+1}\Psi\left(\frac{x}{t+1}\right),
\end{equation}
in $\left\{(x,t):\,|t|<\frac{1}{4}\,{\rm dist}\left(x,\,\{\ell_j,\,j\in \mathcal{J}_1\}\right),\,|t|<T_1,\,|x|<2\right\}.$ Indeed
$$\frac{\partial}{\partial t} \left( (t+1)^{\frac d2-1} U\left( (t+1)y,t \right)\right)=0,$$
in $\left\{(y,t):\,|t|<{\rm dist}\left( (t+1)y,\,\{\ell_j,\,j\in \mathcal{J}_1\}\right),\,\,|(t+1)y|<2,\,\,|t|<T_1\right\}$ and the results follow.
\\

 Since ${\rm supp}\,\OR{U}(0)\subseteq \overline{B_1}$, we conclude that 
\begin{equation}\label{eq:equalinitialdata}
\OR{U}(0)=\left(\Psi,\,-x\cdot\nabla \Psi-\left(\frac{d}{2}-1\right)\Psi\right)\,\,\,{\rm in}\,\,R^d.
\end{equation}
Since $\OR{U}(0)$ is supported in $\overline{B_1}$ and satisfies equation (\ref{eq:main}) in $R^d\times [-T,\,T]$, we conclude that ${\rm supp}\,\Psi\subseteq \overline{B_1}$ and
\begin{equation}\label{eq:degenerateelliptic1}
-\Delta \Psi+y\cdot\nabla (y\cdot\nabla \Psi)+(d-1)y\cdot\nabla \Psi+\frac{d(d-2)}{4}\Psi=|\Psi|^{\dual-2}\Psi,
\end{equation} 
in $R^d\backslash\{\ell_j,\,j\in \mathcal{J}_1\}$. By (\ref{eq:equalinitialdata}) and the fact that $\OR{U}(0)\in\HL$, we conclude that $\Psi\in \dot{H}^1$, and hence the equation (\ref{eq:degenerateelliptic1}) is in fact satisfied ``across" the points $\ell_j$, i.e., $\ell_j$ are removable singularities. Thus equation (\ref{eq:degenerateelliptic1}) holds in $R^d$. Standard elliptic regularity theory shows that $\Psi\in C^{2}(B_1)$. Hence  
$$\widetilde{U}(x,t):=(t+1)^{-\frac{d}{2}+1}\Psi\left(\frac{x}{t+1}\right)$$
is a classical solution to equation (\ref{eq:main}) for $|x|<1-|t|$, $|t|<T$. By (\ref{eq:equalinitialdata}), we see that $\OR{U}(0)=\OR{\widetilde{U}}(0)$. Hence, by finite speed of propagation, $\widetilde{U}\equiv U$ for $|x|<1-|t|$, $|t|<T$. Thus, $\widetilde{U}\in \Snorm\left(|x|<(1-|t|),\,|t|<T\right)$. By the support property of $\Psi$, thus $\widetilde{U}\in\Snorm(R^d\times [-a,a])$ with $a=\min(\frac{1}{2},\,T)$, and is a solution to equation (\ref{eq:main}) with the same initial data as that of $U$. Therefore, $U\equiv \widetilde{U}$ is a compactly supported finite energy exact self similar solution, as in Proposition 5.7 of \cite{KMacta}. By \cite[section 6]{KMacta} we see that $U$ must be trivial. \\

\smallskip
\noindent
{\it Case II:} Now we consider Case II, for a fixed $j$, $t^j_n\equiv 0$, $\lambda^j_n\ll t_n$. To simplify notation, we assume that the profile is $U_1^L$ with parameters $\lambda^1_n$, $c^1_n$ satisfying $\lim\limits_{n\to\infty}\frac{c^1_n}{t_n}=\ell_1$ and $t^1_n\equiv 0$. 
Suppose that the nonlinear profile $U_1$ associated with $U_1^L$, $\lambda^1_n$ exists in $R^d\times [-T,\,T]$ with 
\begin{equation*}
\|U_1\|_{\Snorm(R^d\times[-T,T])}\leq M<\infty.
\end{equation*}

The idea in the characterization of $U_1$ is similar to Case I. We shall still remove the profiles which contain more energy than the threshold energy provided by small data theory and have nontrivial interactions with $U_1$ (there are only finitely many such profiles), and then use perturbative arguments to deal with other small profiles. After a rescaling, we can again pass to the limit in a region with several lightcones removed. Using the vanishing condition (\ref{eq:vanishingcharacterization}), we then obtain a first order equation for $U_1$ which will enable us to classify $U_1$. The difference with Case I is due to the fact that here $\lambda^1_n\ll t_n$. As a consequence, the first order equation we obtain in the end is different from the self similar case, and $U_1$ has to be a traveling wave, instead of a self similar solution. \\

Once more, fix $\epsilon=\epsilon(d)>0$ sufficiently small, and let $J(\epsilon)$ be such that (\ref{eq:allsmallprofiles}) holds. The profiles $U^L_j$, $j\geq J(\epsilon)$, will be controlled perturbatively. For $1<j<J(\epsilon)$, we split the profiles into five categories.
\begin{itemize}
\item i) The $j$'s for which $\lambda^j_n\ll \left|t^j_n\right|$, and which, passing to a subsequence if necessary, $\left|t^j_n\right|\gtrsim \lambda^1_n$ for all $n$. These profiles have negligible interaction with the profile $U_1^L$, with scaling parameter $\lambda^1_n$, as can be easily seen by rescaling, at least for a short time comparable to $\lambda^1_n$. More precisely, we have, as can be easily seen by rescaling,
\begin{equation}\label{eq:largeprofileok}
\lim_{n\to\infty}\,\left\|(\lambda^j_n)^{-\frac{d}{2}+1}U^L_j\left(\frac{x-c^j_n}{\lambda^j_n},\,\frac{t-t^j_n}{\lambda^j_n}\right)\right\|_{\Snorm\left(R^d\times [-T_1\lambda^1_n,\,T_1\lambda^1_n]\right)}\to 0,
\end{equation}
as $n\to\infty$, if $T_1>0$ is sufficiently small and hence, as we will see, these profiles can be controlled perturbatively. 
\item ii) The profiles $U^L_j$ such that $\lambda^j_n\ll \left|t^j_n\right|$, $t^j_n=o(\lambda^1_n)$ as $n\to\infty$. We denote this set of $j$ as $\mathcal{J}_1$. These profiles have nontrivial interaction with the profile $U_1^L$, $\lambda^1_n$, and thus will be removed.
\item iii) The profiles $U^L_j$ that have the property that $t^j_n\equiv 0$ and passing to a subsequence if necessary, that $\lambda^j_n\ll \lambda^1_n$. These profiles have nontrivial interaction with the profile $U_1^L$, $\lambda^1_n$, and thus will be removed. We denote this set of $j$ as $\mathcal{J}_2$. 
\item iv) The $j$'s for which $t^j_n\equiv 0$, and after passing to a subsequence if necessary, $\lambda^j_n\sim \lambda^1_n$. Then (by the pseudo-orthogonality of parameters), $\frac{\left|c^j_n-c^1_n\right|}{\lambda^1_n}\to\infty$. These profiles have no interaction in the limit with the first profile, locally in space, and will be dealt with perturbatively.
\item v)  The $j$'s for which $t^j_n\equiv 0$, $\lambda^1_n=o(\lambda^j_n)$. These profiles have asymptotically vanishing interaction with the first profile for a time interval of the size of $ \lambda^1_n$, and will be dealt with perturbatively.
\end{itemize}

Define now, for $J\in\mathbb{N}$,
\begin{eqnarray*}
&&\hspace{1.3in}(h_{0n},\,h_{1n})=\\
\\
&&\left((\lambda^1_n)^{-\frac{d}{2}+1}U^L_1\left(\frac{x-c^1_n}{\lambda^1_n},\,0\right),\,(\lambda^1_n)^{-\frac{d}{2}}\partial_tU_1^L\left(\frac{x-c^1_n}{\lambda^1_n},\,0\right)\right)\\
\\
&&\quad\quad+\sum_{1<j\leq J,\,\,\,j\not\in \mathcal{J}_1\cup \mathcal{J}_2}\bigg(\left(\lambda^j_n\right)^{-\frac{d}{2}+1}U_j^L\left(\frac{x-c^j_n}{\lambda^j_n}, \,-\frac{t^j_n}{\lambda^j_n}\right),\\
&&\quad\quad\quad\hspace{4cm}\,\,\,\,\,\,\,\,\,\,\left(\lambda^j_n\right)^{-\frac{d}{2}}\partial_tU_j^L\left(\frac{x-c^j_n}{\lambda^j_n}, \,-\frac{t^j_n}{\lambda^j_n}\right)\bigg)\\
&&\\
&&\quad\quad\quad\quad+\,(w^J_{0n},\,w^J_{1n}).
\end{eqnarray*}

It is simple to check that
\begin{eqnarray}
&&\|(h_{0n},\,h_{1n})-(u_{0n},\,u_{1n})\|_{\HL\left(B_{M\lambda^1_n}(c^1_n)\backslash\, \bigcup_{j\in \mathcal{J}_1\cup \mathcal{J}_2} B_{\eta\, \lambda^1_n}(c^j_n)\right)}\nonumber\\
&&\quad\quad\quad\quad+\,\|h_{0n}-u_{0n}\|_{L^{\dual}\left(B_{M\lambda^1_n}(c^1_n)\backslash\, \bigcup_{j\in \mathcal{J}_1\cup \mathcal{J}_2} B_{\eta\, \lambda^1_n}(c^j_n)\right)}\,\to\, 0,\label{eq:errorvanishesinitial}
\end{eqnarray}
as $n\to\infty$, for any $M>1$ and $\eta>0$. \\

Consider the rescaled initial data sequences
\begin{eqnarray*}
&&(\widetilde{u}_{0n},\,\widetilde{u}_{1n})=\left(t_n^{\frac{d}{2}-1}u_{0n}(t_n\,x),\,t_n^{\frac{d}{2}}u_{1n}(t_n\,x)\right);\\
&&(\widetilde{h}_{0n},\,\widetilde{h}_{1n})=\left(t_n^{\frac{d}{2}-1}h_{0n}(t_n\,x),\,t_n^{\frac{d}{2}}h_{1n}(t_n\,x)\right),
\end{eqnarray*}
and let $\widetilde{u}_n$, $\widetilde{h}_n$ be the solution to equation (\ref{eq:main}) with initial data $(\widetilde{u}_{0n},\,\widetilde{u}_{1n})$ and $(\widetilde{h}_{0n},\,\widetilde{h}_{1n})$, respectively. By the principle of finite speed of propagation and rescaling, we get that
\begin{equation}
\widetilde{u}_n(x,t)=t_n^{\frac{d}{2}-1}\,u(t_nx,\,t_n(t+1)),\,\,{\rm for}\,\,|x|<2,\,\,|t|<\frac{1}{2}.
\end{equation}

By Lemma \ref{lm:nonlinearprofiledecomposition}, $\widetilde{h}_n$ admits the following decomposition for $|t|<T_1\lambda^1_nt_n^{-1}$, $T_1$ small
\begin{eqnarray*}
\widetilde{h}_n(x,t)&=&\bigg((t_n^{-1}\lambda^1_n)^{-\frac{d}{2}+1}U_1\left(\frac{x-t_n^{-1}c^1_n}{t_n^{-1}\lambda^1_n},\,\frac{t}{t_n^{-1}\lambda^1_n}\right),\\
&&\hspace{4cm}\,(t_n^{-1}\lambda^1_n)^{-\frac{d}{2}}\partial_tU_1\left(\frac{x-t_n^{-1}c^1_n}{t_n^{-1}\lambda^1_n},\,\frac{t}{t_n^{-1}\lambda^1_n}\right)\bigg)\\
\\
&&\quad\quad\quad+\sum_{1<j\leq J,\,\,j\not\in \mathcal{J}_1\cup \mathcal{J}_2}\left(\left(t_n^{-1}\lambda^j_n\right)^{-\frac{d}{2}+1}U_j\left(\frac{x-t_n^{-1}c^j_n}{t_n^{-1}\lambda^j_n}, \,\frac{t-t_n^{-1}t^j_n}{t_n^{-1}\lambda^j_n}\right),
\right.\\
&&\quad\quad\quad\quad\quad\quad\quad\quad\quad\quad\quad\quad\,\left.\left(t_n^{-1}\lambda^j_n\right)^{-\frac{d}{2}}\partial_tU_j\left(\frac{x-t_n^{-1}c^j_n}{t_n^{-1}\lambda^j_n}, \,\frac{t-t_n^{-1}t^j_n}{t_n^{-1}\lambda^j_n}\right)\right)\\
&&\\
&&\quad\quad\quad\quad\quad+\left(t_n^{\frac{d}{2}-1}w_n^J(t_nx,t_nt),\,t_n^{\frac{d}{2}}\partial_tw^J_n(t_nx,t_nt)\right)+\OR{r}^J_n,
\end{eqnarray*}
where 
\begin{equation}
\lim_{J\to\infty}\limsup_{n\to\infty}\,\sup_{|t|<T_1\lambda^1_nt_n^{-1}}\,\left\|\OR{r}^J_n(t)\right\|_{\HL}=0;
\end{equation}
and $\widetilde{h}_n$ verifies the bound
\begin{equation}
\left\|\widetilde{h}_n\right\|_{\Snorm\left(R^d\times[-T_1t_n^{-1}\lambda^1_n,\,T_1t_n^{-1}\lambda^1_n]\right)}\leq M_1<\infty,
\end{equation}
for $J,\,n$ sufficiently large.\\

Consider the rescaled and translated $\widetilde{u}_n$ and $\widetilde{h}_n$ as follows:
\begin{eqnarray*}
\OR{\widetilde{\widetilde{u}}}_n(t)&:=&\left((t_n^{-1}\lambda^1_n)^{\frac{d}{2}-1}\,\widetilde{u}_n\left(t_n^{-1}\lambda^1_n\,(x+t_n^{-1}c^1_n),\,t_n^{-1}\lambda^1_n\,t\right),\right.\\
&&\quad\quad\quad\quad\quad\,\left.(t_n^{-1}\lambda^1_n)^{\frac{d}{2}}\,\partial_t\widetilde{u}_n\left(t_n^{-1}\lambda^1_n\,(x+t_n^{-1}c^1_n),\,t_n^{-1}\lambda^1_n\,t\right)\right),\\
\OR{\widetilde{\widetilde{h}}}_n(t)&:=&\left((t_n^{-1}\lambda^1_n)^{\frac{d}{2}-1}\,\widetilde{h}_n\left(t_n^{-1}\lambda^1_n\,(x+t_n^{-1}c^1_n),\,t_n^{-1}\lambda^1_n\,t\right),\right.\\
&&\quad\quad\quad\quad\quad\,\left.(t_n^{-1}\lambda^1_n)^{\frac{d}{2}}\,\partial_t\widetilde{h}_n\left(t_n^{-1}\lambda^1_n\,(x+t_n^{-1}c^1_n),\,t_n^{-1}\lambda^1_n\,t\right)\right),
\end{eqnarray*}
with initial data $(\widetilde{\widetilde{u}}_{0n},\widetilde{\widetilde{u}}_{1n})$ and $(\widetilde{\widetilde{h}}_{0n},\widetilde{\widetilde{h}}_{1n}).$

By rescaling and translation, (\ref{eq:errorvanishesinitial}) implies that for $M>1$, $\eta>0$, we have
\begin{eqnarray}
&&\left\|\left(\widetilde{\widetilde{h}}_{0n},\,\widetilde{\widetilde{h}}_{1n}\right)-\left(\widetilde{\widetilde{u}}_{0n},\,\widetilde{\widetilde{u}}_{1n}\right)\right\|_{\HL\left(B_M(0)\backslash\,\bigcup_{j\in \mathcal{J}_1\cup \mathcal{J}_2,\,j>1}B_{\eta}(\frac{c^j_n-c^1_n}{\lambda^1_n})\right)}\nonumber\\
&&\quad\quad\quad\quad+\,\left\|\widetilde{\widetilde{h}}_{0n}-\widetilde{\widetilde{u}}_{0n}\right\|_{L^{\dual}\left(B_M(0)\backslash\,\bigcup_{j\in \mathcal{J}_1\cup \mathcal{J}_2,\,j>1}B_{\eta}(\frac{c^j_n-c^1_n}{\lambda^1_n})\right)}\,\rightarrow\, 0,\label{eq:errorvanishesinitial2}
\end{eqnarray}
as $n\to\infty$.\\

Denote 
\begin{equation}
\mathcal{J}_{\ast}:=\left\{j\in \mathcal{J}_1\cup \mathcal{J}_2:\,\,\frac{c^j_n-c^1_n}{\lambda^1_n}\,\,{\rm is\,\,bounded\,\,in\,\,}n\right\}.
\end{equation}
Passing to a subsequence, we can assume that $\frac{c^j_n-c^1_n}{\lambda^1_n}\to x_j$ as $n\to\infty$ for each $j\in \mathcal{J}_{\ast}$. For any $R>0$, $\sigma>0$, $y\in R^d$ with $|y|\leq M-R$ and
\begin{equation*}
{\rm dist}\left(y,\,\{x_j,\,j\in \mathcal{J}_{\ast}\}\right)>R+\sigma,
\end{equation*}
we see by Lemma \ref{lm:localinspaceapproximation} and (\ref{eq:errorvanishesinitial2}) that
\begin{eqnarray}
&&\sup_{|t|<T_1\wedge R}\,\left\|\OR{\widetilde{\widetilde{u}}}_n(t)-\OR{\widetilde{\widetilde{h}}}_n(t)\right\|_{\HL(\{|x-y|<R-|t|\})}\nonumber\\
&&\hspace{1.5cm}+\,\left\|\widetilde{\widetilde{u}}_n-\widetilde{\widetilde{h}}_n\right\|_{\Snorm\left(\left\{(x,t):\,|x-y|<R-|t|,\,|t|<T_1\wedge R\right\}\right)}\rightarrow 0,\label{eq:vanishesimportant}
\end{eqnarray}
as $n\to\infty$.\\

The vanishing condition (\ref{eq:vanishingcharacterization}) at $\tau=\lambda^1_n$, with a simple rescaling argument, implies that, for any $M>1$
\begin{eqnarray*}
&&\frac{1}{t_n^{-1}\lambda^1_n}\int_{-t_n^{-1}\lambda^1_n}^{t_n^{-1}\lambda^1_n}\int_{\left|x-\frac{c^1_n}{t_n}\right|<Mt_n^{-1}\lambda^1_n}\,\left|\partial_t\widetilde{u}_n+\frac{x}{t+1}\cdot\nabla \widetilde{u}_n+\left(\frac{d}{2}-1\right)\frac{\widetilde{u}_n}{t+1}\right|^2\,dxdt\\
\\
&&\to 0, \,\,{\rm as}\,\,n\to\infty.
\end{eqnarray*}

Note that, by H\"older's inequality,
\begin{eqnarray*}
&&\frac{1}{t_n^{-1}\lambda^1_n}\int_{-t_n^{-1}\lambda^1_n}^{t_n^{-1}\lambda^1_n}\int_{\left|x-\frac{c^1_n}{t_n}\right|<Mt_n^{-1}\lambda^1_n}\,\left|\frac{\widetilde{u}_n}{t+1}\right|^2\,dxdt\\
&&\quad\lesssim\, \left[\sup_{t\in (-t_n^{-1}\lambda^1_n,\,t_n^{-1}\lambda^1_n)}\|\widetilde{u}_n(\cdot,t)\|^2_{L^{\dual}(R^d)}\right]\cdot\left(t_n^{-1}\lambda^1_n\right)^{2}\to 0,
\end{eqnarray*}
as $n\to\infty$. Noting also that 
\begin{equation*}
\lim_{n\to\infty}\,\frac{c^1_n}{t_n}=\ell_1,\,\,\,{\rm and}\,\,\lambda^1_n\ll t_n,
\end{equation*}
we get that 
\begin{equation}\label{eq:vanishingcharacterization5}
\lim_{n\to\infty}\,\frac{1}{t_n^{-1}\lambda^1_n}\int_{-t_n^{-1}\lambda^1_n}^{t_n^{-1}\lambda^1_n}\int_{\left|x-\frac{c^1_n}{t_n}\right|<Mt_n^{-1}\lambda^1_n}\,\left|\partial_t\widetilde{u}_n+\ell_1\cdot\nabla \widetilde{u}_n\right|^2\,dxdt=0.
\end{equation}
Rescaling and translating (\ref{eq:vanishingcharacterization5}), we obtain that 
\begin{equation}
\lim_{n\to\infty}\int_{-1}^1\int_{|x|<M}\left|\partial_t\widetilde{\widetilde{u}}_n+\ell_1\cdot\nabla \widetilde{\widetilde{u}}_n\right|^2\,dxdt=0,
\end{equation}
for each $M>1$.

By the expansion of $\widetilde{h}_n$, $\widetilde{\widetilde{h}}_n$ admits the following decomposition:
\begin{eqnarray*}
&&\widetilde{\widetilde{h}}_n=\OR{U}_1(x,t)+\OR{\widetilde{\widetilde{r}}}^J_n+\\
&&\quad\sum_{1< j\leq J,\,j\not\in \mathcal{J}_1\cup \mathcal{J}_2}\left(\left((\lambda^1_n)^{-1}\lambda^j_n\right)^{-\frac{d}{2}+1}U_j\left(\frac{x-\frac{c^j_n-c^1_n}{\lambda^1_n}}{(\lambda^1_n)^{-1}\lambda^j_n},\,\frac{t-(\lambda^1_n)^{-1}t^j_n}{(\lambda^1_n)^{-1}\lambda^j_n}\right),\right.\\
&&\quad\quad\quad\quad\quad\quad\quad\quad\quad\quad\left.\left((\lambda^1_n)^{-1}\lambda^j_n\right)^{-\frac{d}{2}}\partial_tU_j\left(\frac{x-\frac{c^j_n-c^1_n}{\lambda^1_n}}{(\lambda^1_n)^{-1}\lambda^j_n},\,\frac{t-(\lambda^1_n)^{-1}t^j_n}{(\lambda^1_n)^{-1}\lambda^j_n}\right)
\right)\\
&&\quad\quad\quad\quad\quad +\,\left((\lambda^1_n)^{-\frac{d}{2}+1}w^J_n\left(\frac{x+\frac{c^1_n}{\lambda^1_n}}{\lambda^1_n},\,\frac{t}{\lambda^1_n}\right),\,(\lambda^1_n)^{-\frac{d}{2}}\partial_tw^J_n\left(\frac{x+\frac{c^1_n}{\lambda^1_n}}{\lambda^1_n},\,\frac{t}{\lambda^1_n}\right)\right).
\end{eqnarray*}

By the pseudo-orthogonality of the parameters and the weak convergence to $0$ property of $w^J_n$, due to the profile decomposition, we conclude that 
\begin{equation}
\OR{\widetilde{\widetilde{h}}}_n\rightharpoonup \OR{U}_1 \,\,\,\,{\rm in}\,\,\left\{(x,t):\,|x-y|<R\wedge T_1-|t|,\,|t|<R\wedge T_1\right\}.
\end{equation}
By (\ref{eq:vanishesimportant}), we also have
\begin{equation}
\OR{\widetilde{\widetilde{u}}}_n\rightharpoonup \OR{U}_1 \,\,\,\,{\rm in}\,\,\left\{(x,t):\,|x-y|<R\wedge T_1-|t|,\,|t|<R\wedge T_1\right\}.
\end{equation}
 Hence
\begin{equation}\label{eq:firstorder222}
\partial_tU_1+\ell_1\cdot\nabla U_1\equiv 0,
\end{equation}
in $\left\{(x,t):\,|x-y|<R\wedge T_1-|t|,\,|t|<R\wedge T_1\right\}$. \\

Moving $y$ around and shrinking $\sigma$ to $0$, we conclude that
\begin{equation}\label{eq:firstordersoliton}
\partial_tU_1+\ell_1\cdot\nabla U_1\equiv 0,
\end{equation}
in $\left\{(x,t):\,|t|<{\rm dist}\left(x,\,\{x_j,\,j\in \mathcal{J}_{\ast}\}\right)\wedge T_1\right\}$. 
(\ref{eq:firstordersoliton}) implies that $U_1$ satisfies
\begin{equation*}
U_1(x,t)=\Psi(x-\ell_1t),
\end{equation*}
in $\left\{(x,t):\,|t|<{\rm dist}\left(x,\,\{x_j,\,j\in \mathcal{J}_{\ast}\}\right)\wedge T_1\right\}$. Since $U_1$ solves equation (\ref{eq:main}), we conclude that
\begin{equation}\label{eq:ellipticsoliton}
\ell_1\cdot\nabla\left(\ell_1\cdot\nabla \Psi\right)-\Delta\Psi=|\Psi|^{\dual-2}\Psi,
\end{equation}
in $R^d\backslash\,\{x_j,\,j\in \mathcal{J}_{\ast}\}$. Since $\Psi\in\dot{H}^1(R^d)$, equation (\ref{eq:ellipticsoliton}) is actually satisfied across $x_j$, i.e., $x_j$ are removable singularities. Hence $\Psi\in\dot{H}^1(R^d)$ and
\begin{equation}\label{eq:ellipticsoliton1}
\ell_1\cdot\nabla\left(\ell_1\cdot\nabla \Psi\right)-\Delta\Psi=|\Psi|^{\dual-2}\Psi,\,\,{\rm in}\,\,R^d.
\end{equation}
Using finite speed of propagation, one can prove $|\ell_1|\leq 1$. It follows from a Pohozaev identity that $|\ell_1|=1$ is also impossible. We refer to Steps 1 and 2 of the proof of \cite[Lemma 2.6]{DKMsmallnonradial} for the details. Finally, using the change of variable:
$$Q(y)=\Psi\left( y+\frac{y\cdot\ell_1}{|\ell_1|^2}\ell_1\left(\frac{1}{\sqrt{1-|\ell_1|^2}}-1\right) \right),$$
we see that $-\Delta Q=|Q|^{\dual-2}Q$, and that
$\Psi\equiv Q_{\ell_1}$. Hence, in this case the profile is a traveling wave, which is non-trivial since $U_1$ is non-trivial. By the energy expansion, there can only be only a finite number of such profiles. This finishes dealing with the profiles in Case II.\\

For profiles in Case III, $0<\lambda^j_n\ll |t^j_n|$, and hence
\begin{eqnarray*}
&&\lim_{n\to\infty}\left\|(\lambda^j_n)^{-\frac{d}{2}+1}U_j^L\left(\frac{x-c^j_n}{\lambda^j_n},\,-\frac{t^j_n}{\lambda^j_n}\right)\right\|_{L^{\dual}(R^d)}\\
&&=\lim_{n\to\infty}\left\|U_j^L\left(x,\,-\frac{t^j_n}{\lambda^j_n}\right)\right\|_{L^{\dual}(R^d)}\\
&&=\lim_{t\to\pm\infty}\left\|U_j^L(x,t)\right\|_{L^{\dual}(R^d)}\\
&&=0,
\end{eqnarray*}
by a well-known property of linear waves. Thus, these profiles can be absorbed into the error term. Recalling that $v$ is the regular part of $u$ and hence
$$(u_{0n},u_{1n})=\OR{u}(t_n)-\OR{v}(t_n)+o_n(1),$$
$\,\,{\rm in}\,\,\big(\dot{H}^1\cap L^{\dual}\big)\times L^2 (\{|x|<r_0\})$, so that the theorem is proved.\\

\section{virial identity and exclusion of dispersive energy in the region $\{(x,t):\,|x|<\lambda t\}$ for $\lambda\in (0,1)$ along a sequence of times}\label{sec:virial}
In this section, we use a virial identity to obtain a decomposition with better residue term. 
\begin{theorem}\label{th:betterresidue}
Let $u$ be as in the last section. Then there exists a time sequence $t_n\downarrow 0$, such that $\OR{u}(t_n)$ has the asymptotic decomposition (\ref{eq:maindecompositionsec5}) with properties as in Theorem \ref{th:mainpreliminary}. In addition the residue term $(\epsilon_{0n},\,\epsilon_{1n})$ satisfies the following refined vanishing condition:
\begin{eqnarray*}
&&\|\spartial\epsilon_{0n}\|_{L^2(B_{r_{0}})}+\|\nabla \epsilon_{0n}\|_{L^2(B_{\lambda t_n}\cup\, B^c_{t_n})}+\|\epsilon_{1n}\|_{L^2(B_{\lambda t_n}\cup \,B^c_{t_n})}\nonumber\\
&&\nonumber\\
&&\hspace{2cm}+\,\left\|\epsilon_{1n}+\frac{r}{t_n}\partial_r\epsilon_{0n}\right\|_{L^2(B_{t_n})}\to 0,\label{eq:betterresidue}
\end{eqnarray*}
as $n\to\infty$, for any $\lambda\in(0,1)$, where $B^c_{t_n}$ is the complement set of $B_{t_n}$ in $B_{r_0}$.
\end{theorem}

\smallskip
\noindent
{\it Proof.} By Lemma \ref{lm:keyvanishingMorawetz}, there exists a sequence $\mu_n\downarrow 0$, such that
\begin{equation}\label{eq:vanishingcondition6}
\lim_{n\to\infty}\,\frac{1}{\mu_n}\int_{\mu_n}^{2\mu_n}\int_{|x|<4t}\,\left(\partial_tu+\frac{x}{t}\cdot\nabla u+\left(\frac{d}{2}-1\right)\frac{u}{t}\right)^2\,dx\,dt=0,
\end{equation}
and such that along two sequences of times $t_{1n}$ and $t_{2n}$ with $t_{1n}\in (\frac{4}{3}\mu_{n},\,\frac{13}{9}\mu_n)$, $t_{2n}\in (\frac{14}{9}\mu_n,\,\frac{5}{3}\mu_n)$, we have that 
\begin{equation}\label{eq:sequencevanishing1}
\sup_{0<\tau<\frac{t_{1n}}{16}}\,\frac{1}{\tau}\int_{|t-t_{1n}|<\tau}\int_{|x|<4t}\,\left(\partial_tu+\frac{x}{t}\cdot\nabla u+\left(\frac{d}{2}-1\right)\frac{u}{t}\right)^2\,dxdt\to 0,\,\,{\rm as}\,\,n\to\infty,
\end{equation}
and
\begin{equation}\label{eq:sequencevanishing2}
\sup_{0<\tau<\frac{t_{2n}}{16}}\,\frac{1}{\tau}\int_{|t-t_{2n}|<\tau}\int_{|x|<4t}\,\left(\partial_tu+\frac{x}{t}\cdot\nabla u+\left(\frac{d}{2}-1\right)\frac{u}{t}\right)^2\,dxdt\to 0,\,\,{\rm as}\,\,n\to\infty.
\end{equation}
By Theorem \ref{th:mainpreliminary}, along $t_{1n}$ and $t_{2n}$ approaching zero, we have the following decomposition
\begin{eqnarray}
&&\hspace{3cm}\OR{u}(t_{\iota n})=\OR{v}(t_{\iota n})+\OR{\epsilon}_{\iota n}+\nonumber\\
&&\sum_{j=1}^{J_{\iota}}\,\left((\lambda^j_{\iota n})^{-\frac{d}{2}+1}\, Q_{\ell_{\iota j}}^{\iota j}\left(\frac{x-c^j_{\iota n}}{\lambda^j_{\iota n}},\,0\right),\,(\lambda^j_{\iota n})^{-\frac{d}{2}}\, \partial_tQ_{\ell_{\iota j}}^{\iota j}\left(\frac{x-c^j_{\iota n}}{\lambda^j_{\iota n}},\,0\right)\right),\label{eq:maindecomposition5}
\end{eqnarray}
where $\iota\in\{1,\,2\}$, and $\lambda^j_{\iota n}\ll t_{\iota n}$, $c^j_{\iota n}\in B_{\beta t_{\iota n}}$ for some $\beta\in (0,1)$ and $\ell_{\iota j}=\lim\limits_{n\to\infty}\frac{c^j_{\iota n}}{t_{\iota n}}$. $\lambda^j_{\iota n},\,c^j_{\iota n}$ satisfy the pseudo-orthogonality condition (\ref{eq:pseudo}). In the above, with slight abuse of notation, we denote $\OR{\epsilon}_{\iota n}:=(\epsilon_{\iota n},\,\partial_t\epsilon_{\iota n})\in\HL$. Moreover, $\epsilon_{\iota n}$ vanishes in the sense that for some fixed $r_0>0$
\begin{equation}\label{eq:vanishingdispersivetermLp}
\|\epsilon_{\iota n}\|_{L^{\dual}(|x|\leq r_0)}\to 0,\quad {\rm as}\,\,n\to\infty.
\end{equation}
We observe that by the decomposition (\ref{eq:maindecomposition5}), for any $\epsilon>0$ small, it holds that
\begin{eqnarray*}
&&\|u(\cdot,t_{\iota n})\|_{L^2(B_{t_{\iota n}})}\\
\\
&&\lesssim \|v(\cdot,\,t_{\iota n})\|_{L^2(|x|<t_{\iota n})}+\|\epsilon_{\iota n}\|_{L^2(|x|<t_{\iota n})}\\
                                                                 &&\quad\quad+\,\left\|\sum_{j=1}^{J_{\iota}}\,(\lambda^j_{\iota n})^{-\frac{d}{2}+1}\, Q_{\ell_{\iota j}}^{\iota j}\left(\frac{x-c^j_{\iota n}}{\lambda^j_{\iota n}},\,0\right)\right\|_{L^2\left(\bigcup_{j=1}^{J_{\iota}}B_{\epsilon\, t_{\iota n}(c^j_{\iota n})}\right)}\\
&&\quad\quad+\,\left\|\sum_{j=1}^{J_{\iota}}\,(\lambda^j_{\iota n})^{-\frac{d}{2}+1}\, Q_{\ell_{\iota j}}^{\iota j}\left(\frac{x-c^j_{\iota n}}{\lambda^j_{\iota n}},\,0\right)\right\|_{L^2\left(B_{t_{\iota n}}\backslash \left(\bigcup_{j=1}^{J_{\iota}}B_{\epsilon\, t_{\iota n}}(c^j_{\iota n})\right)\right)}\\
&&\lesssim\|v(\cdot,t_{\iota n})\|_{L^{\dual}(B_{t_{\iota n}})}\,t_{\iota n}+\|\epsilon_{\iota n}\|_{L^{\dual}(B_{t_{\iota n}})}\,t_{\iota n}\\
&&\quad\quad+\left\|\sum_{j=1}^{J_{\iota}}\,(\lambda^j_{\iota n})^{-\frac{d}{2}+1}\, Q_{\ell_{\iota j}}^{\iota j}\left(\frac{x-c^j_{\iota n}}{\lambda^j_{\iota n}},\,0\right)\right\|_{L^{\dual}\left(\bigcup_{j=1}^{J_{\iota}}B_{\epsilon\, t_{\iota n}(c^j_{\iota n})}\right)}\,\epsilon\,t_{\iota n}\\
&&\quad\quad+\left\|\sum_{j=1}^{J_{\iota}}\,(\lambda^j_{\iota n})^{-\frac{d}{2}+1}\, Q_{\ell_{\iota j}}^{\iota j}\left(\frac{x-c^j_{\iota n}}{\lambda^j_{\iota n}},\,0\right)\right\|_{L^{\dual}\left(B_{t_{\iota n}}\backslash \left(\bigcup_{j=1}^{J_{\iota}}B_{\epsilon\, t_{\iota n}}(c^j_{\iota n})\right)\right)}\,t_{\iota n}\\
&&\lesssim o(t_{\iota n})+\epsilon \,t_{\iota n}, \quad {\rm as}\,\,n\to\infty.
\end{eqnarray*}
Since this is true for all small $\epsilon>0$, we obtain that
\begin{equation}\label{eq:secondvanishing}
\lim_{n\to\infty}\,\frac{1}{t_{\iota n}}\,\|u(\cdot,t_{\iota n})\|_{L^2(B_{t_{\iota n}})}=0.
\end{equation}

Thanks to the vanishing condition (\ref{eq:secondvanishing}), we can use the following virial identity.  Multiply equation (\ref{eq:main}) with $u$, and integrate over the region $\{(x,t):\,|x|<t,\,t\in (t_{1n},\,t_{2n})\}$. We obtain
\begin{eqnarray*}
0&=&\int_{t_{1n}}^{t_{2n}}\int_{|x|<t}\,\partial_t\left(\partial_tu\,u\right)-(\partial_tu)^2-{\rm div}\,\left(\nabla u \,u\right)+|\nabla u|^2-|u|^{\dual}\,dxdt\\
&=&\frac{1}{\sqrt{2}}\,\int_{t_{1n}}^{t_{2n}}\int_{|x|=t}\,\left(-\partial_tu-\frac{x}{|x|}\cdot\nabla u\right)\,u\,d\sigma dt\\
&&\quad\quad\quad+\int_{|x|<t_{2n}}\,[\partial_tu\,u](x,t_{2n})\,dx-\int_{|x|<t_{1n}}\,[\partial_tu\,u](x,t_{1n})\,dx\\
&&\quad\quad\quad+\int_{t_{1n}}^{t_{2n}}\int_{|x|<t}\,-(\partial_tu)^2+|\nabla u|^2-|u|^{\dual}\,dxdt.
\end{eqnarray*}
Noting that $t_{2n}\sim t_{1n}$ and $t_{2n}-t_{1n}\sim t_{1n}$, with the help of the control on the energy flux (Proposition \ref{L:trace}), we can estimate
\begin{eqnarray*}
&&\left|\int_{t_{1n}}^{t_{2n}}\int_{|x|=t}\,\left(-\partial_tu-\frac{x}{|x|}\cdot\nabla u\right)\,u\,d\sigma dt\right|\\
&&\lesssim\,\left(\int_{t_{1n}}^{t_{2n}}\int_{|x|=t}\,\left(\partial_tu+\frac{x}{|x|}\cdot\nabla u\right)^2\,d\sigma dt\right)^{\frac{1}{2}}\,\left(\int_{t_{1n}}^{t_{2n}}\int_{|x|=t}\,u^2\,d\sigma dt\right)^{\frac{1}{2}}\\
&&\lesssim o(1)\cdot \left(\int_{t_{1n}}^{t_{2n}}\int_{|x|=t}\,|u|^{\dual}\,d\sigma dt\right)^{\frac{1}{\dual}}\,|t_{2n}-t_{1n}|\\
&&\lesssim o\left(t_{2n}-t_{1n}\right),\quad{\rm as}\,\,n\to\infty.
\end{eqnarray*}

Using crucially (\ref{eq:secondvanishing}), we can also estimate
\begin{eqnarray*}
\left|\int_{|x|<t_{\iota n}}\,\partial_tu\,u(x,t_{\iota n})\,dx\right|&\lesssim&\,\|\partial_tu\|_{L^2(B_{t_{\iota n}})}\,\left(\int_{|x|<t_{\iota n}}\,u^2(x,t_{\iota n})\,dx\right)^{\frac{1}{2}}\\
&\lesssim&o(t_{\iota n})\lesssim o(t_{2n}-t_{1n}).
\end{eqnarray*}
Summarizing the above, we conclude that
\begin{equation*}
\frac{1}{t_{2n}-t_{1n}}\int_{t_{1n}}^{t_{2n}}\int_{|x|<t}\,|\nabla u|^2-|\partial_tu|^2-|u|^{\dual}\,dxdt\to 0,\quad{\rm as}\,\,n\to\infty.
\end{equation*}
We emphasize the important fact that $t_{2n}-t_{1n}\sim \mu_n$, where $\mu_n \downarrow 0$ is as in (\ref{eq:vanishingcondition6}).\\

Passing to a subsequence and renumbering the indices, we can assume that
\begin{eqnarray}
&&\left|\frac{1}{t_{2n}-t_{1n}}\int_{t_{1n}}^{t_{2n}}\int_{|x|<t}\,|\nabla u|^2-|\partial_tu|^2-|u|^{\dual}\,dxdt\right|\leq 4^{-n};\label{eq:crucialdecay2}\\
&&\frac{1}{\mu_n}\int_{\mu_n}^{2\mu_n}\int_{|x|<4t}\,\left(\partial_tu+\frac{x}{t}\cdot\nabla u+\left(\frac{d}{2}-1\right)\frac{u}{t}\right)^2(x,t)\,dx\,dt\leq 8^{-n}.
\end{eqnarray}
Denoting 
\begin{equation*}
g(t)=\int_{|x|<4t}\,\left(\partial_tu+\frac{x}{t}\cdot\nabla u+\left(\frac{d}{2}-1\right)\frac{u}{t}\right)^2(x,t)\,dx.
\end{equation*}
Then using the boundedness of Hardy-Littlewood maximal functions from $L^1$ to $L^{1,\infty}$, we get that
\begin{equation}\label{eq:smallmeasure1}
\left|\left\{t\in (t_{1n},\,t_{2n}):\,M(g\chi_{(\mu_{n},\,2\mu_{n})})(t)>2^{-n}\right\}\right|\lesssim 4^{-n}|t_{2n}-t_{1n}|.
\end{equation}
Since $\int_{R^d}\,|\nabla_{x,t}u|^2+|u|^{\dual}(x,t)\,dx\leq K$ for some constant $K>0$, we see by  (\ref{eq:crucialdecay2}) that
\begin{eqnarray}
&&\left|\left\{t\in (t_{1n},\,t_{2n}):\,\int_{|x|<t}\,\left[|\nabla u|^2-|\partial_tu|^2-|u|^{\dual}\right](x,t)\,dx<2^{-n}\right\}\right|\nonumber\\
\nonumber\\
&&\gtrsim K^{-1}2^{-n}|t_{2n}-t_{1n}|.\label{eq:smallmeasure2}
\end{eqnarray}

Hence we can find a sequence $t_n\in (t_{1n},\,t_{2n})$, such that along the sequence $t_n$, we have that 
\begin{eqnarray}
&&\lim_{n\to\infty}\,M(g\chi_{(\mu_n,\,2\mu_n)})(t_n)=0,\label{eq:decayno1}\\
&&\limsup_{n\to\infty}\,\int_{|x|<t_n}\,\left[|\nabla u|^2-|\partial_tu|^2-|u|^{\dual}\right](x,t_n)\,dx\leq 0.
\label{eq:decayno2}
\end{eqnarray}

The vanishing condition (\ref{eq:decayno1}) implies by Theorem \ref{th:mainpreliminary}  that there exist $r_0>0$, $J_0\ge 0$, $\lambda^j_{n}\ll t_n$, $c^j_n\in B_{\beta t_n}$ for some $\beta\in (0,1)$ with $\ell_j=\lim\limits_{n\to\infty}\frac{c^j_n}{t_n}$ well defined, for $1\leq j\leq J_0<\infty$, such that in $B_{r_0}$
\begin{eqnarray}
&&\hspace{3cm}\OR{u}(t_n)=\OR{v}(t_n)+\OR{\epsilon}_n+\nonumber\\
&&+\,\sum_{j=1}^{J_0}\left((\lambda^j_{n})^{-\frac{d}{2}+1}Q_{\ell_j}^j\left(\frac{x-c^j_{ n}}{\lambda^j_{n}},\,0\right),\,(\lambda^j_{n})^{-\frac{d}{2}}\partial_tQ_{\ell_j}^j\left(\frac{x-c^j_{ n}}{\lambda^j_{n}},\,0\right)\right),\label{eq:decomposition8}
\end{eqnarray}
with $\OR{\epsilon}_n$ vanishes asymptotically in the sense that 
\begin{equation}\label{goodvanishing13}
\|\epsilon_n\|_{L^{\dual}(B_{r_0})}\to 0,\quad {\rm as}\,\,n\to\infty.
\end{equation}
In addition, the parameters $\lambda_n^j,\,c^j_n$ satisfy the pseudo-orthogonality condition (\ref{eq:pseudo}).\\

We now use (\ref{eq:decayno2}) to obtain refined vanishing properties of $\OR{\epsilon}_n$. We remark that (\ref{eq:decayno2}) is not coercive for general solutions. However it is zero for traveling waves (see Claim \ref{claim:vanishonsoliton}) and coercive for the dispersive part $\OR{\epsilon_n}$ of the solution. Hence it is particularly suited for controlling $\OR{\epsilon}_n$ along the sequence $t_n$. The idea of using an additional virial type quantity such as (\ref{eq:decayno2}) to eliminate dispersive energy was first introduced in \cite{JiaKenig}.  \\

We claim
\begin{claim}\label{claim:vanishonsoliton}
Let $\ell\in R^d$ with $|\ell |<1$, and let $Q_{\ell}$ be a traveling wave with velocity $\ell$. Then $Q_{\ell}(x,t)$ satisfies
\begin{equation}
\int_{R^d}\,\left[|\nabla Q_{\ell}|^2-|\partial_tQ_{\ell}|^2-|Q_{\ell}|^{\dual}\right](x,t)\,dx=0,\quad{\rm for\,\,all}\,\,t.
\end{equation}
\end{claim}
The proof will be given at the end of this section.\\

For the regular part $\OR{v}$, we clearly have
\begin{equation}
\lim_{t\to0+}\,\int_{|x|<t}\left[|\nabla v|^2-|\partial_tv|^2-|v|^{\dual}\right](x,t)\,dx=0.
\end{equation}
By Claim \ref{claim:vanishonsoliton}, and the pseudo-orthogonality of the profiles in the decomposition (\ref{eq:decomposition8}), we can conclude from (\ref{eq:decayno2}) (see the argument after (\ref{eq:vanishingtravelingwave6})) that
\begin{equation}\label{eq:decayepsilon1}
\limsup_{n\to\infty}\,\int_{|x|<t_n}\,|\nabla \epsilon_n|^2-|\partial_t\epsilon_n|^2-|\epsilon_n|^{\dual}\,dx\leq 0.
\end{equation}
By the vanishing condition $\|\epsilon_n\|_{L^{\dual}(|x|<t_n)}\to 0$, we see that
\begin{equation}\label{eq:crucialnegative19}
\limsup_{n\to\infty}\,\int_{|x|<t_n}|\nabla \epsilon_n|^2-|\partial_t\epsilon_n|^2\,dx\leq 0.
\end{equation}
On the other hand, (\ref{eq:decayno1}) implies that
\begin{equation*}
\int_{|x|<t_n}\,\left(\partial_tu+\frac{x}{t_n}\cdot\nabla u+\left(\frac{d}{2}-1\right)\frac{u}{t_n}\right)^2(x,t_n)\,dx\to 0,\,\,{\rm as}\,\,n\to\infty.
\end{equation*}
By the estimates $\int_{|x|<t_n}|u|^2\,dx=o(t_n^2)$ as $n\to\infty$, which follows from similar arguments as in the proof of (\ref{eq:secondvanishing}), we also have
\begin{equation}
\int_{|x|<t_n}\,\left(\partial_tu+\frac{x}{t_n}\cdot\nabla u\right)^2(x,t_n)\,dx\to 0,\,\,{\rm as}\,\,n\to\infty.
\end{equation}
Using the fact that
\begin{equation}\label{eq:vanishingtravelingwave6}
\int_{R^d}\,\left|\partial_tQ_{\ell}+\ell\cdot\nabla Q_{\ell}\right|^2(x,t)\,dx\equiv 0,
\end{equation}
together with $\ell_j=\lim\limits_{n\to\infty}\frac{c^j_n}{t_n}$ and $\lambda^j_n=o(t_n)$ which imply that $$\left((\lambda^j_{n})^{-\frac{d}{2}+1}Q_{\ell_j}^j\left(\frac{x-c^j_{ n}}{\lambda^j_{n}},\,0\right),\,(\lambda^j_{n})^{-\frac{d}{2}}\partial_tQ_{\ell_j}^j\left(\frac{x-c^j_{ n}}{\lambda^j_{n}},\,0\right)\right)$$ is concentrated where $\frac{x}{t_n}=\ell_j+o_n(1)$, we see that
\begin{eqnarray*}
&&\int_{|x|<t_n}\left((\lambda^j_{n})^{-\frac{d}{2}}\partial_tQ_{\ell_j}^j\left(\frac{x-c^j_{ n}}{\lambda^j_{n}},\,0\right)+(\lambda^j_{n})^{-\frac{d}{2}}\frac{x}{t_n}\cdot\nabla Q_{\ell_j}^j\left(\frac{x-c^j_{ n}}{\lambda^j_{n}},\,0\right) \right)^2\,dx\\
&&\\
&&=\int_{|x|<t_n}\left((\lambda^j_{n})^{-\frac{d}{2}}\partial_tQ_{\ell_j}^j\left(\frac{x-c^j_{ n}}{\lambda^j_{n}},\,0\right)+(\lambda^j_{n})^{-\frac{d}{2}}\ell_j\cdot\nabla Q_{\ell_j}^j\left(\frac{x-c^j_{ n}}{\lambda^j_{n}},\,0\right) \right)^2\,dx\\
&&\\
&&\hspace{3.5in}+\,o_n(1)\\
&&\\
&&=\int_{R^d}\left((\lambda^j_{n})^{-\frac{d}{2}}\partial_tQ_{\ell_j}^j\left(\frac{x-c^j_{ n}}{\lambda^j_{n}},\,0\right)+(\lambda^j_{n})^{-\frac{d}{2}}\ell_j\cdot\nabla Q_{\ell_j}^j\left(\frac{x-c^j_{ n}}{\lambda^j_{n}},\,0\right) \right)^2\,dx\\
&&\\
&&\hspace{3.5in}+\,o_n(1)\\
&&\\
&&=\int_{R^d}\,\left|\partial_tQ^j_{\ell_j}+\ell_j\cdot\nabla Q^j_{\ell_j}\right|^2(x,0)\,dx+o_n(1)\\
&&\\
&&=o_n(1).
\end{eqnarray*}
Note that a similar argument was implicitly used in the proof of (\ref{eq:decayepsilon1}).
Hence, from the decomposition (\ref{eq:decomposition8}), we conclude that
\begin{equation}\label{eq:decayepsilon2}
\lim_{n\to\infty}\,\int_{|x|<t_n}\,\left|\partial_t\epsilon_n+\frac{x}{t_n}\cdot\nabla \epsilon_n\right|^2\,dx=0.
\end{equation}
Combining (\ref{eq:decayepsilon2}) and (\ref{eq:crucialnegative19}), we see that
\begin{equation}\label{eq:gooddecayepsilon}
\limsup_{n\to\infty}\,\int_{|x|<t_n}\,|\nabla \epsilon_n|^2-\left|\frac{x}{t_n}\cdot\nabla \epsilon_n\right|^2\,dx\leq 0.
\end{equation}
(\ref{eq:gooddecayepsilon}) implies that for any $\lambda\in(0,1)$
\begin{eqnarray}
&&\int_{|x|<\lambda t_n}|\nabla\epsilon_n|^2\,dx\to 0,\,\,{\rm as}\,\,n\to\infty,\label{eq:epsilondecay2}\\
&&\int_{|x|<t_n}|\spartial \epsilon_n|^2\,dx\to 0,\,\,{\rm as}\,\,n\to\infty.\label{eq:epsilondecay1}
\end{eqnarray}
Combining (\ref{eq:epsilondecay2}) with (\ref{eq:decayepsilon2}), we see that in fact for any $\lambda\in(0,1)$
\begin{equation}
\int_{|x|<\lambda t_n}\,|\nabla_{x,t}\epsilon_n|^2\,dx\to 0,\,\,{\rm as}\,\,n\to\infty.
\end{equation}
The theorem is proved.\\

\smallskip
\noindent
{\it Proof of Claim \ref{claim:vanishonsoliton}.} 
We can assume without loss of generality that $\ell=le_1$, where $e_1=(1,0,\dots,0)\in R^d$. Then 
\begin{equation*}
Q_{\ell}=Q\left(\frac{x_1-lt}{\sqrt{1-l^2}},x_2,\dots,x_d\right).
\end{equation*}
Then direct calculations imply that (see the calculations after (4.35) in \cite{DKM3})
\begin{eqnarray*}
&&\int_{R^d}\left(\,|\nabla Q_{\ell}|^2-|\partial_tQ_{\ell}|^2-|Q_{\ell}|^{\dual}\right)(x,t)\,dx\\
&&=\int_{R^d}\,\left(|\nabla Q|^2-|Q|^{\dual}\right)\left(\frac{x_1-lt}{\sqrt{1-l^2}},x_2,\dots,x_d\right)\,dx\\
&&=\sqrt{1-l^2}\int_{R^d}\,\left(|\nabla Q|^2-|Q|^{\dual}\right)(x)\,dx\\
&&=0.
\end{eqnarray*}
In the last identity we have used the fact that $Q$ is a steady state to equation (\ref{eq:main}). The claim is proved.

\section{Ruling out profiles from time infinity}
In this section, we use the refined vanishing conditions in Theorem \ref{th:betterresidue} to rule out any remaining nontrivial profiles, and obtain better characterizations for the residue term. The arguments are similar to the ones in Section 6 of \cite{DKMsp}.\\

\hspace{.6cm} Let us firstly prove the following result, which describes the profile decomposition of a sequence of special initial data.
\begin{lemma}\label{lm:refinedprofile}
Suppose that $(u_{0n},\,u_{1n})$ is a bounded sequence in $\HL$. In addition, assume that there exists $\lambda_n\uparrow 1$, such that
\begin{equation}\label{refinedvan1}
\|u_{0n}\|_{L^{\dual}(R^d)}+\|u_{1n}\|_{L^2(R^d\backslash B_1)}+\|\spartial\, u_{0n}\|_{L^2(R^d)}\to 0,
\end{equation}
and
\begin{equation}\label{refinedvan2}
\|u_{1n}+\partial_r u_{0n}\|_{L^2(B_1)}+\|\partial_r u_{0n}\|_{L^2\left(B_{|x|<\lambda_n}\bigcup\{x:\,|x|>1\}\right)}\to 0,
\end{equation}
as $n\to\infty$. Then passing to a subsequence, $(u_{0n},\,u_{1n})$ has the following profile decomposition
\begin{eqnarray}
&&\hspace{3cm}(u_{0n},\,u_{1n})=\nonumber\\
&&\sum_{j=1}^J\,\left(\left(\lambda^j_n\right)^{-\frac{d}{2}+1}U_j^L\left(\frac{x-c^j_n}{\lambda^j_n}, \,-\frac{t^j_n}{\lambda^j_n}\right),\,\left(\lambda^j_n\right)^{-\frac{d}{2}}\partial_tU_j^L\left(\frac{x-c^j_n}{\lambda^j_n}, \,-\frac{t^j_n}{\lambda^j_n}\right)\right)\nonumber\\
&&\nonumber\\
&&\quad\quad\quad+\left(w^J_{0n},\,w^J_{1n}\right),\label{eq:refinedprofiledecomposition}
\end{eqnarray}
where the profiles and parameters satisfy, in addition to the usual orthogonality conditions (\ref{eq:profileorthogonality1}) (\ref{eq:profileorthogonality2}), that
\begin{eqnarray}
&&\lim_{n\to\infty}\frac{t^j_n}{\lambda^j_n}\in\{\pm \infty\};\label{eq:refined1}\\
&&\lim_{n\to\infty}|t^j_n|=1;\label{eq:refined2}\\
&&\lim_{n\to\infty}c^j_n=0,\label{eq:refined3}
\end{eqnarray}
for each $j$.
\end{lemma}

\smallskip
\noindent
{\it Proof.} Passing to a subsequence, we can assume that 
$(u_{0n},\,u_{1n})$ has the following profile decomposition
\begin{eqnarray}
&&\hspace{3cm}(u_{0n},\,u_{1n})=\nonumber\\
&&\sum_{j=1}^J\,\left(\left(\lambda^j_n\right)^{-\frac{d}{2}+1}U_j^L\left(\frac{x-c^j_n}{\lambda^j_n}, \,-\frac{t^j_n}{\lambda^j_n}\right),\,\left(\lambda^j_n\right)^{-\frac{d}{2}}\partial_tU_j^L\left(\frac{x-c^j_n}{\lambda^j_n}, \,-\frac{t^j_n}{\lambda^j_n}\right)\right)\nonumber\\
&&\nonumber\\
&&\quad\quad\quad+\left(w^J_{0n},\,w^J_{1n}\right).\label{eq:oldprofiledecomposition}
\end{eqnarray}
Let us firstly show that $\lim\limits_{n\to\infty}\frac{t^j_n}{\lambda^j_n}\in\{\pm\infty\}$. Assume for some $j$, $t^j_n\equiv 0$. By the concentration property of $(u_{0n},\,u_{1n})$, $c^j_n$ is bounded. By the assumption that $\|u_{0n}\|_{L^{\dual}(R^d)}\to 0$ as $n\to\infty$, the profile must satisfy $U^L_j(\cdot,0)\equiv 0$. Hence
\begin{equation}
\int_{R^d}|\partial_tU^L_j+\partial_rU^L_j|^2(x,0)\,dx=\int_{R^d}|\partial_tU^L_j|^2(x,0)\,dx>0.
\end{equation}
By orthogonality of profiles (see Lemma 3.7 in \cite{DKMReview}), we then have
\begin{equation}
\liminf_{n\to\infty}\|u_{1n}+\partial_ru_{0n}\|^2_{L^2(B_1)}\ge \left\|\left(\partial_tU^L_j+\partial_rU^L_j\right)(\cdot,0)\right\|^2_{L^2(R^d)}>0.
\end{equation}
This is a contradiction with (\ref{refinedvan2}). Hence there are no profiles with $t^j_n\equiv 0$.\\

\hspace{.6cm}For the remaining profiles, we divide into several cases, depending on the value of 
$\tau=\lim\limits_{n\to\infty}t^j_n$ (which always exists in $[0,\infty]$ after passing to a subsequence).\\

{\it Case I.} Let us firstly rule out a profile $\OR{U}_j^L$ with $|\tau|\in (0,1)\cup (1,\infty)$. From (\ref{eq:refined1}), we see that $\lambda^j_n\to 0$. 
A moment of reflection, using elementary geometry, shows that the set
\begin{equation}
E_n:=\{x:\,\lambda_n<|x|<1\}
\end{equation}
satisfies for all $M>0$ that
\begin{equation}
\left|\left\{x\in R^d:\,\left||x-c^j_n|-|t^j_n|\right|<M\lambda^j_n\right\}\cap\, E_n\right|=o\left(\left|t^j_n\right|^{d-1}\lambda^j_n\right),
\end{equation}
as $n\to\infty$. Thus we can apply Lemma \ref{lm:removedcone} and conclude that for any $\beta>1$ and sufficiently large $n$,
\begin{equation}
\|(u_{0n},\,u_{1n})\|_{\HL\left(\left\{x:\,\frac{|t^j_n|}{\beta}<|x-c^j_n|<\beta |t^j_n|\right\}\backslash\,E_n\right)}\ge \|\OR{U}^L_j\|_{\HL}>0.
\end{equation}
This is a contradiction with the property of $(u_{0n},\,u_{1n})$ that the energy is concentrated in $E_n$. \\

{\it Case II.} Now we rule out the profile $\OR{U}_j^L$, with $\tau \in\{\pm\infty\}$. This case is similar to {\it Case I.} We omit the routine details.\\

{\it Case III.} Consider the case $\tau=0$. We note that by the concentration property of $(u_{0n},\,u_{1n})$, $\left|c^j_n\right|\to 1$. Passing to a subsequence and rotating the coordinate system if necessary, we can assume that $c^j_n\to e_1:=(1,0\dots,0)\in R^d$. The vanishing conditions (\ref{refinedvan1}) and (\ref{refinedvan2}) of $(u_{0n},\,u_{1n})$ imply that
\begin{equation}\label{goodv12}
\int_{B_{2t^j_n}(c^j_n)}\left|u_{1n}+\partial_1u_{0n}\right|^2+|\nabla_{x'}u_{0n}|^2\,dx\to 0,
\end{equation}
as $n\to\infty$, where $x'=(x_2,\ldots,x_d)$.
By Lemma \ref{lm:localizedorthogonality}, we have that for any $\beta>1$
\begin{eqnarray}
&&\liminf_{n\to\infty}\int_{\frac{t^j_n}{\beta}<|x-c^j_n|<\beta t^j_n}\,(u_{1n}+\partial_1u_{0n})^2+|\nabla_{x'}u_{0n}|^2\,dx\nonumber\\
&&\ge\int_{R^d}\,\left[\left(\partial_tU^L_j+\partial_1U^L_j\right)^2+|\nabla_{x'}U^L_j|^2\right](x,0)\,dx>0.\label{eq:hoholocalizedexpansion}
\end{eqnarray}
(\ref{eq:hoholocalizedexpansion}) contradicts (\ref{goodv12}).\\

{\it Case IV.} We have $|\tau|=1$. Using the same arguments as in {\it Case I}, we also have $\lim\limits_{n\to\infty}c^j_n=0$. These are the profiles that appear in the decomposition (\ref{eq:refinedprofiledecomposition}). \\

The lemma is proved.\\

Now we are ready to prove the following theorem ruling out profiles from time infinity.
\begin{theorem}\label{th:decompbetter}
Let $\OR{u}\in C((0,1],\HL(R^d))$, with $u\in \Snorm(R^d\times (\epsilon,1])$ for any $\epsilon>0$, be a type II blow-up solution to equation (\ref{eq:main}) with blow-up time $t=0$. Assume that $0\in\mathcal{S}$ is a singular point. Then there exist an integer $J_0\ge 1$, $r_0>0$, a time sequence $t_n\downarrow 0$,  scales $\lambda_n^j$ with $0<\lambda_n^j\ll t_n$, positions $c_n^j\in R^d$ satisfying $c_n^j\in B_{\beta t_n}$ for some $\beta\in(0,1)$ with $\ell_j=\lim\limits_{n\to\infty}\frac{c_n^j}{t_n}$ well defined, and traveling waves $Q_{\ell_j}^j$, for $1\leq j\leq J_0$, such that inside the ball $B_{r_0}$ we have
\begin{eqnarray}
&&\hspace{3cm}\OR{u}(t_n)=\OR{v}(t_n)+(\epsilon_{0n},\,\epsilon_{1n})+\nonumber\\
&&+\sum_{j=1}^{J_0}\,\left((\lambda_n^j)^{-\frac{d}{2}+1}\, Q_{\ell_j}^j\left(\frac{x-c_n^j}{\lambda_n^j},\,0\right),\,(\lambda_n^j)^{-\frac{d}{2}}\, \partial_tQ_{\ell_j}^j\left(\frac{x-c_n^j}{\lambda_n^j},\,0\right)\right).\label{eq:maindecomposition}
\end{eqnarray}
Moreover, let $\OR{\epsilon}^L_n$ be the solution to linear wave equation with initial data $(\epsilon_{0n},\,\epsilon_{1n})\in\HL$, then the following vanishing condition holds:
\begin{eqnarray}
&&\|\epsilon^L_n\|_{\Snorm(R^d\times R)}+\|( \epsilon_{0n},\,\epsilon_{1n})\|_{\HL(\{x:\,|x|<\lambda t_n,\,{\rm or}\,|x|>t_n\})}+\nonumber\\
\nonumber\\
&&\quad\quad\,\,+\,\left\|\epsilon_{1n}+\frac{r}{t_n}\partial_r\epsilon_{0n}\right\|_{L^2(|x|\leq t_n)}+\|\spartial\,\epsilon_{0n}\|_{L^2(R^d)}\to 0,\label{eq:refinedintro}
\end{eqnarray}
as $n\to\infty$ for any $\lambda\in(0,1)$. In the above $\spartial$ denotes the tangential derivative.
In addition, the parameters $\lambda_n^j,\,c_n^j$ satisfy the pseudo-orthogonality condition (\ref{eq:profileorthogonality2})
for each $1\leq j\neq j'\leq J_0$.
\end{theorem}

\smallskip
\noindent
{\it Proof.} By Theorem \ref{th:betterresidue}, we only need to show that $J_0\ge 1$ and 
\begin{equation}\label{eq:lastclaim}
\lim_{n\to\infty}\|\epsilon^L_n\|_{\Snorm(R^d\times R)}=0.
\end{equation}
By Lemma \ref{lm:refinedprofile} with a rescaling, passing to a subsequence, $(\epsilon_{0n},\,\epsilon_{1n})$ has the following profile decomposition
\begin{eqnarray}
&&\hspace{3cm}(\epsilon_{0n},\,\epsilon_{1n})=\left(w^J_{0n},\,w^J_{1n}\right)+\label{eq:lastrefinedprofiledecomposition}\\
&&\sum_{j=1}^J\,\left(\left(\lambda^j_n\right)^{-\frac{d}{2}+1}U_j^L\left(\frac{x-c^j_n}{\lambda^j_n}, \,-\frac{t^j_n}{\lambda^j_n}\right),\,\left(\lambda^j_n\right)^{-\frac{d}{2}}\partial_tU_j^L\left(\frac{x-c^j_n}{\lambda^j_n}, -\frac{t^j_n}{\lambda^j_n}\right)\right),\nonumber
\end{eqnarray}
where the profiles and parameters satisfy, in addition to the usual conditions (\ref{eq:profileorthogonality1}), (\ref{eq:profileorthogonality2}) on the parameters and vanishing condition for $w^J_n$, that
\begin{eqnarray}
&&\lim_{n\to\infty}\frac{t^j_n}{\lambda^j_n}\in\{\pm \infty\};\label{eq:lastrefined1}\\
&&\lim_{n\to\infty}\frac{|t^j_n|}{t_n}=1;\label{eq:lastrefined2}\\
&&\lim_{n\to\infty}\frac{c^j_n}{t_n}=0,\label{eq:lastrefined3}
\end{eqnarray}
for each $j$. \\

We shall firstly rule out the profiles with $\lim\limits_{n\to\infty}\frac{t^j_n}{t_n}=1$. Consider the following sequence of initial data
\begin{eqnarray}
&&\hspace{3cm}(\widetilde{u}_{0n},\,\widetilde{u}_{1n})=\nonumber\\
&&\sum_{k=1}^J\,\left(\left(\lambda^k_n\right)^{-\frac{d}{2}+1}U_k^L\left(\frac{x-c^k_n}{\lambda^k_n}, \,-\frac{t^k_n}{\lambda^k_n}\right),\,\left(\lambda^k_n\right)^{-\frac{d}{2}}\partial_tU_k^L\left(\frac{x-c^k_n}{\lambda^k_n}, \,-\frac{t^k_n}{\lambda^k_n}\right)\right)\nonumber\\
&&\quad\quad\quad+\OR{v}(t_n)+(w^J_{0n},\,w^J_{1n}).\label{eq:modifiedinitialdata}
\end{eqnarray}
Hence we have removed any possible solitons from $\OR{u}(t_n)$, in order to be able to control the evolution of solution $\OR{\widetilde{u}}$ to equation (\ref{eq:main}) with $\OR{\widetilde{u}}(t_n)=(\widetilde{u}_{0n},\,\widetilde{u}_{1n})$. Note that
\begin{equation}\label{lastdifference}
\lim_{n\to\infty}\,\|(u_{0n},\,u_{1n})-(\widetilde{u}_{0n},\,\widetilde{u}_{1n})\|_{\HL(|x|>\frac{1+\beta}{2} t_n)}=0.
\end{equation}
Lemma \ref{lm:nonlinearprofiledecomposition} and the definition of nonlinear profiles imply that for $t\in [\frac{t_n}{2},\,t_n]$, $\OR{\widetilde{u}}$ has the following asymptotic expansion
\begin{eqnarray}
\OR{\widetilde{u}}(t)&=&\OR{v}(t)+\OR{w}^J_n(t)+\OR{r}^J_n(t)+\nonumber\\
&&\quad\quad+\sum_{k=1}^J\,\bigg(\left(\lambda^k_n\right)^{-\frac{d}{2}+1}U_k^L\left(\frac{x-c^k_n}{\lambda^k_n}, \,\frac{t-t_n-t^k_n}{\lambda^k_n}\right),\nonumber\\
&&\hspace{4cm}\,\left(\lambda^k_n\right)^{-\frac{d}{2}}\partial_tU_k^L\left(\frac{x-c^k_n}{\lambda^k_n}, \,\frac{t-t_n-t^k_n}{\lambda^k_n}\right)\bigg),\nonumber
\end{eqnarray}
where we recall that $\OR{w}^J_n(t)$ is the solution to the linear wave equation with $\OR{w}^J_n(t_n)=(w^J_{0n},\,w^J_{1n})$, and
\begin{equation*}
\lim_{J\to\infty}\,\limsup_{n\to\infty}\,\sup_{t\in[t_n/2,t_n]}\|\OR{r}^J_n\|_{\HL}=0.
\end{equation*}
At $t=\frac{t_n}{2}$, if $\lim\limits_{n\to\infty}\frac{t^j_n}{t_n}=1$, then by Lemma \ref{lm:concentrationoffreeradiation} and (\ref{eq:lastrefined1},\,\ref{eq:lastrefined2},\,\ref{eq:lastrefined3}), the energy of
$$\left(\left(\lambda^j_n\right)^{-\frac{d}{2}+1}U_j^L\left(\frac{x-c^j_n}{\lambda^j_n}, \,\frac{t-t_n-t^j_n}{\lambda^j_n}\right),\,\left(\lambda^j_n\right)^{-\frac{d}{2}}\partial_tU_j^L\left(\frac{x-c^j_n}{\lambda^j_n}, \,\frac{t-t_n-t^j_n}{\lambda^j_n}\right)\right)$$
is concentrated in $\{x:\,\frac{3t_n}{2\gamma}<|x|<\frac{3\gamma t_n}{2}\}$, for any $\gamma>1$. Hence by Lemma \ref{lm:localizedorthogonality}, $\OR{\widetilde{u}}(\frac{t_n}{2})$ has a uniform amount of nontrivial energy accumulation in $\{x:\,\frac{3t_n}{2\gamma}<|x|<\frac{3\gamma t_n}{2}\}$. Note also by (\ref{lastdifference}) and Lemma \ref{lm:localinspaceapproximation}, we have 
\begin{equation}
\left\|\OR{u}(\frac{t_n}{2})-\OR{\widetilde{u}}(\frac{t_n}{2})\right\|_{\HL(\{x:\,\frac{3t_n}{2\gamma}<|x|<\frac{3\gamma t_n}{2}\})}\to 0,
\end{equation}
as $n\to\infty$ if $\gamma$ is chose close to $1$. Thus $\OR{u}(\frac{t_n}{2})$ also has a uniform amount of nontrivial energy accumulation in  $\{x:\,\frac{3t_n}{2\gamma}<|x|<\frac{3\gamma t_n}{2}\}$, as $n\to\infty$. This is a contradiction with the fact that outside the lightcone $|x|>t$, $\OR{u}$ is regular. \\

Next we rule out profiles with $\lim\limits_{n\to\infty}\frac{t^j_n}{t_n}=-1$, keeping in mind that the profiles with $\lim\limits_{n\to\infty}\frac{t^j_n}{t_n}=1$ have been ruled out. Introduce $\OR{\widetilde{u}}$ as in the first part of the proof. We fix a small time $\delta>0$. By Lemma \ref{lm:nonlinearprofiledecomposition} and the definition of nonlinear profiles, $\OR{\widetilde{u}}$ has the following asymptotic expansion for $t\in [t_n,\delta]$:
\begin{eqnarray}
\OR{\widetilde{u}}(t)&=&\OR{v}(t)+\OR{w}^J_n(t)+\OR{r}^J_n(t)\nonumber\\
&&\quad\quad+\sum_{k=1}^J\,\bigg(\left(\lambda^k_n\right)^{-\frac{d}{2}+1}U_k^L\left(\frac{x-c^k_n}{\lambda^k_n}\,\frac{t-t_n-t^k_n}{\lambda^k_n}\right),\nonumber\\
&&\hspace{3cm}\,\left(\lambda^k_n\right)^{-\frac{d}{2}}\partial_tU_k^L\left(\frac{x-c^k_n}{\lambda^k_n}, \,\frac{t-t_n-t^k_n}{\lambda^k_n}\right)\bigg),\label{lastapproximation}
\end{eqnarray}
where $\OR{w}^J_n(t_n)=(w^J_{0n},\,w^J_{1n})$, and
\begin{equation*}
\lim_{J\to\infty}\,\limsup_{n\to\infty}\,\sup_{t\in[t_n,\delta]}\|\OR{r}^J_n\|_{\HL}=0.
\end{equation*}
Fix $\beta<\beta_1<\beta_2<1$. Define the modified sequence of initial data 
\begin{equation}
(u_{0n}',\,u_{1n}'):=\left\{\begin{array}{ll}
                                        \OR{u}(t_n)&\,\,{\rm if}\,\,|x|>\beta_2t_n,\\
   &\\
\frac{t_n^{-1}|x|-\beta_1}{\beta_2-\beta_1}\OR{u}(t_n)+\frac{\beta_2-t_n^{-1}|x|}{\beta_2-\beta_1}\OR{\widetilde{u}}(t_n)&\,\,{\rm if}\,\,\beta_1t_n<|x|<\beta_2t_n,\\
&\\
                       \OR{\widetilde{u}}(t_n) &  \,\,{\rm if}\,\,|x|<\beta_1t_n.
                               \end{array}\right.
\end{equation}
It is easy to verify that for $r_0>0$ sufficiently small,
$$\lim_{n\to\infty}\|(u_{0n}',\,u_{1n}')-\OR{\widetilde{u}}(t_n)\|_{\HL(B_{r_0})}=0.$$
Let $\OR{u'}$ be the solution to equation (\ref{eq:main}) with $\OR{u'}(t_n)=(u_{0n}',\,u_{1n}')$. By finite speed of propagation, we see that
\begin{equation}\label{eq:identical}
u(x,t)=u'(x,t),\,\,{\rm for}\,\,|x|>\beta_2\,t_n+t-t_n,\,t\in [t_n,\delta].
\end{equation}
By Lemma \ref{lm:longtimeperturbation},  and by choosing $\delta$ smaller if necessary, we see that 
\begin{equation}\label{thesame}
\sup_{t\in[t_n,\,\delta]}\,\|\OR{u'}(t)-\OR{\widetilde{u}}(t)\|_{\HL(B_{r_0-\delta})}\to 0,\,\,{\rm as}\,\,n\to\infty.
\end{equation}
By (\ref{lastapproximation}), Lemma \ref{lm:concentrationoffreeradiation}, and the pseudo-orthogonality of profiles, we get that at time $t=\delta$, $$\liminf_{n\to\infty}\left\|\OR{\widetilde{u}}(\delta)\right\|_{\HL(\delta+\beta_2t_n-t_n<|x|<\delta+t_n)}\ge\|\OR{U}^L_j(0)\|_{\HL(R^d)}.$$
Hence by (\ref{thesame}) and (\ref{eq:identical}), $\OR{u}(\delta)$ has a uniform amount of nontrivial energy concentration in the region $\delta+\beta_2t_n-t_n<|x|<\delta+t_n$, as $n\to\infty$. This is a contradiction with the fact that $\OR{u}(\delta)\in\HL$. 

\hspace{.6cm}In summary, there are no nontrivial profiles in the decomposition (\ref{eq:lastrefinedprofiledecomposition}). Thus, $\epsilon^L$ verifies (\ref{eq:lastclaim}). Then $J_0\ge 1$ follows since $0\in \mathcal{S}$. Theorem \ref{th:decompbetter} is proved.\\

\section{A linear channel of energy inequality and conclusion of the proof}
It remains to show that the dispersive term actually vanishes asymptotically in the energy space, not only in the sense of (\ref{eq:refinedintro}). The key ingredient is the following ``channel of energy inequality" for linear wave equations. \\

\begin{lemma}\label{lm:linearchannel}
Suppose that $(u_0,u_1)\in\HL(R^d)$ satisfies ${\rm supp}\,(u_0,u_1)\subset \overline{B_R(0)}$, for some $R>0$. Let $u$ be the solution to the linear wave equation $$\partial_{tt}u-\Delta u=0,\,\,{\rm in}\,\,R^d\times R,$$
with initial data $\OR{u}(0)=(u_0,u_1)$. Denote
\begin{eqnarray}
E_0&=&\int_{R^d}\left[\frac{|\nabla u|^2}{2}+\frac{|\partial_tu|^2}{2}\right](x,t)\,dx;\label{eq:energy1}\\
E_{-1}&=&\int_{R^d}\left[\frac{| u|^2}{2}+\frac{\left|\partial_t|\nabla|^{-1}u\right|^2}{2}\right](x,t)\,dx.\label{eq:energy0}
\end{eqnarray}
$E_0$ and $E_{-1}$ are conserved quantities.
Fix any $\eta_0\in(0,R)$. Then we have for all $t\ge 0$,
\begin{eqnarray}
&&(R-\eta_0)\int_{|x|\ge t+\eta_0}\left[\frac{|\nabla u|^2}{2}+\frac{|\partial_tu|^2}{2}\right](x,t)\,dx\nonumber\\
&&\quad\ge -\int_{R^d}x\cdot\nabla u_0\,u_1\,dx-E_0\eta_0-C(d)E_0^{\frac{1}{2}}E_{-1}^{\frac{1}{2}}.\label{eq:linearchannel}
\end{eqnarray}
\end{lemma}

\smallskip
\noindent
{\it Proof.} The fact that $E_0$ and $E_{-1}$ are conserved quantities follows from standard energy conservation for linear wave equations. Note that $\OR{u}(t)\in \dot{H}^1\times L^2$ is compactly supported for all $t$ and $d\ge 3$, so that $E_{-1}$ is well-defined. 

Let us prove (\ref{eq:linearchannel}).  Direct calculation shows that
\begin{equation*}
\frac{d}{dt}\int_{R^d}\left[u_t\,\left(x\cdot\nabla u+\frac{d-1}{2}u\right)\right](x,t)\,dx=-E_0.
\end{equation*}
Thus,
\begin{equation}\label{eq:virial113}
\int_{R^d}\left[u_t\,\left(x\cdot\nabla u+\frac{d-1}{2}u\right)\right](x,t)\,dx=\int_{R^d}u_1\,\left(x\cdot\nabla u_0+\frac{d-1}{2}u_0\right)\,dx-tE_0.
\end{equation}
Note that by Cauchy-Schwarz inequality and the fact that $E_0$ and $E_{-1}$ are conserved,  
\begin{equation}\label{colf}
\int_{R^d}|\partial_tu|\,|u|\,dx\leq C(d)E_0^{\frac{1}{2}}E_{-1}^{\frac{1}{2}}.
\end{equation}
Combining (\ref{eq:virial113}) and (\ref{colf}), we obtain that
\begin{equation}
\int_{R^d}\left[x\cdot\nabla u\,\partial_tu\right](x,t)\,dx\leq \int_{R^d}x\cdot\nabla u_0\,u_1\,dx-tE_0+C(d)E_0^{\frac{1}{2}}E_{-1}^{\frac{1}{2}}.
\end{equation}
By splitting the integral into the regions $|x|>t+\eta_0$ and $|x|<t+\eta_0$, we get that
\begin{eqnarray}
\int_{|x|>t+\eta_0}\left[x\cdot\nabla u\,\partial_tu\right](x,t)\,dx&\leq& \int_{R^d}x\cdot\nabla u_0\,u_1\,dx-tE_0\nonumber\\
      &&\quad+(t+\eta_0)\int_{|x|<t+\eta_0}\frac{|\nabla_{x,t}u|^2}{2}\,dx+C(d)E_0^{\frac{1}{2}}E_{-1}^{\frac{1}{2}}\nonumber\\
\nonumber\\
&\leq&-t\int_{|x|>t+\eta_0}\frac{|\nabla_{x,t}u|^2}{2}\,dx+E_0\eta_0\nonumber\\
&&\quad\quad-\,\eta_0\int_{|x|>t+\eta_0}\frac{|\nabla_{x,t}u|^2}{2}\,dx+\nonumber\\
\nonumber\\
&&\quad\quad\,\,+\,\int_{R^d}x\cdot\nabla u_0\,u_1\,dx+C(d)E_0^{\frac{1}{2}}E_{-1}^{\frac{1}{2}}.\label{eq:upperbound8}
\end{eqnarray}
By the finite speed of propagation, ${\rm supp}\,(u(t),\,\partial_tu(t))\subset \overline{B_{t+R}(0)}$. Hence,
\begin{equation}\label{eq:alowerbound8}
\int_{|x|>t+\eta_0}\left[x\cdot\nabla u\,\partial_tu\right](x,t)\,dx\ge -(t+R)\int_{|x|>t+\eta_0}\frac{|\nabla_{x,t}u|^2}{2}.
\end{equation}
From (\ref{eq:upperbound8}) and (\ref{eq:alowerbound8}), we get that
\begin{equation*}
-(R-\eta_0)\int_{|x|>t+\eta_0}\frac{|\nabla_{x,t}u|^2}{2}(x,t)\,dx\leq E_0\eta_0+\int_{R^d}x\cdot\nabla u_0\,u_1\,dx+C(d)E_0^{\frac{1}{2}}E_{-1}^{\frac{1}{2}}.
\end{equation*}
This is exactly (\ref{eq:linearchannel}). The lemma is proved.\\

In general the right hand side of (\ref{eq:linearchannel}) may not be positive and the estimate may not be useful. But if the initial data is in some sense ``well prepared", such as in the case of the dispersive error verifying the vanishing condition (\ref{eq:refinedintro}), (\ref{eq:linearchannel}) becomes very powerful.
More precisely, we have the following corollary.

\begin{corollary}\label{cor:channel}
Suppose that $(\epsilon_{0n},\epsilon_{1n})$ is a bounded sequence in $\HL(R^d)$, which vanishes asymptotically in the sense that for any $\lambda\in(0,1)$,
\begin{eqnarray}
&&\|( \epsilon_{0n},\,\epsilon_{1n})\|_{\HL(\{x:\,|x|<\lambda,\,\,{\rm or}\,\,|x|>1\})}+\nonumber\\
\nonumber\\
&&\quad\quad\,\,+\,\left\|\epsilon_{1n}+\partial_r\epsilon_{0n}\right\|_{L^2(R^d)}+\|\spartial\,\epsilon_{0n}\|_{L^2(R^d)}\to 0,\label{eq:vanish8}
\end{eqnarray}
as $n\to\infty$. Let $\epsilon^L_n$ be the solution of the linear wave equation with initial data $(\epsilon_{0n},\epsilon_{1n})$. Assume that 
$$\inf_{n}\|(\epsilon_{0n},\,\epsilon_{1n})\|_{\HL}>0.$$
Then for any $\eta_0\in(0,1)$ and sufficiently large $n$ depending on $\eta_0$ and the inf above, we have the following channel of energy inequality
\begin{equation}\label{eq:truechannel}
\int_{|x|>t+\eta_0}\left[\frac{|\nabla \epsilon^L_n|^2}{2}+\frac{|\partial_t\epsilon^L_n|^2}{2}\right](x,t)\,dx\ge \frac{\eta_0}{2}\|(\epsilon_{0n},\epsilon_{1n})\|^2_{\HL},
\end{equation}
for all $t\ge0$.
\end{corollary}

\smallskip
\noindent
{\it Proof.} Let us first localize $(\epsilon_{0n},\epsilon_{1n})$, in order to apply Lemma \ref{lm:linearchannel}. Fix $R>1$ to be determined below. Let $\eta\in C_c^{\infty}(B_R)$ with $\eta|_{B_1}\equiv 1$. Denote 
$$(\widetilde{\epsilon}_{0n},\widetilde{\epsilon}_{1n}):=\eta \,(\epsilon_{0n},\epsilon_{1n}).$$
By the vanishing condition (\ref{eq:vanish8}), we have that
\begin{equation}\label{same88}
\left\|(\widetilde{\epsilon}_{0n},\widetilde{\epsilon}_{1n})-(\epsilon_{0n},\epsilon_{1n})\right\|_{\HL(R^d)}\to 0,\,\,\,\,{\rm as}\,\,n\to\infty.
\end{equation}
By the vanishing condition (\ref{eq:vanish8}) and (\ref{same88}), it is straightforward to verify that
\begin{equation}\label{eq:nonde8}
\int_{R^d}x\cdot\nabla \widetilde{\epsilon}_{0n}\,\widetilde{\epsilon}_{1n}= -\int_{R^d}\left[\frac{|\nabla \widetilde{\epsilon}_{0n}|^2}{2}+\frac{|\widetilde{\epsilon}_{1n}|^2}{2}\right]\,dx+o_n(1),\,\,{\rm as}\,\,n\to\infty.
\end{equation}
In addition, (\ref{eq:vanish8}) implies that $(\widetilde{\epsilon}_{0n},\widetilde{\epsilon}_{1n})\rightharpoonup 0$ as $n\to\infty$. Hence by the regularity and support property of $(\widetilde{\epsilon}_{0n},\widetilde{\epsilon}_{1n})$, we conclude that
\begin{equation}\label{eq:de8}
\|(\widetilde{\epsilon}_{0n},\widetilde{\epsilon}_{1n})\|_{L^2\times\dot{H}^{-1}(R^d)}\to 0,\,\,{\rm as}\,\,n\to\infty.
\end{equation}
Let $\widetilde{\epsilon}^L_n$ be the free radiation with initial data $(\widetilde{\epsilon}_{0n},\widetilde{\epsilon}_{1n})$. Applying Lemma \ref{lm:linearchannel} we get that for sufficiently large $n$ and all $t\ge 0$,
\begin{eqnarray*}
&&(R-\eta_0)\int_{|x|\ge t+\eta_0}\left[\frac{|\nabla \widetilde{\epsilon}^L_n|^2}{2}+\frac{|\widetilde{\epsilon}^L_n|^2}{2}\right](x,t)\,dx\nonumber\\
&&\quad\ge -\int_{R^d}x\cdot\nabla \widetilde{\epsilon}_{0n}\,\widetilde{\epsilon}_{1n}\,dx-E_{0n}\eta_0-C(d)E_{0n}^{\frac{1}{2}}E_{-1n}^{\frac{1}{2}}\\
&&\quad\ge (1-\eta_0)E_{0n}+o_n(1),\,\,{\rm as}\,\,n\to\infty.
\end{eqnarray*}
In the above, $E_{0n}$ and $E_{-1n}$ are defined as in (\ref{eq:energy1}) and (\ref{eq:energy0}) for the sequence $(\widetilde{\epsilon}_{0n},\widetilde{\epsilon}_{1n})$. Take $R$ sufficiently close to $1$. By (\ref{same88}), for sufficiently large $n$, we conclude that
\begin{equation}\label{prechannel}
\int_{|x|>t+\eta_0}\left[\frac{|\nabla \widetilde{\epsilon}^L_n|^2}{2}+\frac{|\partial_t\widetilde{\epsilon}^L_n|^2}{2}\right](x,t)\,dx\ge \frac{\eta_0+1}{4}\|(\epsilon_{0n},\epsilon_{1n})\|^2_{\HL},
\end{equation}
for all $t\ge0$. (\ref{same88}) implies
$$\sup_{t\ge 0}\left\|(\widetilde{\epsilon}^L_n,\,\partial_t\widetilde{\epsilon}^L_n)-(\epsilon^L_n,\,\partial_t\epsilon^L_n)\right\|_{\HL}\to 0,$$
as $n\to\infty$.
The corollary follows from (\ref{prechannel}) and the above bound.\\

Now we are ready to prove the main theorem in the finite time blow up case.

\begin{theorem}\label{th:mainFinal}
 Let $\OR{u}$ satisfy the assumptions of Theorem \ref{th:decompbetter}. Then the conclusion of Theorem \ref{th:decompbetter} holds with the following additional property:
\begin{equation}
 \label{final_decay}
\left\|(\epsilon_{0n},\epsilon_{1n})\right\|^2_{\HL}\to 0,\,\,{\rm as}\,\,n\to\infty.
 \end{equation} 
 \end{theorem}
%

\smallskip
\noindent
{\it Proof.} By Theorem \ref{th:decompbetter}, along a sequence of times $t_n\to 0+$, we already have the decomposition (\ref{eq:maindecomposition}), with the residue term $(\epsilon_{0n},\epsilon_{1n})$ vanishing in the sense that 
\begin{eqnarray}
&&\|\epsilon^L_n\|_{\Snorm(R^d\times R)}+\|( \epsilon_{0n},\,\epsilon_{1n})\|_{\HL(\{x:\,|x|<\lambda t_n,\,{\rm or}\,|x|>t_n\})}+\nonumber\\
\nonumber\\
&&\quad\quad\,\,+\,\left\|\epsilon_{1n}+\frac{r}{t_n}\partial_r\epsilon_{0n}\right\|_{L^2(B_{t_n})}+\|\spartial\,\epsilon_{0n}\|_{L^2(R^d)}\to 0,\label{eq:refinedintroFinal}
\end{eqnarray}
as $n\to\infty$ for any $\lambda\in(0,1)$. We need to show that (\ref{final_decay}) holds.
%

Assuming that along a subsequence we have $$\left\|(\epsilon_{0n},\epsilon_{1n})\right\|^2_{\HL}\ge \mu_0>0,$$
we will arrive at a contradiction. 

By a rescaled and time-translated version of Corollary \ref{cor:channel}, for any $\eta_0\in(0,1)$ and sufficiently large $n$, the linear solution $\epsilon^L_n$ verifies the following channel of energy inequality 
\begin{equation}\label{eq:channelforepsilon}
\int_{|x|>t+\eta_0t_n}\frac{\left|\nabla_{x,t}\epsilon^L_n\right|^2}{2}(x,t)\,dx\ge \frac{\eta_0}{2}\,\left\|(\epsilon_{0n},\epsilon_{1n})\right\|^2_{\HL}\ge \frac{\eta_0\mu_0}{4},\,\,{\rm for\,\,all}\,\,t\ge 0.
\end{equation}

Choose $\eta_0>\beta>0$, with $\beta$ being given by Theorem \ref{th:decompbetter}.
Define a new sequence of initial data $(w_{0n},w_{1n})\in\HL$ as follows
\begin{equation}
(w_{0n},w_{1n}):=\OR{v}(t_n)+(\epsilon_{0n},\epsilon_{1n}).
\end{equation}
By the decomposition (\ref{eq:maindecomposition}) and the choice $\eta_0>\beta$, we see that
\begin{equation}\label{eq:exteriorapp}
\|(w_{0n},w_{1n})-\OR{u}(t_n)\|_{\HL(\{|x|\ge \eta_0t_n,\,\,{\rm and}\,\,|x|<r_0\})}\to 0,\,\,{\rm as}\,\,n\to\infty.
\end{equation}
Let $w_n$ be the solution to equation (\ref{eq:main}) with $\OR{w}_n(t_n)=(w_{0n},w_{1n})$. We shall use $r_0,\,\delta>0$ from Theorem \ref{th:decompbetter} and its proof. Shrinking $\delta$ if necessary, we assume $4\delta<r_0$. Then by the vanishing condition (\ref{eq:refinedintroFinal}) and a simple case of Lemma \ref{lm:nonlinearprofiledecomposition}, we can conclude that
\begin{equation}\label{eq:goodapp18}
\OR{w}_n(t)=\OR{v}(t)+\left(\epsilon^L_n,\partial_t\epsilon^L_n\right)(t-t_n)+\OR{r}_n(t),
\end{equation}
where 
\begin{equation}
\sup_{t\in(t_n,\delta)}\|\OR{r}_n(t)\|_{\HL}\to 0,\,\,{\rm as}\,\,n\to\infty.
\end{equation}
By (\ref{eq:exteriorapp}) and finite speed of propagation, we get that
\begin{equation}\label{eq:almostid18}
\sup_{t\in(t_n,\delta)}\|\OR{u}(t)-\OR{w}_n(t)\|_{\HL(\delta\ge|x|\ge t-t_n+\eta_0t_n)}\to 0,
\end{equation}
as $n\to\infty$.

(\ref{eq:channelforepsilon}) and (\ref{eq:goodapp18}), together with the concentration property of $(\epsilon_{0n},\epsilon_{1n})$ and finite speed of propagation, imply that for sufficiently large $n$ and all $t\in(t_n,\delta)$,
\begin{equation}
\int_{t-t_n+\eta_0t_n\leq |x|\leq t}\frac{\left|\nabla_{x,t}w_n\right|^2}{2}(x,t)\,dx\ge \frac{\mu_0\eta_0}{6}>0.
\end{equation}
Hence by (\ref{eq:almostid18}), we also have for sufficiently large $n$,
\begin{equation}\label{eq:channelforu}
\int_{t-t_n+\eta_0t_n\leq |x|\leq t}\frac{\left|\nabla_{x,t}u\right|^2}{2}(x,t)\,dx\ge \frac{\mu_0\eta_0}{8}>0.
\end{equation}
Applying (\ref{eq:channelforu}) at $t=\frac{\delta}{2}$, we thus conclude that
$$\int_{\frac{\delta}{2}-t_n+\eta_0t_n\leq |x|\leq \frac{\delta}{2}}\frac{\left|\nabla_{x,t}u\right|^2}{2}\left(x,\delta/2\right)\,dx\ge \frac{\mu_0\eta_0}{8}>0.$$
Since $t_n\to 0+$ as $n\to\infty$, we have a contradiction with the finiteness of $\left\|\OR{u}(\frac{\delta}{2})\right\|_{\HL}$, which precludes concentration of energy in arbitrarily small regions. The theorem is proved.

\section{Soliton resolution along a sequence of times for global, bounded solutions}
In this section, we briefly outline the main steps in the proof of Theorem \ref{th:main} in the global existence case. The proof follows similar arguments as in the finite time blow up case. We will therefore only highlight the new features.\\

Let us first recall that for a global bounded solution $\OR{u}$ to equation (\ref{eq:main}), we can extract the {\it radiation term} from $u$, see \cite{DKMsp}.
\begin{lemma}\label{lm:radiation}
Let $\OR{u}$ be a solution to equation (\ref{eq:main}) on $R^d\times [0,\infty)$ such that 
\begin{equation}
\label{bounded_global}
\sup_{t\ge 0}\,\|\OR{u}(t)\|_{\HL}<\infty.
\end{equation}
Then there exists a unique finite energy solution $\OR{u}^L$ to the linear wave equation $$\partial_{tt}u^L-\Delta u^L=0,$$
such that
\begin{equation}
\lim_{t\to\infty}\int_{|x|\ge t-A} \left[|\nabla (u-u^L)|^2+|\partial_t(u-u^L)|^2\right](x,t)\,dx=0,
\end{equation}
for any $A>0$.
\end{lemma}

Let us consider a solution $\OR{u}$ to equation (\ref{eq:main}) on $R^d\times [0,\infty)$ such that (\ref{bounded_global}) holds and let $u^L$ be its radiation term. By Lemma \ref{lm:radiation}, for a sufficiently large $T>0$, we have
\begin{eqnarray}
&&\|u^L\|_{\Snorm(R^d\times (T,\infty))}+\noindent\\
\nonumber\\
&&\quad\quad+\,\sup_{t\ge T}\int_{|x|\ge t-1}\left[|\nabla (u-u^L)|^2+|\partial_t(u-u^L)|^2\right](x,t)\,dx\leq \epsilon_0.\label{eq:small11}
\end{eqnarray}
In the above, $\epsilon_0=\epsilon_0(d)$ is a small positive number, chosen so that we can apply small data global existence and scattering result for equation (\ref{eq:main}). If $\epsilon_0$ is chosen sufficiently small, then by finite speed of propagation and small data theory for equation (\ref{eq:main}), we can conclude that
\begin{equation}
\|u\|_{\Snorm\left(\{(x,t):\,|x|\ge t,\,t\ge T\}\right)}\leq C(d)\epsilon_0.
\end{equation}
Since $u\in \Snorm(R^d\times (0,T))$, it follows that 
\begin{equation}\label{eq:finiteexterior}
u\in \Snorm(\{(x,t):\,|x|\ge t>0\}).
\end{equation}
Then by the assumption (\ref{bounded_global}), the bound (\ref{eq:finiteexterior}) and H\"older inequality, we get that $\partial_{x,t}|u|^{\dual}\in L^1_{x,t}\left(\{(x,t):\,|x|\ge t>0\}\right)$. Thus, by simple real analysis arguments, we conclude that $|u|^{\dual}$ can be restricted to the lightcone $\{(x,t):\,|x|=t\}$. That is, 
\begin{equation}\label{eq:controlnegative}
\int_0^{\infty}\int_{|x|=t}|u|^{\dual}\,d\sigma dt<\infty.
\end{equation}
As in the case of finite blow up case, the bound (\ref{eq:controlnegative}) implies that the main part 
\begin{equation}\label{eq:mainpartofenergyflux11}
\int_0^{\infty}\int_{|x|=t}\,\left[\frac{|\nabla u|^2}{2}+\frac{|\partial_tu|^2}{2}-\frac{x}{|x|}\cdot\nabla u\,\partial_tu\right](x,t)\,d\sigma dt
\end{equation}
of the energy flux
\begin{equation*}
\int_0^{\infty}\int_{|x|=t}\,\left[\frac{|\nabla u|^2}{2}+\frac{|\partial_tu|^2}{2}-\frac{x}{|x|}\cdot\nabla u\,\partial_tu-\frac{|u|^{\dual}}{\dual}\right](x,t)\,d\sigma dt
\end{equation*}
can be controlled. The control of the main flux term (\ref{eq:mainpartofenergyflux11}) implies the following Morawetz estimate for $u$:\\

for any $0<10\,t_1<t_2<\infty$, we have
\begin{equation}\label{eq:morawetz11}
\int_{t_1}^{t_2}\int_{|x|<t}\left(\partial_tu+\frac{x}{t}\cdot\nabla u+\left(\frac{d}{2}-1\right)\frac{u}{t}\right)^2\,dx\,\frac{dt}{t}\leq C\left(\log{\frac{t_2}{t_1}}\right)^{\frac{1}{2}},
\end{equation}
for some $C$ independent of $t_1,\,t_2$.\\

Again, the estimate (\ref{eq:morawetz11}) implies that the integral 
$$\int_{|x|<t}\left(\partial_tu+\frac{x}{t}\cdot\nabla u+\left(\frac{d}{2}-1\right)\frac{u}{t}\right)^2\,dx$$
vanishes on average asymptotically as $t\to\infty$. By the same arguments as in the finite time blow up case, and the property for linear waves that for all $C>0$ (see \cite{DKMsp}, Theorem 2.1)
$$\lim_{t\to\infty}\int_{|x|<Ct}\left(\partial_tu^L+\frac{x}{t}\cdot\nabla u^L+\left(\frac{d}{2}-1\right)\frac{u^L}{t}\right)^2\,dx=0,$$
we can show that there exist two sequences of times $t_{1n}$ and $t_{2n}$ approaching $\infty$ with $t_{1n}\sim t_{2n}\sim t_{2n}-t_{1n}$, such that for all $C>0$
\begin{eqnarray}
&&\sup_{0<\tau<\frac{t_{\iota n}}{16}}\frac{1}{\tau}\int_{|t_{\iota n}-t|<\tau}\,\int_{|x|<Ct}\,\left(\partial_tu+\frac{x}{t}\cdot\nabla u+\left(\frac{d}{2}-1\right)\frac{u}{t}\right)^2\,dx\,dt\nonumber\\
\nonumber\\
&&\quad\quad\quad\,\,\,\,\to 0,\,\,{\rm as}\,\,n\to \infty.\label{eq:vanishingM13}
\end{eqnarray}
for $\iota=1,\,2$.\\

Using exactly the same arguments as in the finite time blow up case, we can show that along each time sequence $t_{\iota n}$, we have the following preliminary decomposition: \footnote{More precisely, these decompositions depend on $\iota=1,2$, but we suppress the dependence for ease of notation.}\\

There exist scales $\lambda_n^j$ with $0<\lambda_n^j\ll t_n$, positions $c_n^j\in R^d$ satisfying $c_n^j\in B_{\beta\, t_n}$ for some $\beta\in(0,1)$, with $\ell_j=\lim\limits_{n\to\infty}\frac{c_n^j}{t_n}$ well defined, and traveling waves $Q_{\ell_j}^j$, for $1\leq j\leq J_0$, such that
\begin{eqnarray}
\OR{u}(t_n)&=&\sum_{j=1}^{J_0}\,\left((\lambda_n^j)^{-\frac{d}{2}+1}\, Q_{\ell_j}^j\left(\frac{x-c_n^j}{\lambda_n^j},\,0\right),\,(\lambda_n^j)^{-\frac{d}{2}}\, \partial_tQ_{\ell_j}^j\left(\frac{x-c_n^j}{\lambda_n^j},\,0\right)\right)+\nonumber\\
\nonumber\\
&&\quad\quad\quad\,\,+\,\OR{u}^L(t_n)+(\epsilon_{0n},\epsilon_{1n}),\label{eq:maindecomposition16}
\end{eqnarray}
where the residue term $(\epsilon_{0n},\epsilon_{1n})$ vanishes asymptotically in the following sense:
\begin{equation}\label{eq:vanL6}
\|\epsilon_{0n}\|_{L^{\dual}(R^d)}\to 0,\quad {\rm as}\,\,n\to\infty;
\end{equation}
if we write as
\begin{equation}
(\epsilon_{0n},\,\epsilon_{1n})(x)=\left((\lambda^j_n)^{-\frac{d}{2}+1}\widetilde{\epsilon}^{\,j}_{0n}\left(\frac{x-c^j_n}{\lambda^j_n}\right),\,(\lambda^j_n)^{-\frac{d}{2}}\widetilde{\epsilon}^{\,j}_{1n}\left(\frac{x-c^j_n}{\lambda^j_n}\right)\right),
\end{equation}
then $(\widetilde{\epsilon}^j_{0n},\widetilde{\epsilon}^j_{1n})\rightharpoonup 0$ as $n\to\infty$, for each $j\leq J_0.$
In addition, the parameters $\lambda_n^j,\,c^j_n$ satisfy the pseudo-orthogonality condition (\ref{eq:profileorthogonality2}), for $1\leq j\neq j'\leq J_0.$\\

Indeed, along a subsequence of $t_{\iota n}$, we have the profile decomposition
\begin{eqnarray}
&&\hspace{4cm}\OR{u}(t_{\iota n})=\OR{u}^L(t_{\iota n})+\nonumber\\
&&+\,\sum_{j=1}^J\,\left(\left(\lambda^j_n\right)^{-\frac{d}{2}+1}U_j^L\left(\frac{x-c^j_n}{\lambda^j_n}, \,-\frac{t^j_n}{\lambda^j_n}\right),\,\left(\lambda^j_n\right)^{-\frac{d}{2}}\partial_tU_j^L\left(\frac{x-c^j_n}{\lambda^j_n}, \,-\frac{t^j_n}{\lambda^j_n}\right)\right)\nonumber\\
&&\nonumber\\
&&\quad\quad\quad\hspace{3cm}+\left(w^J_{0n},\,w^J_{1n}\right),\label{eq:oldprofiledecompositionglobalcase}
\end{eqnarray}
where the parameters satisfy the usual orthogonality conditions.
By Lemma \ref{lm:radiation} and the discussion on pages 144-145 of \cite{BaGe}, we have $\lambda^j_n\lesssim t_n$, and $\lim\limits_{n\to\infty}\frac{\lambda^j_n}{t_n}\neq 0$ for at most one $j$, say $j=j_0$, $\left|c^j_n\right|\lesssim t_n$, $t^j_n\lesssim t_n$. By passing to a subsequence, if such $j_0$ exists, we can choose $\lambda^{j_0}_n=t_n$, and let $\ell_j=\lim\limits_{n\to\infty}\frac{c^j_n}{t_n}$ be well-defined for each $j$. \\

We can again divide the profiles into three cases, as in the finite time blow up case, and claim\\

\begin{itemize}

\item {\it Case I:} \,\,$t^{j_0}_n\equiv 0$, $\lambda^{j_0}_n\equiv t_n$. In this case, the nonlinear profile $U_{j_0}$ is a compactly supported self similar solution, a case ruled out by \cite{KMacta};\\

\item {\it Case II:} \,\, $t^j_n\equiv 0$, $\lambda^j_n\ll t_n$. In this case, $|\ell_j|<1$ and $\OR{U}_j^L(\cdot,0)=\OR{Q}^j_{\ell_j}$;\\

\item {\it Case III:}\,\, $\lambda^j_n\ll \left|t^j_n\right|$. In this case, one can absorb these profiles into the error term. \\

\end{itemize}
The proofs of the above claim follow exactly the same lines as in the finite blow up case, using (\ref{eq:vanishingM13}), perturbative techniques to control profiles with small energy or with asymptotically vanishing interaction with the profile under consideration, and removing  profiles with large energy and nontrivial interaction with the profile under consideration, to obtain the first order equation (\ref{eq:firstorderselfsimilar}) for $U_{j_0}$ in Case I and (\ref{eq:firstorder222}) for $U^L_j$ in Case II. The first order equations, together with equation (\ref{eq:main}), show that $U_{j_0}$ must be trivial and that $U^L_j$ must be a traveling wave, exactly as in the finite blow up case. We omit the details.\\

The decomposition (\ref{eq:maindecomposition16}) and the vanishing condition (\ref{eq:vanL6}) imply that along the time sequences $t_{\iota n}$, we have
\begin{equation}\label{eq:secondvanishing16}
\|u(\cdot,t_n)\|_{L^2(|x|<t_n)}=o(t_n),\,\,{\rm as}\,\,n\to\infty.
\end{equation}
Thus we can use a second virial identity as in the finite time blow up case, and obtain 
\begin{equation*}
\frac{1}{t_{2n}-t_{1n}}\int_{t_{1n}}^{t_{2n}}\int_{|x|<t}\,|\nabla u|^2-|\partial_tu|^2-|u|^{\dual}\,dxdt\to 0,\quad{\rm as}\,\,n\to\infty.
\end{equation*}

Then along a new sequence of times $t_n\to\infty$ as $n\to\infty$, we can achieve two vanishing conditions simultaneously as in Section \ref{sec:virial}, which allows us to obtain a preliminary decomposition of $\OR{u}(t_n)$, but with a better residue term:
\begin{eqnarray}
\OR{u}(t_n)&=&\sum_{j=1}^{J_0}\,\left((\lambda_n^j)^{-\frac{d}{2}+1}\, Q_{\ell_j}^j\left(\frac{x-c_n^j}{\lambda_n^j},\,0\right),\,(\lambda_n^j)^{-\frac{d}{2}}\, \partial_tQ_{\ell_j}^j\left(\frac{x-c_n^j}{\lambda_n^j},\,0\right)\right)+\nonumber\\
\nonumber\\
&&\quad\quad\quad\,\,+\,\OR{u}^L(t_n)+(\epsilon_{0n},\epsilon_{1n}).\label{eq:maindecomposition19}
\end{eqnarray}
Let $\OR{\epsilon}^L_n$ be the solution to the linear wave equation with initial data $(\epsilon_{0n},\,\epsilon_{1n})\in\HL$, then the following vanishing condition holds:
\begin{eqnarray}
&&\|\epsilon^L_n\|_{\Snorm(R^d\times R)}+\|( \epsilon_{0n},\,\epsilon_{1n})\|_{\HL(\{x:\,|x|<\lambda t_n,\,{\rm or}\,|x|>t_n\})}+\nonumber\\
\nonumber\\
&&\quad\quad\,\,+\,\left\|\epsilon_{1n}+\frac{r}{t_n}\partial_r\epsilon_{0n}\right\|_{L^2(|x|\leq t_n)}+\|\spartial\,\epsilon_{0n}\|_{L^2(R^d)}\to 0,\label{eq:refined19}
\end{eqnarray}
as $n\to\infty$ for any $\lambda\in(0,1)$.\\

The final task is to show that actually $\|(\epsilon_{0n},\epsilon_{1n})\|_{\HL(R^d)}\to 0$, as $n\to\infty$. For this purpose, we appeal to the channel of energy inequality from Lemma \ref{lm:linearchannel} and Corollary \ref{cor:channel}. The proof in this case is slightly different from the finite time blow up case, and we thus sketch some of the details below.\\

Assume that for a subsequence of $t_n$, we have that
$$\|(\epsilon_{0n},\epsilon_{1n})\|_{\HL(R^d)}^2\ge \mu_0>0.$$
We will arrive at a contradiction. By a rescaled and time-translated version of Corollary \ref{cor:channel}, the linear solution $\epsilon^L_n$ verifies the following channel of energy inequality, for any $\eta_0\in(0,1)$ and sufficiently large $n$,
\begin{equation}\label{eq:channelforepsilon2}
\int_{|x|>t+\eta_0t_n}\frac{\left|\nabla_{x,t}\epsilon^L_n\right|^2}{2}(x,t)\,dx\ge \frac{\eta_0}{2}\,\left\|(\epsilon_{0n},\epsilon_{1n})\right\|^2_{\HL}\ge \frac{\eta_0\mu_0}{4},\,\,{\rm for\,\,all\,\,}t\ge 0.
\end{equation}
Choose $\eta_0>\beta>0$, where $\beta$ is as in the paragraph before (\ref{eq:maindecomposition16}).
Define a new sequence of initial data $(w_{0n},w_{1n})\in\HL$ as follows
\begin{equation}
(w_{0n},w_{1n}):=\OR{u}^L(t_n)+(\epsilon_{0n},\epsilon_{1n}).
\end{equation}
By the decomposition (\ref{eq:maindecomposition19}) and the choice that $\eta_0>\beta$, we see that
\begin{equation}\label{eq:exteriorapp2}
\|(w_{0n},w_{1n})-\OR{u}(t_n)\|_{\HL(|x|\ge \eta_0t_n)}\to 0,\,\,{\rm as}\,\,n\to\infty.
\end{equation}
Let $w_n$ be the solution to equation (\ref{eq:main}) with $\OR{w}_n(t_n)=(w_{0n},w_{1n})$. Then by the vanishing condition (\ref{eq:refined19}) and a simple case of Lemma \ref{lm:nonlinearprofiledecomposition}, we can conclude that
\begin{equation}\label{eq:goodapp19}
\OR{w}_n(t)=\OR{u}^L(t)+\left(\epsilon^L_n,\partial_t\epsilon^L_n\right)(t-t_n)+\OR{r}_n(t),
\end{equation}
where 
\begin{equation}\label{eq:rvanishes12}
\sup_{t\ge t_n}\|\OR{r}_n(t)\|_{\HL}\to 0,\,\,{\rm as}\,\,n\to\infty.
\end{equation}
By (\ref{eq:exteriorapp2}) and finite speed of propagation, we get that
\begin{equation}\label{eq:almostid19}
\sup_{t\ge t_n}\|\OR{u}(t)-\OR{w}_n(t)\|_{\HL(|x|\ge t-t_n+\eta_0t_n)}\to 0,
\end{equation}
as $n\to\infty$.

(\ref{eq:channelforepsilon2}), (\ref{eq:goodapp19}) and (\ref{eq:rvanishes12}), together with the concentration properties of $(\epsilon_{0n},\epsilon_{1n})$ and finite speed of propagation, imply that for sufficiently large $n$ and all $t\ge t_n$,
\begin{equation}
\int_{t-t_n+\eta_0t_n\leq |x|\leq t}\frac{\left|\nabla_{x,t}(w_n-u^L)\right|^2}{2}(x,t)\,dx\ge \frac{\eta_0\mu_0}{6}>0.
\end{equation}
Hence by (\ref{eq:almostid19}), for sufficiently large $n$ and all $t\ge t_n$ we also have
\begin{equation}\label{eq:channelforu2}
\int_{t-t_n+\eta_0t_n\leq |x|\leq t}\frac{\left|\nabla_{x,t}(u-u^L)\right|^2}{2}(x,t)\,dx\ge \frac{\eta_0\mu_0}{8}>0.
\end{equation}
Fix now $n_0$ so that (\ref{eq:channelforu2}) holds, and let $$A=-t_{n_0}+\eta_0t_{n_0}.$$
Then, the lower bound (\ref{eq:channelforu2}) contradicts Lemma \ref{lm:radiation}. Hence 
$$\|(\epsilon_{0n},\epsilon_{1n})\|_{\HL(R^d)}\to 0, \,\,{\rm as}\,\, n\to\infty,$$
and this concludes the proof of Theorem \ref{th:main} in the global existence case.

\appendix

\section{Appendix: Residue term for the profile decomposition}
\label{A:residue}
It was proved in \cite{Bulut} that in the profile decomposition for the wave equation in $R^6\times R$ in the energy space, the residue term $w^J_n$ vanishes asymptotically in the sense that
$$\lim_{J\to\infty}\limsup_{n\to\infty}\|w^J_n\|_{L^{\frac{7}{2}}(R^6\times R)}=0.$$
Here we consider the two end point cases, not considered in \cite{Bulut}, that is, whether $w^J_n$ vanishes asymptotically in the Strichartz norms $L^{\infty}_tL^3(R^6\times R)$  and $L^2_tL^4_x(R^6\times R)$. 

 We first show that the residue term $w^J_n(x,t)$ vanishes asymptotically in the sense that 
\begin{equation}\label{eq:vanishing1}
\lim_{J\to\infty}\limsup_{n\to\infty}\left\|w^J_n\right\|_{L^{\infty}_tL^3_x(R^6\times R)}=0.
\end{equation}
The key ingredient is the following lemma.

\begin{lemma}[Inverse Sobolev inequality]\label{lm:concentration}
Suppose that $f\in\dot{H}^1(R^6)$ with $$\|f\|_{\dot{H}^1(R^6)}\leq M<\infty$$ and $$\|f\|_{L^3(R^6)}\ge \delta>0,$$ then there exists $\epsilon=\epsilon(M,\delta)>0$, $x_0\in R^6$ and $r_0>0$, such that
\begin{equation}\label{eq:concentration}
\|f\|_{L^2(B_{r_0}(x_0))}\ge \epsilon \,r_0\,.\\
\end{equation}
\end{lemma}
(See \cite{Gerard} or Chapter 4 of \cite{Killip}.)

 Now we give a sketch of proof for (\ref{eq:vanishing1}). Assume that (\ref{eq:vanishing1}) does not hold. Then there exist a sequence $J_k\to \infty$ and $\delta>0$ such that 
$$\limsup_{n\to\infty}\left\|w^{J_k}_n\right\|_{L^{\infty}_tL^3_x(R^6\times R)}\ge \delta.$$
Using 
$$
\lim_{k\to\infty}\limsup_{n\to\infty}\left\|w^{J_k}_n\right\|_{L^{7/2}(R^6\times R)}=0,
$$
we obtain a sequence $n_k\to\infty$ such that
\begin{equation}
\label{limsup1}
\limsup_{k\to\infty}\left\|w^{J_k}_{n_k}\right\|_{L^{\infty}_tL^3_x(R^6\times R)}\ge \delta/2
\end{equation} 
and
\begin{equation}
\label{limsup2}
 \lim_{k\to\infty}\left\|w^{J_k}_{n_k}\right\|_{L^{7/2}(R^6\times R)}=0.
\end{equation} 
By (\ref{limsup1}), we obtain that at some $t_k$, $\|w^{J_k}_{n_k}(\cdot,t_k)\|_{L^3(R^6)}\gtrsim \delta$. Since by energy conservation and Pythagorean expansion, $\|w^{J_k}_{n_k}(\cdot,t_k)\|_{\dot{H}^1(R^6)}\leq M$, we can apply Lemma 10.1 and conclude that there exists $r_k>0$, $x_k\in R^6$ and $\epsilon=\epsilon(M,\delta)>0$, such that
\begin{equation}\label{eq:concentrationforn}
\left\|w^{J_k}_{n_k}(\cdot,t_k)\right\|_{L^2(B_{r_k}(x_k))}\ge \epsilon\,r_k.
\end{equation}
By translation and rescaling, we see that 
\begin{equation}\label{eq:nontrivial111}
\left\|r_k^2w^{J_k}_{n_k}\big(x_k+r_k\cdot,\,t_k\big)\right\|_{L^2(B_1)}\ge \epsilon.
\end{equation}
Extracting a subsequence, we obtain that for some $(U_0,U_1)\in \dot{H}^1\times L^2$,
$$ \left(r_k^2w_{n_k}^{J_k}\big(x_k+r_k\cdot,\,t_k\big),r_k\partial_t w_{n_k}^{J_k}\big(x_k+r_k\cdot,\,t_k\big)\right)\rightharpoonup (U_0,U_1),$$
 as $k\to\infty$, weakly in $\dot{H}^1\times L^2(R^6)$. By \eqref{eq:nontrivial111} and the compactness of the embedding $H^1_{{\rm loc}}(R^6)\subset L^2_{{\rm loc}}(R^6)$, we obtain that $U_0\neq 0$. By Strichartz estimates, denoting by $U^L$ the solution of the linear wave equation with data $(U_0,U_1)$,
$$r_k^2w_{n_k}^{J_k}\big(x_k+r_k\cdot,\,t_k+r_k\cdot\big)\rightharpoonup U^L,$$
 as $k\to\infty$, weakly in $L^{7/2}(R^6\times R)$.
By a classical property of the weak convergence,
\begin{multline*}
\liminf_{k\to\infty}\left\| w_{n_k}^{J_k}\right\|_{L^{7/2}(R^6\times R)}=\liminf_{k\to\infty}\left\| r_k^2w_{n_k}^{J_k}\big(x_k+r_k\cdot,t_k+r_k\cdot\big)\right\|_{L^{7/2}(R^6\times R)}\\
\ge \left\|U^L\right\|_{L^{7/2}(R^6\times R)}>0, 
\end{multline*}
contradicting (\ref{limsup2}).\\

We remark that the other endpoint vanishing condition 
\begin{equation}\label{eq:vanishing2}
\lim_{J\to\infty}\limsup_{n\to\infty}\left\|w^J_n\right\|_{L^{2}_tL^4_x(R^6\times R)}=0
\end{equation}
{\it does not} hold in general. A simple counter example is the following. Take a nontrivial smooth finite energy solution $U^L(x,t)$ to the linear wave equation. For simplicity, we also assume that 
$${\rm supp}\,(U^L(\cdot,0),\partial_tU^L(\cdot,0))\subset B_1.$$
Then 
\begin{equation}\label{eq:suppPRO}
{\rm supp}\, U^L\subset \{(x,t)\in R^6\times R:\,|x|<1+|t|\},
\end{equation}
 and by Theorem 1.1 in \cite{weightedestimate}, 
\begin{equation}\label{eq:decayC1}
|\partial^{\alpha} U^L(x,t)|\lesssim \frac{1}{(1+|t|)^{\frac{5}{2}}(1+||t|-|x||)^{\frac{5}{2}}}\chi_{|x|<1+|t|},
\end{equation}
${\rm for}\,\,t\in R\,\,\,{\rm and\,\,all\,\,}\alpha\in\big(\mathbb{N}\cup \{0\}\big)^6.$
 Fix $N>1$ and let $t^N_n>0$ for $1\leq n\leq N$ be such that for $n\neq n'$
$$\left|t^N_n-t^N_{n'}\right|\ge \epsilon_{N}^{-1}>8$$
for a tiny number $\epsilon_{N}$ to be chosen below. Consider now the initial data 
$$(u_{0N},u_{1N})=\sum_{n=1}^N\,\frac{1}{\sqrt{N}}\OR{U}^L(x,-t^N_n).$$
If we choose $\epsilon_{N}$ sufficiently small (in this case the profiles are almost orthogonal, because $t^N_n$ are far apart from each other), then by the support property (\ref{eq:suppPRO}) and decay property (\ref{eq:decayC1}), we get that
$$\left\|(u_{0N},u_{1N})\right\|_{\HL}\approx 1.$$
The solution to the linear wave equation with initial data $(u_{0N},u_{1N})$ is
$$u_N^L(x,t)=\sum_{n=1}^N\,\frac{1}{\sqrt{N}}U^L(x,t-t^N_n).$$
Assume again that the $\epsilon_{N}$ is small enough, then by (\ref{eq:suppPRO}) and (\ref{eq:decayC1}) we have
\begin{eqnarray*} 
&&\left\|u_N^L\right\|_{L^2_tL^4_x(R^6\times R)}^2\\
&&\gtrsim\sum_{n=1}^N\int_{|t-t_n^N|\leq \frac{1}{3\epsilon_N}}\left\|\frac{1}{\sqrt{N}}U^L(\cdot,t-t_n^N)\right\|_{L^4(R^6)}^2\,dt\\
&&\approx \sum_{n=1}^N \frac{1}{N}\|U^L\|^2_{L^2_tL^4_x(R^6\times R)} \approx 1,
\end{eqnarray*}
but 
$$\left\|u_N^L\right\|_{L^{\infty}_tL^3_x(R^6\times R)}\approx \frac{1}{\sqrt{N}}\|U^L\|_{L^{\infty}_tL^3_x(R^6\times R)},$$
and similarly
$$\left\|u_N^L\right\|_{L^{\frac{7}{2}}(R^6\times R)}^{\frac{7}{2}}\approx \sum_{n=1}^N\left(\frac{1}{\sqrt{N}}\right)^{\frac{7}{2}}\|U^L\|_{L^{\frac{7}{2}}(R^6\times R)}^{\frac{7}{2}}\approx \frac{1}{N^{\frac{3}{4}}}\|U^L\|_{L^{\frac{7}{2}}(R^6\times R)}^{\frac{7}{2}},$$
if $\epsilon_{N}$ is chosen sufficiently small, depending on $N$. Fix such choices of $\epsilon_{N},\,t^N_n$ (again, depending on $N$) for each $N>1$.  For the bounded sequence $(u_{0N},u_{1N})\in \HL$, with associated linear evolution $u^L_N$, we have
$$\|(u_{0N},u_{1N})\|_{\HL}\approx 1,\,\,\,\left\|u^L_N\right\|_{L^2_tL^4_x(R^6\times R)}\approx 1,\,\,\left\|u^L_N\right\|_{L^{\infty}_tL^3_x(R^6\times R)}=o_N(1).$$
 Hence, the sequence $(u^N_0,u^N_1)$ provides a counter example to (\ref{eq:vanishing2}).\\

\smallskip
\noindent
{\it Remark:} The analogous results hold for all $d\ge 3$. Moreover, the corresponding results are also valid for the Schr\"{o}dinger equation.

\vspace{1cm}

\center{{\bf Acknowledgement}\\
We thank Guixiang Xu for valuable comments on a preliminary version of the paper.}

\end{document}